\documentclass[12pt,english]{amsart}
\usepackage[T1]{fontenc}
\usepackage[latin9]{inputenc}
\usepackage[letterpaper]{geometry}
\usepackage{amsthm}
\usepackage{amsmath}
\usepackage{amssymb}

\makeatletter
\numberwithin{equation}{section}
\numberwithin{figure}{section}
\theoremstyle{plain}
\newtheorem{thm}{Theorem}[section]
  \theoremstyle{definition}
  \newtheorem{defn}[thm]{Definition}
  \theoremstyle{remark}
  \newtheorem{rem}[thm]{Remark}
  \theoremstyle{plain}
  \newtheorem{lem}[thm]{Lemma}
  \theoremstyle{plain}
  \newtheorem{prop}[thm]{Proposition}
  \theoremstyle{plain}
  \newtheorem{cor}[thm]{Corollary}
 \theoremstyle{definition}
  \newtheorem{example}[thm]{Example}

\usepackage{color}

\usepackage{babel}

\makeatother

\usepackage{babel}

\begin{document}

\title[The John-Nirenberg and John-Str\"omberg theorems for BMO]{A new look at the John-Nirenberg and \\ John-Strömberg theorems
for BMO
\vskip 5mm  
Lecture Notes
}

\author{Michael Cwikel, Yoram Sagher and Pavel Shvartsman}

\address{Cwikel and Shvartsman: Department of Mathematics, Technion - Israel
Institute of Technology, Haifa 32000, Israel }

\address{Sagher: Department of Mathematics, Florida Atlantic University, 777
Glades Road, Boca Raton, FL 33431, USA }

\email{mcwikel@math.technion.ac.il }

\email{yoram.sagher@gmail.com}

\email{pshv@tx.technion.ac.il}

\thanks{The research of the first and third named authors was 
supported by funding from the Martin and Sima Jelin Chair in 
Mathematics, by the Technion V.P.R.\ Fund and by the Fund for Promotion 
of Research at the Technion.}

\begin{abstract}
We develop some techniques for studying various versions of the function
space $BMO$. Special cases of one of our results give alternative
proofs of the celebrated John-Nirenberg inequality and of related
inequalities due to John and to Wik. Our approach enables us to pose
a simply formulated {}``geometric'' question, for which an affirmative
answer would lead to a version of the John-Nirenberg inequality with
dimension free constants.\\
A more detailed summary of the main ideas and results of this paper
can be found at \texttt{http://www.math.technion.ac.il/\textasciitilde{}mcwikel/bmo/CwikSaghShvaSummary.pdf}
\end{abstract}
\maketitle
\section{\label{sec:intro}Introduction. Our main question.}
We begin by inviting the reader to consider and hopefully even answer
the following question. We will subsequently refer to it as {}``Question
A''.\bigskip{}

\textit{Do there exist two absolute constants $\tau\in(0,1/2)$ and
$s>0$ which have the following property?}

\textbf{\textit{For every positive integer $d$ and for every cube
$Q$ in $\mathbb{R}^{d}$, whenever $E_{+}$ and $E_{-}$ are two
disjoint measurable subsets of $Q$ whose $d$-dimensional Lebesgue
measures satisfy \[
\min\left\{ \lambda(E_{+}),\lambda(E_{-})\right\} >\tau\lambda(Q\setminus E_{+}\setminus E_{-})\,,\]
 then there exists some cube $W$ contained in $Q$ for which\[
\min\left\{ \lambda(W\cap E_{+}),\lambda(W\cap E_{-})\right\} \ge s\lambda(W)\,\]
}} 

\bigskip{}

We are led to consider this question because of our interest in the
space $BMO$ of functions of bounded mean oscillation introduced by
John and Nirenberg \cite{JohnNirenberg}. We recall that these are
the functions $f:\mathbb{R}^{d}\to\mathbb{R}$ which have the property
that \[
\sup_{Q}\frac{1}{\lambda(Q)}\int_{Q}\left|f-f_{Q}\right|d\lambda<\infty\]
where the supremum is taken over all cubes $Q$ in $\mathbb{R}^{d}$
and where $f_{Q}$ is the average of $f$ on $Q$. 

We will show that an affirmative answer to Question A would have very
interesting consequences for the study of a remarkable property of
functions of bounded mean oscillation. It would imply (see Theorem
\ref{thm:mainjs}) that the following {}``dimension free'' version
\begin{equation}
\lambda\left(\left\{ x\in Q:\left|f(x)-m_{f}\right|\ge\alpha\right\} \right)\le\max\left\{ \frac{1}{2\tau},2\sqrt{\frac{1}{2\tau}}\right\} \cdot\lambda(Q)\cdot\exp\left(-\frac{\alpha s\log\frac{1}{2\tau}}{8\left\Vert f\right\Vert _{BMO}}\right)\label{eq:newurp}\end{equation}
of the John-Nirenberg inequality \cite{JohnNirenberg} holds for every
$\alpha\ge0$. It would also imply some slightly stronger inequalities.
(Here $m_{f}$ is any median of the measurable function $f$ on the
cube $Q$ in $\mathbb{R}^{d}$.)

Having formulated our question, let us now state to what extent we
have been able, so far, to answer it or to simplify it.

For each particular value of $d\in\mathbb{N}$ we \textit{can} find
numbers, $\tau\in(0,1/2)$ and $s>0$ which \textit{do} have the property
sought in Question A. Furthermore we can show that their having this
property, implies that the inequality (\ref{eq:newurp}) is satisfied.

We do not yet have an answer to Question A, because at least one of
our constants $\tau$ and $s$ depends on $d$. We can take, for example,
$\tau=\sqrt{2}-1$, but, for that choice of $\tau$, we have only
been able to obtain a value of $s$ which depends on $d$, namely
$s=2^{-d}\left(3-2\sqrt{2}\right)$.

Regardless of whether $\tau$ and $s$ really have to depend on the
dimension, it seems of interest that, in the expression of the form
$C\lambda(Q)\cdot\exp\left({\displaystyle -\frac{c\alpha}{\left\Vert f\right\Vert _{BMO}}}\right)$
on the right hand side of our version (\ref{eq:newurp}) of the John-Nirenberg
inequality, we have revealed a quite explicit connection between the
constants $C$ and $c$ and a geometric property expressed by the
constants $\tau$ and $s$.

It also seems of interest that the {}``geometric'' condition sought
in Question A is, more or less {}``equivalent'' to an analytic condition
which compares certain kinds of $BMO$ {}``norms'' of functions
$f$ on $\mathbb{R}^{d}$ with related kinds of $BMO$ {}``norms''
of their rearrangements $f^{*}$ on $(0,\infty)$. (The implications,
in two opposite directions, which express this {}``sort of equivalence''
are precisely formulated and established in Theorems \ref{thm:Nultra}
and \ref{thm:Ninverse}.)

As we shall show in Section \ref{sec:attempts}, if Question A can
be answered affirmatively in some special cases, then this will suffice
to answer it in general. In particular it would suffice to give an
affirmative answer in the case where the subsets $E_{+}$ and $E_{-}$
are each finite unions of dyadic subcubes of $Q$. Thus Question A
can be considered to be a combinatoric question as much as a geometric
one.

Our results can be expressed in more abstract terms, and they apply
to other versions of the space $BMO$ including the one considered
by Wik \cite{wik}, where cubes are replaced by {}``false cubes''.

This preliminary version of our paper is written
more or less in the style of {}``lecture notes''. We hope that this
will make it helpful for graduate students and that experts will forgive
us for writing more, maybe much more than they need to read about
various things. We have attempted to find some sort of reasonable
middle way between what may suit these two subsets of our audience
by relegating quite a number of better known facts, results and proofs
to the appendices in Section \ref{sec:Appendices}.

We have surely omitted references to some very pertinent papers about
this topic, but we hope to correct at least some of our omissions
in future versions. The reader is invited to draw our attention to
any such omissions.

\section{\label{sec:Notation-and-terminology}Notation, terminology and some
more introduction}

Throughout this paper $d$ will denote a positive integer and $\lambda$
will denote $d$-dimensional Lebesgue measure on $\mathbb{R}^{d}$.
The value of $d$ will always be clear from the context. When $d=1$
we will also often use the notation $\left|E\right|$ instead of $\lambda(E)$
for each measurable subset $E$ of $\mathbb{R}$. By a \textit{cube}
in $\mathbb{R}^{d}$ we will always mean a $d$-dimensional closed
cube with sides parallel to the axes.
\begin{defn}
\label{def:admissibleset}To save tedious repetitions of terminology,
we will say that a set $E$ is \textit{admissible} if it is a measurable
subset (i.e., a Lebesgue measurable subset) of $\mathbb{R}^{d}$ and
its $d$-dimensional Lebesgue measure satisfies $0<\lambda(E)<\infty$. 
\end{defn}
For each admissible set $E$ and each measurable real valued function
$f$ whose domain of definition contains $E$\,, we define the \textit{mean
oscillation} of $f$ on $E$ by \begin{equation}
\mathbf{O}(f,E):=\inf_{c\in\mathbb{R}}\frac{1}{\lambda(E)}\int_{E}\left|f-c\right|d\lambda\,.\label{eq:defv}\end{equation}

It is convenient to fix some notation for two other frequently used
variants of the functional $\mathbf{O}(f,E)$. So we set \[
\mathbf{A}(f,E):=\frac{1}{\lambda(E)}\int_{E}\left|f-f_{E}\right|d\lambda\]
for every function $f$ which is integrable on $E$, and where $f_{E}:=\frac{1}{\lambda(E)}\int_{E}fd\lambda$.
We also set \[
\mathbf{D}(f,E):=\frac{1}{\lambda(E)^{2}}\iint_{E\times E}\left|f(x)-f(y)\right|d\lambda(x)d\lambda(y)\]
({}``\textbf{A}'' and {}``\textbf{D}'' are our abbreviations for
{}``average'' and {}``double integral'' respectively). We recall
that the set of medians of $f$ on $E$ consists of all numbers $c\in\mathbb{R}$
which satisfy \[
\lambda\left(\left\{ x\in E:f(x)<c\right\} \right)\le\frac{1}{2}\lambda(E)\mbox{ and }\lambda\left(\left\{ x\in E:f(x)>c\right\} \right)\le\frac{1}{2}\lambda(E)\,.\]
This set is always non empty, and the infimum in (\ref{eq:defv})
is attained, i.e., \begin{equation}
\mathbf{O}(f,E)=\frac{1}{\lambda(E)}\int_{E}\left|f-c\right|d\lambda\label{eq:newmib}\end{equation}
whenever $c$ is a median of $f$\,. We also recall that \begin{equation}
\mathbf{O}(f,E)\le\mathbf{A}(f,E)\le\mathbf{D}(f,E)\le2\,\mathbf{O}(f,E)\,\label{eq:avao}\end{equation}
for all functions $f$ which are integrable on $E$. For the reader's
convenience we recall the easy proofs of these standard facts in Appendix
\ref{sub:mediansandoscillation}. 
\begin{defn}
\label{def:bmode}Let $D$ be some measurable subset of $\mathbb{R}^{d}$
with positive measure and let $\mathcal{E}$ be some collection of
admissible subsets $E$ of $D$. We define the space $BMO(D,\mathcal{E})$
to consist of all (equivalence classes of) measurable functions $f:D\to\mathbb{R}$
for which the seminorm \begin{equation}
\left\Vert f\right\Vert _{BMO(D,\mathcal{E})}:=\sup_{E\in\mathcal{E}}\mathbf{O}(f,E)\label{eq:osn}\end{equation}
is finite. 
\end{defn}
One may also define $BMO(D,\mathcal{E})$ equivalently via either
one of the seminorms \begin{equation}
\left\Vert f\right\Vert _{BMO(D,\mathcal{E})}^{(\mathbf{A})}:=\sup_{E\in\mathcal{E}}\mathbf{A}(f,E)\label{eq:blonk}\end{equation}
or \begin{equation}
\left\Vert f\right\Vert _{BMO(D,\mathcal{E})}^{(\mathbf{D})}=\sup_{E\in\mathcal{E}}\mathbf{D}(f,E)\label{eq:fronk}\end{equation}
which (cf.~(\ref{eq:avao})) are each equivalent to the seminorm
(\ref{eq:osn}) to within constants of equivalence $1$ and $2$.

Of course if $f$ coincides a.e.~with a constant function then $\left\Vert f\right\Vert _{BMO(D,\mathcal{E})}=0$.
The reverse implication may also be true for suitable choices of $D$
and $\mathcal{E}$. In all cases the seminorm $\left\Vert \cdot\right\Vert _{BMO(D,\mathcal{E})}$
defines a norm on suitable equivalence classes of functions in $BMO(D,\mathcal{E})$
which may, for suitable choices of $D$ and $\mathcal{E}$, be simply
equivalence classes of functions modulo constants. 

The most frequently considered way of choosing $D$ and $\mathcal{E}$
is: 

\begin{equation}
\left\{ \begin{array}{l}
D\mbox{ is either }\mathbb{R}^{d}\mbox{ or some fixed cube in }\mathbb{R}^{d}\mbox{ and }\mathcal{E}\mbox{ is chosen }\\
\mbox{to be }\mathcal{Q}(D),\mbox{ the collection of all cubes contained in }D\mbox{.}\end{array}\right.\label{eq:jnc}\end{equation}
As the reader no doubt recalls, functions of the space $BMO(D,\mathcal{E})$
were first introduced and studied by John and Nirenberg \cite{JohnNirenberg}
for the case where $D$ is a cube in $\mathbb{R}^{d}$ and $\mathcal{E=Q}(D)$.
The original motivation for studying these functions apparently came
from John's study \cite{johnRandS} of problems in the theory of elasticity,
related in particular to the concept of elastic strain. One of the
first applications of \cite{JohnNirenberg} was in a paper \cite{moser}
by Moser extending Harnack's theorem about harmonic functions to functions
which are solutions of elliptic second order PDEs. But the space of
these functions and its analogues have since turned out to also have
many other deep properties and numerous other, sometimes quite surprising
applications in analysis. One particularly notable example of such
an application is the connection with $H^{p}$ spaces revealed in
the paper \cite{FS} of Fefferman and Stein.

The choice of $D$ and $\mathcal{E}$ specified in (\ref{eq:jnc})
is only one among several possible interesting choices, and we will
list four more examples of such choices now, taking the opportunity
to also fix our notation for them, notation which will be used thoughout
the paper. In each of these examples we will take the set $D$ to
either be $\mathbb{R}^{d}$ or some measurable subset of $\mathbb{R}^{d}$
with non empty interior.

\begin{equation}
\left\{ \begin{array}{l}
\mathcal{E}\mbox{ is chosen to be }\mathcal{D}(D),\mbox{ the collection of all dyadic cubes }\\
\mbox{contained in }D\mbox{ .}\end{array}\right.\label{eq:dyadic}\end{equation}

\begin{equation}
\left\{ \begin{array}{l}
\mathcal{E}\mbox{ is chosen to be }\mathcal{B}(D),\mbox{ the collection of all euclidean balls}\\
\mbox{contained in }D\mbox{ .}\end{array}\right.\label{eq:balls}\end{equation}

\begin{equation}
\left\{ \begin{array}{l}
\mathcal{E}\mbox{ is chosen to be }\mathcal{K}(D),\mbox{ the collection of all bounded closed}\\
\mbox{convex subsets of }D\mbox{ \mbox{which have non empty interiors}.}\end{array}\right.\label{eq:fedja}\end{equation}

\begin{equation}
\left\{ \begin{array}{l}
\mathcal{E}\mbox{ is chosen to be }\mathcal{W}(D),\mbox{ the collection of all \textit{special rectangles} }\\
\mbox{contained in }D\mbox{ .}\end{array}\right.\label{eq:wik0}\end{equation}

By special rectangles we mean all those subsets of $\mathbb{R}^{d}$
which are the cartesian products $I_{1}\times I_{2}\times...\times I_{d}$
of $d$ bounded closed intervals of positive length, where, for each
$j=1,2,....,d$, the length $\left|I_{j}\right|$ of $I_{j}$ equals
either $\min_{k=1,2,...,d}\left|I_{k}\right|$ or $2\min_{k=1,2,...,d}\left|I_{k}\right|$\,.
Such sets, and their associated space $BMO(\mathbb{R}^{d},\mathcal{W}(\mathbb{R}^{d}))$
were introduced and studied by Wik in \cite{wik}. He used the terminology
{}``false cubes'' for special rectangles. Below we will describe
his results in more detail.

Of course the seminorm $\left\Vert f\right\Vert _{BMO(D,\mathcal{K}(D))}$
is larger than any of the other seminorms $\left\Vert f\right\Vert _{BMO(D,\mathcal{E}))}$
arising from the other choices of $\mathcal{E}$ listed just above,
and for this reason it will be of less interest for us here for the
particular aims of this paper. However we remark that a result of
Nazarov, Sodin and Vol'berg (\cite{nsv} p.~13 and \cite{nsv2})
shows that every polynomial $P:\mathbb{R}^{d}\to\mathbb{R}$ of degree
$n$ satisfies \begin{equation}
\left\Vert \log\left|P\right|\right\Vert _{BMO(D,\mathcal{K}(\mathbb{R}^{d}))}\le\frac{4+\log4}{2}n\,.\label{eq:nsv}\end{equation}
It is remarkable that there is no dependence on the dimension $d$
in this inequality. We are naturally led to ask whether the left side
of (\ref{eq:nsv}) can also be bounded from below by $cn$ for some
absolute positive constant $c$. If this can be shown to be the case,
then other results in \cite{nsv} would imply that a dimension free
version of John-Nirenberg inequality holds, for $D=\mathbb{R}^{d}$
and $\mathcal{E=K}(\mathbb{R}^{d})$, at least for all functions of
the special form $\log\left|P\right|$. An analogous question with
analogous consequences can be asked for the apparently more difficult
and perhaps more interesting case where $\mathcal{E}$ is chosen to
be $\mathcal{Q}(\mathbb{R}^{d})$. 

There are also other more {}``exotic'' versions of the space $BMO$.
But these seem to be quite beyond the scope of what we will study
in this paper. For example, the measure $\lambda$ may be replaced
by a more general measure, and, furthermore, the underlying set $\mathbb{R}^{d}$
may be replaced by other suitable sets. There is even a version of
$BMO$ in the setting of martingales.

We have already alluded above to the following result in \cite{JohnNirenberg},
in fact the main result of that paper. Let $D$ be a cube in $\mathbb{R}^{d}$
and let $f:D\to\mathbb{R}$ be a function belonging to the space $BMO(D,\mathcal{Q}(D))$.
Then\begin{equation}
\begin{array}{c}
{\displaystyle \lambda\left(\left\{ x\in D:\left|f(x)-f_{D}\right|>\alpha\right\} \right)\le B\lambda(D)\exp\left(-\frac{b\alpha}{\left\Vert f\right\Vert _{BMO(D,\mathcal{Q}(D))}^{(\mathbf{A})}}\right)}\\
\mbox{for every }\alpha>0\,,\phantom{xxxxxxxxxxxxxxxxxxxxxxxxxxxxxxxxxxxxxxx}\end{array}\label{eq:jne}\end{equation}
where $B$ and $b$ are constants which depend only on the dimension
$d$. (In fact in the formulation of the main result (Lemma 1) of
\cite{JohnNirenberg} the number $f_{D}$ appearing on the right side
of (\ref{eq:jne}) is not explicitly chosen to be the average of $f$
on $D$. Nor is the functional $\left\Vert f\right\Vert _{BMO(D,\mathcal{Q}(D))}^{(\mathbf{A})}$
explicitly chosen for estimating the mean oscillation of $f$. However
such choices are made in Lemma 1' of the same paper which is used
to obtain Lemma 1.) 
\begin{rem}
\label{rem:optimallystrong}The result that (\ref{eq:jne}) holds
is an optimally strong result, in the following sense: Suppose that
$f:D\to\mathbb{R}$ is an arbitrary integrable function which satisfies
an estimate like (\ref{eq:jne}) for every subcube $Q$ of $D$, but
with some fixed constant independent of $Q$ in place of $\left\Vert f\right\Vert _{BMO(D,\mathcal{Q}(D))}^{(\mathbf{A})}$.
Then a simple calculation shows that $f\in BMO(D,\mathcal{Q}(D))$.
Conversely, the result of \cite{JohnNirenberg} obviously gives us
that a function $f\in BMO(D,\mathcal{Q}(D))$ satisifes an estimate
like (\ref{eq:jne}) for every subcube $Q$ of $D$.
\end{rem}
The inequality (\ref{eq:jne}), together with various generalizations
and variants of it, will be at once our main motivation and our main
interest in this paper. In fact (\ref{eq:jne}) is the key to obtaining
various other properties of $BMO$ and has been widely studied further
since its original discovery. The proof of (\ref{eq:jne}) in \cite{JohnNirenberg}
uses a famous lemma of Calderón and Zygmund. Among other proofs of
(\ref{eq:jne}), one of the simpler ones is due to Bennett, DeVore
and Sharpley \cite{bds} (see the remark at the end of Section 3 on
p.~607 of \cite{bds}) using a covering lemma which appeared previously
in \cite{bs}. 
\begin{rem}
\label{rem:bge1}It is a simple exercise to show that the constant
$B$ in (\ref{eq:jne}) must necessarily satisfy $B\ge1$. Furthermore,
versions of (\ref{eq:jne}) have been proved in which $B=2$. (See
e.g., \cite{stein} or the results of \cite{wik} which we shall also
discuss below.) 
\end{rem}

\begin{rem}
\label{rem:difver}It seems appropriate to make one specific remark
regarding the differing formats of results in quite a number of papers
which have followed on from \cite{JohnNirenberg}. These are papers
which present other proofs of (\ref{eq:jne}) or similar inequalities,
some of them generalizing to other settings and to other variants
of $BMO$. By {}``similar inequality'' we mean a version of (\ref{eq:jne})
which is modified in one or more of the following ways: For example
the seminorm $\left\Vert f\right\Vert _{BMO(D,\mathcal{Q}(D))}$ or
$\left\Vert f\right\Vert _{BMO(D,\mathcal{Q}(D))}^{(\mathbf{D})}$
may appear in place of $\left\Vert f\right\Vert _{BMO(D,\mathcal{Q}(D))}^{(\mathbf{A})}$.
Or the set \[
\left\{ x\in D:\left|f(x)-f_{D}\right|>\alpha\right\} \]
may be replaced by $\left\{ x\in D:\left|f(x)-f_{D}\right|\ge\alpha\right\} $,
or by $\left\{ x\in D:\left|f(x)-m\right|>\alpha\right\} $ or \[
\left\{ x\in D:\left|f(x)-m\right|\ge\alpha\right\} \]
 where $m$ is a median of $f$ on $D$. Our remark is this: It is
quite straightforward to check that such changes give an essentially
equivalent inequality. More precisely, if an inequality with any or
all of these changes holds for all $\alpha\ge0$ then, first of all,
this implies that $B\ge1$ (cf.~Remark \ref{rem:bge1}). But, furthermore,
this also implies that an inequality exactly of the form (\ref{eq:jne})
holds for all $\alpha\ge0$, but with possibly new positive constants
$b'$ and $B'$ replacing the positive constants $b$ and $B$. Conversely,
if (\ref{eq:jne}) holds for all $\alpha\ge0$, then each of the above
mentioned variants of (\ref{eq:jne}) also holds for all $\alpha\ge0$,
again possibly with new constants $b'$ and $B'$ replacing $b$ and
$B$. Moreover, if it is necessary to change $b$ and/or $B$ when
making any of these transitions from one version of (\ref{eq:jne})
to another, then the new constants $b'$ and $B'$ will always satisfy
$b/2\le b'\le b$ and $1\le B\le B'\le e^{b}B$. For the reader's
convenience, we include the calculations justifying these implications
in Appendix \ref{sub:dfojne}. 
\end{rem}
In this paper we begin the development of some techniques for studying
$BMO$ type spaces which are apparently somewhat different from those
used so far. We will use them to give an alternative proof of a somewhat
abstract result which includes, as special cases, both (\ref{eq:jne})
and an analogous result of Wik for special rectangles, which we will
describe in a moment. Our proof will not be obviously shorter or simpler
than the analogous proofs in \cite{JohnNirenberg} and \cite{bds}
and \cite{wik} and elsewhere. Nor does it give better constants than
the ones obtained by previously published proofs. But it seems distinctly
possible that, with further development and refinement, some elements
of our approach here may be able to give new information about the
behaviour of the constants $B$ and $b$, for large values of the
dimension $d$ and perhaps even ultimately lead to determining whether
these constants can both be taken to be independent of the dimension
$d$. In Wik's analogue of the inequality (\ref{eq:jne}) in the context
(\ref{eq:wik0}) of special rectangles (see (\ref{eq:wik1})), the
relevant constants are independent of $d$. Our proof in the same
context also gives constants independent of $d$, but Wik's constants
are better than ours. (See Remark \ref{rem:specialrectangles} for
details.)

Let us recall some of Wik's results more explicitly. For each measurable
function $f:\mathbb{R}^{d}\to\mathbb{R}$ he defines $\left\Vert f\right\Vert _{BMO}^{\prime}$
to be the seminorm (\ref{eq:fronk}) where $\mathcal{E=\mathcal{W}}(\mathbb{R}^{d})$
is the collection of all special rectangles (or {}``false cubes'')
in $\mathbb{R}^{d}$, i.e., in our notation \begin{equation}
\left\Vert f\right\Vert _{BMO}^{\prime}=\left\Vert f\right\Vert _{BMO(\mathbb{R}^{d},\mathcal{W}(\mathbb{R}^{d}))}^{(\mathbf{D})}\,.\label{eq:WikDef}\end{equation}
He defines $\left\Vert f\right\Vert _{BMO}$ analogously, except that
here $\mathcal{E}$ is $\mathcal{Q}(\mathbb{R}^{d})$, i.e., \[
\left\Vert f\right\Vert _{BMO}=\left\Vert f\right\Vert _{BMO(\mathbb{R}^{d},\mathcal{Q}(\mathbb{R}^{d}))}^{(\mathbf{D})}\,.\]
The particular interest of the seminorm $\left\Vert f\right\Vert _{BMO}^{\prime}$
lies in the following inequality which is proved by Wik:

\begin{equation}
\lambda\left(\left\{ x\in Q:\left|f(x)-m_{Q}\right|\ge\alpha\right\} \right)\le2\lambda(Q)\exp\left(-\frac{\alpha\ln2}{16\left\Vert f\right\Vert _{BMO}^{\prime}}\right)\,.\label{eq:wik1}\end{equation}
(We have written it here using notation slightly different from that
of \cite{wik}, to make it more convenient for comparison with (\ref{eq:jne}).)
This holds for every special rectangle $Q$, and for every $\alpha\ge0$
and for every number $m_{Q}$ which is a median of $f$ on $Q$. As
an immediate consequence of (\ref{eq:wik1}) and another result comparing
the seminorms $\left\Vert f\right\Vert _{BMO}^{\prime}$ and $\left\Vert f\right\Vert _{BMO}$,
Wik obtains a variant of (\ref{eq:jne}) which (in view of considerations
mentioned in the latter part of Remark \ref{rem:difver}, cf.~also
Lemma \ref{lem:av2med}) implies that the original inequality (\ref{eq:jne})
holds in fact for \[
b=\frac{\ln2}{32\left(2+6\sqrt{\frac{d}{\pi}}\right)}\mbox{ and \ensuremath{B=2e^{b}.}}\]

Here we list some of the features of our approach, some of which have
already been discussed or alluded to:

$\bullet$ Probably the most important feature of this paper is that
it provides the framework for posing {}``Question A'', whose positive
resolution would, as we show, give a dimension free John-Nirenberg
inequality.

$\bullet$ The {}``geometrical'' component of our proof (Theorem
\ref{thm:maingt}), or the more abstract hypothesis which can replace
it, has to be applied only once in the course of proving our versions
of (\ref{eq:jne}) (in contrast to some other known proofs of analogous
results). It is difficult to claim that this component is any simpler
than the Calderón-Zymund Lemma, or the covering lemma of \cite{bs}
p.~202. But there is perhaps more hope for strengthening it to a
dimension free version, than there is for doing away with the dependence
on dimension in approaches based on either of those two lemmata.

$\bullet$ In one of the decisive steps of our proof of the John-Nirenberg
inequality (see Theorem \ref{thm:Nultra}), we find that in some sense
we have reduced our argument to the case where we only have to consider
functions which take the three values $0$, $1$ and $2$.

$\bullet$ In some sense, we only have to consider the easy case where
$d=1$ and the relevant set $D$ in (\ref{eq:jne}) is an interval
and $f$ is a non increasing right continuous function on that interval.
The {}``geometrical'' result of Theorem \ref{thm:maingt}, which
gives an affirmative answer to a {}``dimension dependent'' version
of Question A, is the tool for reducing the general case of functions
of $d$ variables to this easy case.

$\bullet$ Instead of working with the seminorms\[
\left\Vert f\right\Vert _{BMO(D,\mathcal{Q}(D))}\mbox{ or \ensuremath{\left\Vert f\right\Vert _{BMO(D,\mathcal{Q}(D))}^{(\mathbf{D})}}or }\left\Vert f\right\Vert _{BMO(D,\mathcal{Q}(D))}^{(\mathbf{A})}\]
we mainly use another functional, which we denote by $\left\Vert f\right\Vert _{BMO_{0,s}}$
or $\left\Vert f\right\Vert _{BMO(D,\mathcal{Q}(D))}^{(\mathbf{J},s)}$.
This functional was introduced by John \cite{john} and then further
studied by Strömberg \cite{stromberg}. The condition sought in Question
A implies an inequality of the form \[
\left\Vert f^{*}\right\Vert _{BMO_{0,\sigma}(\mathbb{R})}\le\left\Vert f\right\Vert _{BMO_{0,s}(\mathbb{R}^{d})}\]
for suitable values of the parameter $\sigma$ and a suitable class
of functions $f$. Conversely (see Theorem \ref{thm:Ninverse}) if
such an inequality holds for some other appropriate value of $\sigma$,
then it implies the condition sought in Question A. 

$\bullet$ Our approach gives a version of the {}``dimension free''
result of Wik \cite{wik} for special rectangles.

\section{\label{sec:SignedRearrangements}Properties of non increasing rearrangements}

In this section we shall recall some properties of the non increasing
rearrangements of measurable functions which are defined on an arbitrary
measure space $\left(\Omega,\Sigma,\mu\right)$. Most, indeed probably
all of these are well known. A detailed discussion of them can be
found, for example, in \cite{hunt}. Among other relevant references
we mention \cite{bl} and \cite{sagher}. 

For each measurable $f:\Omega\to\mathbb{R}$, one first defines the
distribution function $f_{*}:(0,\infty)\to[0,\infty]$ of $f$ by
\[
f_{*}(\alpha)=\mu\left(\left\{ \omega\in\Omega:\left|f(\omega)\right|>\alpha\right\} \right)\,.\]
One can then define the non increasing rearrangement $f^{*}:(0,\infty)\to[0,\infty)$
of $f$, provided $f_{*}(\alpha)$ is finite for some positive $\alpha$.
It is given by the formula \[
f^{*}(t)=\inf\left\{ \alpha>0:f_{*}(\alpha)\le t\right\} \]
 for each $t>0$. It is, roughly speaking, the right continuous {}``inverse''
of the distribution function. 
\begin{rem}
In all our applications here, we will only need to consider the non
increasing rearrangements of functions for which the set on which
they are non zero has finite measure. Thus the required condition
about the finiteness of the distribution function will always be fulfilled. 
\end{rem}
Here are the properties of the non increasing rearrangement that we
will need in this paper. They all hold for any measurable function
$f:\Omega\to\mathbb{R}$ for which $f_{*}(\alpha)<\infty$ for some
$\alpha>0$. Recall that we denote the one dimensional Lebesgue measure
of subsets $G$ of $(0,\infty)$ by $\left|G\right|$.

(i) $f^{*}$ is non negative, non increasing and right continuous
on $\left(0,\infty\right)$.

(ii) $f$ and $f^{*}$ have the same distribution functions, i.e.,
they satisfy

\begin{equation}
\left|\left\{ t>0:f^{*}(t)>\alpha\right\} \right|=\mu\left(\left\{ \omega\in\Omega:\left|f(\omega)\right|>\alpha\right\} \right)\mbox{ for all }\alpha\in[0,\infty)\,.\label{eq:df}\end{equation}

(iii) If $f_{*}(\beta)<\infty$ for some $\beta\ge0$ then, \begin{equation}
\left|\left\{ t>0:f^{*}(t)\ge\alpha\right\} \right|=\mu\left(\left\{ \omega\in\Omega:\left|f(\omega)\right|\ge\alpha\right\} \right)\mbox{ for all }\alpha\in(\beta,\infty)\,\label{eq:df2}\end{equation}
 and \begin{equation}
\left|\left\{ t>0:f^{*}(t)=\alpha\right\} \right|=\mu\left(\left\{ \omega\in\Omega:\left|f(\omega)\right|=\alpha\right\} \right)\mbox{ for all }\alpha\in(\beta,\infty)\,\label{eq:df3}\end{equation}

(iv) The set $\left\{ t>0:f^{*}(t)>\alpha\right\} $ is the open interval
$\left(0,f_{*}(\alpha)\right)$ for each $\alpha\in[0,\infty)$.

(v) If $\mu(\Omega)<\infty$, then a variant of (\ref{eq:df3}) holds
for all $\alpha\in[0,\infty)$, namely \begin{equation}
\mu\left(\left\{ \omega\in\Omega:\left|f(\omega)\right|=\alpha\right\} \right)=\left|\left\{ t\in(0,\mu(\Omega)):f^{*}(t)=\alpha\right\} \right|\,.\label{eq:df4}\end{equation}

We refer, e.g.~to \cite{hunt} for proofs of properties (i) and (ii).
We can easily deduce (\ref{eq:df2}) from (\ref{eq:df}) with the
help of a sequence of numbers $\left\{ \alpha_{n}\right\} _{n\in\mathbb{N}}$
which satisfies $\beta<\alpha_{n}<\alpha_{n+1}<\alpha$ for each $n$
and also $\lim_{n\to\infty}\alpha_{n}=\alpha$. We have \[
\left|\left\{ t>0:f^{*}(t)>\alpha_{n}\right\} \right|=\mu\left(\left\{ \omega\in\Omega:\left|f(\omega)\right|>\alpha_{n}\right\} \right)<\infty\]
for each $n$ and we can apply the contracting sequence theorem. We
can then immediately obtain (\ref{eq:df3}) by subtracting (\ref{eq:df})
from (\ref{eq:df2}).

Now let us check that property (iv) holds. Since $f^{*}$ is non increasing,
the set \[
\left\{ t>0:f^{*}(t)>\alpha\right\} \]
must be an interval whose left endpoint is $0$. In view of property
(ii) the right endpoint of this interval must be $f_{*}(\alpha)$.
If this interval is unbounded then of course it is open. If it is
bounded and if $\alpha>0$ then the continuity from the right of $f^{*}$
implies that this interval cannot contain its right endpoint. If $\alpha=0$
then the interval is the union of the sequence of open intervals $\left\{ t>0:f^{*}(t)>1/n\right\} $,
$n=1,2,....$ and is therefore also open.

It follows from (iv), and the fact that $f_{*}(\alpha)\le\mu(\Omega)$
for all $\alpha\in[0,\infty)$, that \begin{equation}
\left\{ t>0:f^{*}(t)>\alpha\right\} =(0,\mu(\Omega))\cap\left\{ t>0:f^{*}(t)>\alpha\right\} =\left\{ t\in(0,\mu(\Omega)):f^{*}(t)>\alpha\right\} \label{eq:wwn-1}\end{equation}
 also holds for every $\alpha\in[0,\infty)$\,.

Here now is the proof of (v). We have $\mu(\Omega)<\infty$, which
enables us to obtain (\ref{eq:df4}) in the case where $\alpha=0$,
by first using (\ref{eq:df}) and then (\ref{eq:wwn-1}), as follows:
\begin{eqnarray*}
\mu\left(\left\{ \omega\in\Omega:\left|f(\omega)\right|=0\right\} \right) & = & \mu\left(\Omega\right)-\mu\left(\left\{ \omega\in\Omega:\left|f(\omega)\right|>0\right\} \right)\\
 & = & \left|\left(0,\mu(\Omega)\right)\right|-\left|\left\{ t>0:f^{*}(t)>0\right\} \right|\\
 & = & \left|\left(0,\mu(\Omega)\right)\right|-\left|\left\{ t\in(0,\mu(\Omega)):f^{*}(t)>0\right\} \right|\\
 & = & \left|\left\{ t\in(0,\mu(\Omega)):f^{*}(t)=0\right\} \right|\,.\end{eqnarray*}
 We easily obtain (\ref{eq:df4}) in the remaining case where $\alpha>0$
by first applying (\ref{eq:df3}) using $\beta=0$, and then observing
that in this case \[
\left\{ t>0:f^{*}(t)=\alpha\right\} \subset\left\{ t>0:f^{*}(t)>0\right\} \subset\left(0,\mu(\Omega)\right)\,.\]

We close this section with the following lemma which will be needed
in Section \ref{sec:JohnStrombergFunctional}. 
\begin{lem}
\label{lem:psm}Let $Q$ be an admissible subset of $\mathbb{R}^{d}$.
Let $g:Q\to[0,\infty)$ be a measurable function. Then, \begin{eqnarray}
\left|\left\{ t\in(0,\lambda(Q)):\left|g^{*}(t)-c\right|\le\alpha\right\} \right| & = & \lambda\left(\left\{ x\in Q:\left|g(x)-c\right|\le\alpha\right\} \right)\nonumber \\
\label{eq:z}\\\mbox{for all }c\in\mathbb{R}\mbox{ and all }\alpha\ge0\,.\phantom{sss}\,\,\,\,\,\nonumber \end{eqnarray}

\end{lem}
\noindent \textit{Proof.} The set which appears on the left side
of (\ref{eq:z}) coincides with \[
\left\{ t\in(0,\lambda(Q)):c-\alpha\le g^{*}(t)\le c+\alpha\right\} \]
and the set appearing on the right side of (\ref{eq:z}) coincides
with \[
\left\{ x\in Q:c-\alpha\le g(x)\le c+\alpha\right\} \,.\]
If $c+\alpha<0$ then both of these sets are empty. This is because
$g^{*}$ is non negative by definition, and because we have imposed
the condition that $f$ is non negative and therefore so is $g$.
Accordingly, we only have to prove (\ref{eq:z}) in the case where
$c+\alpha\ge0$. In that case the above-mentioned two sets coincide
respectively with $\left\{ t\in(0,\lambda(Q)):\gamma\le g^{*}(t)\le\delta\right\} $
and $\left\{ x\in Q:\gamma\le g(x)\le\delta\right\} $, where $\gamma=\max\left\{ 0,c-\alpha\right\} $
and $\delta=c+\alpha$. This reduces the proof of (\ref{eq:z}) to
showing that \begin{equation}
\left|\left\{ t\in(0,\lambda(Q)):\gamma\le g^{*}(t)\le\delta\right\} \right|=\lambda\left(\left\{ x\in Q:\gamma\le g(x)\le\delta\right\} \right)\mbox{ whenever }0\le\gamma\le\delta\,.\label{eq:gamdel}\end{equation}

In order to check that (\ref{eq:gamdel}) holds we first note that
\[
\left\{ x\in Q:\gamma\le g(x)\le\delta\right\} =\left\{ x\in Q:g(x)\le\delta\right\} \backslash\left\{ x\in Q:g(x)<\gamma\right\} \]
and \[
\left\{ x\in Q:g(x)<\gamma\right\} \subset\left\{ x\in Q:g(x)\le\delta\right\} \,.\]
Therefore, \[
\left|\left\{ x\in Q:\gamma\le g(x)\le\delta\right\} \right|=\left|\left\{ x\in Q:g(x)\le\delta\right\} \right|-\left|\left\{ x\in Q:g(x)<\gamma\right\} \right|\,.\]
Analogously, we have \begin{eqnarray*}
 &  & \left|\left\{ t\in(0,\lambda(Q)):\gamma\le g^{*}(t)\le\delta\right\} \right|\\
 & = & \left|\left\{ t\in(0,\lambda(Q)):g^{*}(t)\le\delta\right\} \right|-\left|\left\{ t\in(0,\lambda(Q)):g^{*}(t)<\gamma\right\} \right|\,.\end{eqnarray*}
 So it will suffice to show that \begin{equation}
\left|\left\{ t\in(0,\lambda(Q)):g^{*}(t)\le\delta\right\} \right|=\lambda\left(\left\{ x\in Q:g(x)\le\delta\right\} \right)\label{eq:fordel}\end{equation}
 and \begin{equation}
\left|\left\{ t\in(0,\lambda(Q)):g^{*}(t)<\gamma\right\} \right|=\lambda\left(\left\{ x\in Q:g(x)<\gamma\right\} \right)\label{eq:forgam}\end{equation}
 We have that \[
\lambda\left(\left\{ x\in Q:g(x)\le\delta\right\} \right)=\lambda(Q)-\lambda\left(\left\{ x\in Q:g(x)>\delta\right\} \right)\,.\]
Then, by property (i) of non increasing rearrangements (one among
several such numbered properties presented at the beginning of this
section), since $g$ is non negative, this last expression equals
$\lambda(Q)-\left|\left\{ t>0:g^{*}(t)>\delta\right\} \right|$ which
in turn equals \begin{eqnarray*}
\lambda(Q)-\left|\left\{ t\in(0,\lambda(Q)):g^{*}(t)>\delta\right\} \right| & = & \left|\left\{ t\in(0,\lambda(Q)):g^{*}(t)\le\delta\right\} \right|\end{eqnarray*}
which establishes (\ref{eq:fordel}). If $\gamma=0$ then we immediately
obtain (\ref{eq:forgam}) since empty sets have measure $0$. In the
remaining case, where $\gamma>0$, we can first use property (iii)
of non increasing rearrangements (and again the non-negativity of
$g$) to obtain that \[
\lambda\left(\left\{ x\in Q:g(x)\ge\gamma\right\} \right)=\left|\left\{ t>0:g^{*}(t)\ge\gamma\right\} \right|=\left|\left\{ t\in(0,\lambda(Q)):g^{*}(t)\ge\gamma\right\} \right|\,.\]
This, combined with the facts that \[
\lambda\left(\left\{ x\in Q:g(x)<\gamma\right\} \right)=\lambda(Q)-\lambda\left(\left\{ x\in Q:g(x)\ge\gamma\right\} \right)\]
 and \[
\left|\left\{ t\in(0,\lambda(Q)):g^{*}(t)<\gamma\right\} \right|=\lambda(Q)=\left|\left\{ t\in(0,\lambda(Q)):g^{*}(t)\ge\gamma\right\} \right|\]
 immediately gives us (\ref{eq:forgam}) and so completes the proof
of (\ref{eq:gamdel}), which, as already explained, also completes
the proof of (\ref{eq:z}) and therefore of the lemma. $\qed$

\section{\label{sec:JohnStrombergFunctional}The functional of John and Strömberg
for characterizing $BMO$.}

Given an admissible subset $E$ of $\mathbb{R}^{d}$, a real valued
function $f$ which is defined and measurable on $E$, and a number
$s\in(0,1)$, it is convenient to introduce the notation $\mathbf{J}(f,E,s)$
for a special functional which was introduced and studied in\cite{john}
and then considered in greater generality in\cite{stromberg}. Thus
we set \begin{equation}
\mathbf{J}(f,E,s)=\inf_{c\in\mathbb{R}}\left(\inf\left\{ \alpha\ge0:\lambda\left(\left\{ x\in E:\left|f(x)-c\right|>\alpha\right\} \right)<s\lambda(E)\right\} \right)\,.\label{eq:OrigDefn}\end{equation}
 (Here, as always in this paper, $\lambda$ denotes $d$-dimensional
Lebesgue measure on $\mathbb{R}^{d}$.) In \cite{john} and \cite{stromberg}
the set $E$ is always taken to be a cube, and the functional $\mathbf{J}(f,E,s)$
is shown to be a kind of counterpart, a very useful counterpart, of
the functionals $\mathbf{O}(f,E)$, $\mathbf{A}(f,E)$ and $\mathbf{D}(f,E)$.

There is another, perhaps more convenient formula for $\mathbf{J}(f,E,s)$,
namely \[
\mathbf{J}(f,E,s)=\inf_{c\in\mathbb{R}}\left((f-c)\chi_{E}\right)^{*(L)}(s\lambda(E))\,.\]
 (Here $u^{*(L)}$ denotes the \textit{left} continuous rearrangement
of a measurable function $u$.) This formula is mentioned e.g., in
\cite{lerner98}, \cite{YoramPasha98} and \cite{pasha99}. In this
section we shall obtain yet another formula for $\mathbf{J}(f,E,s)$
in terms of rearrangements. (See Proposition \ref{pro:prop}).

In our case $E$ will often be a cube or special rectangle, or, more
generally, a member of some collection $\mathcal{E}$ of admissible
subsets which is used, as in (\ref{eq:osn}), together with some measurable
set $D\subset\mathbb{R}^{d}$, to define a seminorm for some version
of the space $BMO$. Indeed for such $\mathcal{E}$ and $D$, following
the model of \cite{john} and \cite{stromberg}, and analogously to
the seminorm defined by (\ref{eq:osn}), we consider the functional
\begin{equation}
\left\Vert f\right\Vert _{BMO(D,\mathcal{E})}^{(\mathbf{J},s)}:=\sup_{E\in\mathcal{E}}\mathbf{J}(f,E,s)\,.\label{eq:rtp}\end{equation}

\begin{rem}
In the case where $D=\mathbb{R}^{d}$ and $\mathcal{E=}\mathcal{Q}(\mathbb{R}^{d})$
it is known \cite{JawerthTorchinsky,pasha99} that this quantity is
equivalent to a certain $K$-functional. More explicitly, \[
\left\Vert f\right\Vert _{BMO(D,\mathcal{E})}^{(\mathbf{J},e^{-t})}\sim K(t,f;L^{\infty}(\mathbb{R}^{d}),BMO(\mathbb{R}^{d},\mathcal{Q}(\mathbb{R}^{d}))\,.\]

\end{rem}
Despite the choice of notation in (\ref{eq:rtp}), $\left\Vert f\right\Vert _{BMO(D,\mathcal{E})}^{(\mathbf{J},s)}$
is not a norm nor even a seminorm. At least it is homogeneous, i.e.,
as follows almost immediately from the definition, \begin{equation}
\mathbf{J}(rf,E,s)=\left|r\right|\mathbf{J}(f,E,s)\mbox{ and so }\left\Vert rf\right\Vert _{BMO(D,\mathcal{E})}^{(\mathbf{J},s)}=\left|r\right|\left\Vert f\right\Vert _{BMO(D,\mathcal{E})}^{(\mathbf{J},s)}\mbox{ for each }r\in\mathbb{R}\,.\label{eq:homog}\end{equation}
Let us note another simple property of these functionals: If $T$
is an invertible affine transformation of $\mathbb{R}^{d}$, i.e.,
if $Tx=rx+x_{0}$ for some non zero $r\in\mathbb{R}$ and $x_{0}\in\mathbb{R}^{d}$,
and if $g(x)=$$f(rx+x_{0})$, then a simple routine calculation (see
Appendix \ref{sub:affinecommute}) shows that \begin{equation}
\mathbf{J}(g,E,s)=\mathbf{J}(f,rE+x_{0},s)\label{eq:preaffine}\end{equation}
 for each admissible set $E$ contained in the domain of $g$. Consequently
\begin{equation}
\left\Vert g\right\Vert _{BMO(D,\mathcal{E})}^{(\mathbf{J},s)}=\left\Vert f\right\Vert _{BMO(T(D),T(\mathcal{E}))}^{(\mathbf{J},s)}\label{eq:affinestuff}\end{equation}
 where $T\left(\mathcal{E}\right)$ is of course the collection of
sets $\left\{ T(E):E\in\mathcal{E}\right\} $. In various natural
examples, where $D=\mathbb{R}^{d}$ and $\mathcal{E}$ is any one
of the collections $\mathcal{Q}(\mathbb{R}^{d})$, $\mathcal{D}(\mathbb{R}^{d})$,
$\mathcal{B}(\mathbb{R}^{d})$, $\mathcal{K}(\mathbb{R}^{d})$ or
$\mathcal{W}(\mathbb{R}^{d})$ we of course have $T(D)=D$ and $T(\mathcal{E})=\mathcal{E}$\,.

Suppose that $D=\mathbb{R}^{d}$ and (as in (\ref{eq:jnc})) $\mathcal{E}$
is the collection $\mathcal{Q}(\mathbb{R}^{d})$ of all cubes in $\mathbb{R}^{d}$.
In this case it will sometimes be convenient to adopt the notation
of \cite{stromberg} and write\begin{equation}
\left\Vert f\right\Vert _{BMO_{0,s}}=\left\Vert f\right\Vert _{BMO(\mathbb{R}^{d},\mathcal{Q}(\mathbb{R}^{d}))}^{(\mathbf{J},s)}\label{eq:pone}\end{equation}
 and also \[
\left\Vert f\right\Vert _{BMO}=\left\Vert f\right\Vert _{BMO(\mathbb{R}^{d},\mathcal{Q}(\mathbb{R}^{d}))}\,.\]
 It is known that \begin{equation}
s\left\Vert f\right\Vert _{BMO_{0,s}}\le\left\Vert f\right\Vert _{BMO}\le C_{d}\left\Vert f\right\Vert _{BMO_{0,s}}\label{eq:sbjo}\end{equation}
 for every measurable $f:\mathbb{R}^{d}\to\mathbb{R}$, whenever $0<s\le\frac{1}{2}$,
where $C_{d}$ is a constant depending only on the dimension $d$.
This result was originally obtained by John \cite{john} for $0<s<\frac{1}{2}$,
and then extended by Strömberg \cite{stromberg} to include the case
$s=\frac{1}{2}$. Thus the functional $\mathbf{J}(f,Q,s)$ enables
one to characterize $BMO$ functions in an alternative way.

The result (\ref{eq:sbjo}) is false for $s>1/2,$ although the definition
(\ref{eq:OrigDefn}) is valid for all $s\in(0,1)$. This is because
$\mathbf{J}(f,E,s)=0$ and $\left\Vert f\right\Vert _{BMO_{0,s}}=0$
for certain non constant functions $f$ whenever $s>1/2$. (Cf.~the
remark on p.~522 of \cite{stromberg}.) 
\begin{rem}
\label{rem:Cheb}The essential content of (\ref{eq:sbjo}) is the
second inequality. Let us recall the elementary proof of (a more general
version of) the first inequality of (\ref{eq:sbjo}). By Chebyshev's
inequality we have \[
\lambda\left(\left\{ x\in E:\left|f(x)-c\right|>\alpha\right\} \right)\le\frac{1}{\alpha}\int_{E}\left|f-c\right|d\lambda=\frac{\lambda(E)}{\alpha}\mathbf{O}(f,E)\]
for each admissible $E$, each $\alpha>0$, each $f$ which is measurable
on $E$, and each median $c$ of $f$ on $E$. Thus, every $\alpha$
satisfying $\alpha>\frac{1}{s}\mathbf{O}(f,E)$ also satisfies \[
\lambda\left(\left\{ x\in E:\left|f(x)-c\right|>\alpha\right\} \right)<s\lambda(E)\]
for some $c\in\mathbb{R}$. Accordingly, \[
\mathbf{J}(f,E,s)\le\frac{1}{s}\mathbf{O}(f,E)\]
which immediately implies that \begin{equation}
\left\Vert f\right\Vert _{BMO(D,\mathcal{E})}^{(\mathbf{J},s)}\le\frac{1}{s}\left\Vert f\right\Vert _{BMO(D,\mathcal{E})}.\label{eq:talpt}\end{equation}
The first inequality in (\ref{eq:sbjo}) is a special case of (\ref{eq:talpt}). 
\end{rem}

The method which we develop in this paper will obviously imply an
alternative proof of (\ref{eq:sbjo}), but (so far) only for quite
small values of $s$.

We will sometimes need to use the following very simple result. 
\begin{lem}
\label{lem:JLIP}Suppose that the function $\varphi:\mathbb{R}\to\mathbb{R}$
satisfies \[
\left|\varphi(s)-\varphi(t)\right|\le\left|s-t\right|\]
for all $s,t\in\mathbb{R}$. Then

\[
\mathbf{J}(\varphi\circ f,E,s)\le\mathbf{J}(f,E,s)\]
for every admissible set $E$, every $s\in(0,1)$, and every real
valued function $f$ which is defined and measurable on $E$. 
\end{lem}
\noindent \textit{Proof.} This follows immediately from the obvious
inclusion \[
\left\{ \left\{ x\in E:\left|\varphi\left(f(x)\right)-\varphi(c)\right|>\alpha\right\} \right\} \subset\left\{ x\in E:\left|f(x)-c\right|>\alpha\right\} \]
 and the definition of $\mathbf{J}(f,E,s)$. $\qed$ 
\begin{rem}
An analogous result holds for the usual $BMO$ seminorm and functional
$\mathbf{O}(f,E)$. See Lemma \ref{lem:lipz} in Appendix \ref{sub:LipschitzAndBMO}. 
\end{rem}
We remark that, for each $f$ , $E$ and $s$ as above, and for each
$c\in\mathbb{R}$ and each $\alpha\ge0$, the condition \[
\lambda\left(\left\{ x\in E:\left|f(x)-c\right|>\alpha\right\} \right)<s\lambda(E)\]
 is equivalent to \[
\lambda\left(\left\{ x\in E:\left|f(x)-c\right|\le\alpha\right\} \right)>(1-s)\lambda(E)\,.\]
So we also have \begin{equation}
\mathbf{J}(f,E,s)=\inf_{c\in\mathbb{R}}\left(\inf\left\{ \alpha\ge0:\lambda\left(\left\{ x\in E:\left|f(x)-c\right|\le\alpha\right\} \right)>(1-s)\lambda(E)\right\} \right)\,.\label{eq:altdef}\end{equation}

The following proposition gives us another way to calculate and {}``visualize''
$\mathbf{J}(f,E,s)$, at least for functions which are either univariate
and monotone, or non negative. This other way, for some purposes,
seems to be an easier alternative than working with the original definition.
It enables us to work with just one variable (here denoted by $u$),
instead of having to deal with the two variables $\alpha$ and $c$
in the original definition. 
\begin{prop}
\label{pro:prop}

(i) For each $q>0$ and each non increasing right continuous function
$h:(0,q)\to\mathbb{R}$, the formula \begin{equation}
\mathbf{J}(h,(0,q),s)=\frac{1}{2}\inf\left\{ h(u)-h\left(u+(1-s)q\right):0<u<sq)\right\} \label{eq:preffnn}\end{equation}
 holds for each $s\in(0,1)$.

(ii) Furthermore, the formula \begin{equation}
\mathbf{J}(f,Q,s)=\frac{1}{2}\inf\left\{ \left(f\chi_{Q}\right)^{*}(u)-\left(f\chi_{Q}\right)^{*}\left(u+(1-s)\lambda(Q)\right):0<u<s\lambda(Q)\right\} \label{eq:ffnn}\end{equation}
holds for each admissible subset $Q$ of $\mathbb{R}^{d}$, each $s\in(0,1)$
and each \textbf{non negative} real valued function $f$ which is
defined and measurable on $Q$. \end{prop}
\begin{rem}
In our main applications of this proposition the set $Q$ will be
a cube or a special rectangle. But we stress that, despite the choice
of letter, the set $Q$ in (\ref{eq:ffnn}) can be an \textit{arbitrary}
admissible subset. 
\end{rem}

\begin{rem}
\label{rem:mtjm}Restated informally, part (ii) of this proposition
tells us that $2\mathbf{J}(f,Q,s)$ is the {}``minimum'' amount
that $\left(f\chi_{Q}\right)^{*}$ can decrease on any closed subinterval
of $\left(0,\lambda(Q)\right)$ of length exactly $(1-s)\lambda(Q)$. 
\end{rem}

\begin{rem}
It is easy to see from the original definition or from the formula
(\ref{eq:ffnn}), that, for each fixed $Q$ and $f$ the function
$s\mapsto\mathbf{J}(f,Q,s)$ is non increasing. As is explained in
Appendix \ref{sub:lrconj} (but is not needed for any other purposes
in this paper), $s\mapsto\mathbf{J}(f,Q,s)$ is also left continuous,
but in general not right continuous. 
\end{rem}
\noindent \textit{Proof of Proposition \ref{pro:prop}.} We will first
deal with part (i). (To understand the rather simple ideas behind
our proof of (\ref{eq:preffnn}), the reader may care to first look
at the rather shorter and simpler proof given below in Remark \ref{rem:fyed}
for the special case where $h$ is strictly decreasing and uniformly
continuous on $(0,q)$, and to draw some relevant pictures of the
graph of $h$.)

Let $\beta$ equal the right side of (\ref{eq:preffnn}). We will
now prove one {}``half'' of (\ref{eq:preffnn}), namely that $\mathbf{J}(h,(0,q),s)\le\beta$.
Obviously $\beta\ge0$ and there exists a non increasing sequence
$\left\{ \beta_{n}\right\} _{n\in\mathbb{N}}$ which tends to $\beta$
and a sequence $\left\{ u_{n}\right\} _{n\in\mathbb{N}}$ of numbers
satisfying $0<u_{n}<sq$ such that \[
\beta_{n}=\frac{1}{2}\left(h(u_{n})-h(u_{n}+(1-s)q)\right)\,.\]
Since $u_{n}+(1-s)q<q$ and $h$ is right continuous, there exists
$v_{n}$ such that $u_{n}+(1-s)q<v_{n}<q$ and \[
0\le h\left(u_{n}+(1-s)\lambda(Q)\right)-h(v_{n})\le\frac{1}{n}\,.\]
If we set $c_{n}=\frac{1}{2}\left(h(u_{n})+h(v_{n})\right)$ and $\alpha_{n}=\frac{1}{2}\left(h(u_{n})-h(v_{n})\right)$
then \begin{eqnarray*}
[u_{n},v_{n}] & \subset & \left\{ t\in(0,q):h(v_{n})\le h(t)\le h(u_{n})\right\} \\
 & = & \left\{ t\in(0,q):c_{n}-\alpha_{n}\le h(t)\le c_{n}+\alpha_{n}\right\} \\
 & = & \left\{ t\in(0,q):\left|h(t)-c_{n}\right|\le\alpha_{n}\right\} .\end{eqnarray*}
 It follows that \[
\left|\left\{ t\in(0,q):\left|h(t)-c_{n}\right|\le\alpha_{n}\right\} \right|\ge v_{n}-u_{n}>(1-s)q\,.\]
 Consequently (by (\ref{eq:altdef})) we have $\mathbf{J}(h,(0,q),s)\le\alpha_{n}$
for each $n$. Since \[
\lim_{n\to\infty}\alpha_{n}=\lim_{n\to\infty}\beta_{n}=\beta\]
this shows that $\mathbf{J}(h,(0,q),s)\le\beta$.

Next we shall prove the reverse of the preceding inequality, namely
that $\beta\le\mathbf{J}(h,(0,q),s)$. Here again we will use sequences
denoted by $\left\{ \alpha_{n}\right\} _{n\in\mathbb{N}}$, $\left\{ c_{n}\right\} _{n\in\mathbb{N}}$,
$\left\{ u_{n}\right\} _{n\in\mathbb{N}}$ and $\left\{ v_{n}\right\} _{n\in\mathbb{N}}$.
But they will be defined differently from their definitions in the
preceding part of the proof. By (\ref{eq:altdef}), there exists a
non increasing sequence $\left\{ \alpha_{n}\right\} _{n\in\mathbb{N}}$
of non negative numbers which tends to $\mathbf{J}(h,(0,q),s)$ and
a sequence $\left\{ c_{n}\right\} _{n\in\mathbb{N}}$ of real numbers
such that \begin{equation}
\left|\left\{ t\in(0,q):c_{n}-\alpha_{n}\le h(t)\le c_{n}+\alpha_{n}\right\} \right|>(1-s)q\,.\label{eq:feup}\end{equation}

Let us define \[
u_{n}:=\inf\left\{ t\in\left(0,q\right):h(t)\le c_{n}+\alpha_{n}\right\} \]
 and \[
v_{n}:=\sup\left\{ t\in\left(0,q\right):h(t)\ge c_{n}-\alpha_{n}\right\} \,.\]
 Then, by definition, for each $m\in\mathbb{N}$, we have that \[
\left[u_{n}+1/m,v_{n}-1/m\right]\subset\left\{ t\in\left(0,q\right):c_{n}-\alpha_{n}\le h(t)\le c_{n}+\alpha_{n}\right\} \subset[u_{n},v_{n}]\cap\left(0,q\right)\,.\]
 Since we can choose $m$ arbitrarily large, this implies that the
intervals \[
\left(u_{n},v_{n}\right)\mbox{ and }[u_{n},v_{n}]\cap\left(0,q\right)\]
 must have the same length as the interval \[
\left\{ t\in\left(0,q\right):c_{n}-\alpha_{n}\le h(t)\le c_{n}+\alpha_{n}\right\} \,.\]
 In view of (\ref{eq:feup}), this gives us that $v_{n}-u_{n}>(1-s)q$.
Furthermore, $0\le u_{n}<v_{n}\le q$. Therefore, for some sufficiently
small $\varepsilon_{n}>0$, we have $v_{n}-u_{n}>2\varepsilon_{n}+(1-s)q$
and \[
0\le u_{n}<u_{n}+\varepsilon_{n}<u_{n}+\varepsilon_{n}+(1-s)q<u_{n}+\varepsilon_{n}+v_{n}-u_{n}-2\varepsilon_{n}=v_{n}-\varepsilon_{n}<v_{n}\le q\,.\]
 Since the two points $u_{n}+\varepsilon_{n}$ and $u_{n}+\varepsilon_{n}+(1-s)q$
are both in $(0,q)$, the number $\beta$ defined above satisfies
\begin{eqnarray*}
2\beta & \le & h\left(u_{n}+\varepsilon_{n}\right)-h\left(u_{n}+\varepsilon_{n}+(1-s)q\right)\\
 & \le & h\left(u_{n}+\varepsilon_{n}\right)-h\left(v_{n}-\varepsilon_{n}\right).\end{eqnarray*}
By the definitions of $u_{n}$ and $v_{n}$ this last expression is
dominated by \[
c_{n}+\alpha_{n}-\left(c_{n}-\alpha_{n}\right)=2\alpha_{n}\,.\]
Thus $\beta\le\alpha_{n}$ for all $n$. This gives us the remaining
required inequality $\beta\le\mathbf{J}(h,(0,q),s)$ and completes
the proof of (\ref{eq:preffnn}) and part (i) of the proposition.

Now we turn to part (ii) and the proof of the formula (\ref{eq:ffnn}).
We will see that in fact (\ref{eq:ffnn}) can be deduced from (\ref{eq:preffnn}),
essentially by a careful application of the fact that the functions
$f$ and $\left(f\chi_{Q}\right)^{*}$, when restricted to $Q$ and
to $\left(0,\lambda(Q)\right)$ respectively, have the same distribution
function.

The function $\left(f\chi_{Q}\right)^{*}$ is non increasing and right
continuous on $\left(0,\infty\right)$ and therefore also on the subinterval
$\left(0,\lambda(Q)\right)$. So, we can set $q=\lambda(Q)$ and $h=\left(f\chi_{Q}\right)^{*}$
and apply (\ref{eq:preffnn}) to obtain that \begin{eqnarray}
 &  & \mathbf{J}\left(\left(f\chi_{Q}\right)^{*},\left(0,\lambda(Q)\right),s\right)=\nonumber \\
\label{eq:did-1}\\ &  & \frac{1}{2}\inf\left\{ \left(f\chi_{Q}\right)^{*}(u)-\left(f\chi_{Q}\right)^{*}\left(u+(1-s)\lambda(Q)\right):0<u<s\lambda(Q)\right\} \,.\nonumber \end{eqnarray}

We remark that we have used the notation $\left(f\chi_{Q}\right)^{*}$
rather than $f^{*}$ in (\ref{eq:ffnn}) because, in future applications
of this proposition, $f$ might possibly be defined on all of $\mathbb{R}^{d}$
or on some other set which is strictly larger that $Q$. (Indeed the
statement of the proposition explicitly allows for this possibility.)
To simplify the notation in the rest of our proof we will let \[
g=f\mid_{Q}\,,\]
i.e., $g:Q\to[0,\infty)$ will denote the function defined \textit{only}
on $Q$ which is the restriction of $f$ to $Q$. Thus we can unambiguously
write $g^{*}$ instead of $\left(f\chi_{Q}\right)^{*}$, and of course
$\mathbf{J}(f,Q,s)=\mathbf{J}(g,Q,s)$. In view of (\ref{eq:did-1}),
in order to complete the proof of (\ref{eq:ffnn}) and part (ii) of
Proposition \ref{pro:prop}, it will suffice to show that \begin{equation}
\mathbf{J}(g,Q,s)=\mathbf{J}\left(g^{*},\left(0,\lambda(Q)\right),s\right)\,.\label{eq:gqgst}\end{equation}

In view of (\ref{eq:altdef}), we can immediately obtain (\ref{eq:gqgst})
if we know that\[
\left|\left\{ t\in(0,\lambda(Q)):\left|g^{*}(t)-c\right|\le\alpha\right\} \right|=\lambda\left(\left\{ x\in Q:\left|g(x)-c\right|\le\alpha\right\} \right)\]
for all $c\in\mathbb{R}$ and all $\alpha\ge0$. This is exactly the
result which was proved in Lemma \ref{lem:psm} and therefore the
proof of part (ii) of Proposition \ref{pro:prop} is complete. $\qed$
\begin{rem}
\label{rem:fyed}Here, as promised above, is the simpler proof of
(\ref{eq:preffnn}) for the case where $h$ is uniformly continuous,
and strictly decreasing. In this case $h$ has a unique extension
to a continuous function on $\left[0,q\right]$ which we will also
denote by $h$. For each pair of numbers $c\in\mathbb{R}$ and $\alpha\ge0$,
let \[
E(c,\alpha):=\left\{ t\in[0,q]:\left|h(t)-c|\le\alpha\right|\right\} =\left\{ t\in[0,q]:c-\alpha\le h(t)\le c+\alpha\right\} \,.\]
This set is clearly a closed interval $\left[u,u+r\right]$ contained
in $\left[0,q\right]$, on which $h$ attains a minimum value $m$
(at $u+r)$ and a maximum value $M$ (at $u$), and these values both
lie in the interval $\left[c-\alpha,c+\alpha\right]$. If we set $c'=\frac{1}{2}(M+m)$
and $\alpha'=\frac{1}{2}(M-m)$ then of course $E(c',\alpha')=E(c,\alpha)$
and $0\le\alpha'\le\alpha$. Of course the length $r$ of the interval
$E(c,\alpha)$ is the same as the length of the not necessarily closed
interval $\left\{ t\in(0,q):\left|h(t)-c|\le\alpha\right|\right\} $.
So, in order to calculate $\mathbf{J}(h,(0,q),s)$, we have to consider
all intervals $E(c,\alpha)$ which have length exceeding $(1-s)q$
and find the infimum of all values of $\alpha$ which they can have.
If, as above, we write $E(c,\alpha)$ as $[u,u+r]$, then $M=h(u)$
and $m=h(u+r)$ and \[
\alpha'=\frac{1}{2}(M-m)=\frac{1}{2}\left(h(u)-h(u+r)\right)\,.\]
Thus $\mathbf{J}(h,(0,q),s)$ is the infimum of the set $\Omega$
of all numbers $\frac{1}{2}\left(h(u)-h(u+r)\right)$ for which $r>(1-s)q$
and $0\le u\le u+r\le q$. In view of the continuity and monotonicity
of $h$, we can optimally choose $r=(1-s)q$, so that the above infimum
is equal to the infimum of the set $\Omega_{1}$ of all numbers $\frac{1}{2}\left(h(u)-h\left(u+(1-s)q\right)\right)$
for which $0\le u\le u+(1-s)q\le q$. (The infimum is of course attained
for some particular $u\in[0,q]$.) Again by continuity, this infimum
is also equal to the infimum of the subset $\Omega_{2}$ of $\Omega_{1}$
\[
\Omega_{2}=\left\{ \frac{1}{2}\left(h(u)-h\left(u+(1-s)q\right)\right):0<u,\, u+(1-s)q<q\right\} \,.\]
This last fact is exactly what is expressed by the formula (\ref{eq:preffnn}),
and so completes the proof. 
\end{rem}
We conclude this section by mentioning two more results, consequences
of the formula (\ref{eq:ffnn}), which each apply in the {}``limiting''
case $s=1/2$ to any given admissible $Q$ and to each measurable
real function $f$ defined on $Q$. We will not actually need to use
them further in the current version of this paper. The first of these
is the inequality\begin{eqnarray*}
 &  & \mathbf{J}(f,Q,1/2)\\
 & \ge & \frac{1}{2}\min\left\{ \left(f\chi_{Q}\right)^{*}\left(\frac{\lambda(Q)}{4}\right)-\left(f\chi_{Q}\right)^{*}\left(\frac{\lambda(Q)}{2}\right),\left(f\chi_{Q}\right)^{*}\left(\frac{\lambda(Q)}{2}\right)-\left(f\chi_{Q}\right)^{*}\left(\frac{3\lambda(Q)}{4}\right)\right\} \end{eqnarray*}
which follows from (\ref{eq:ffnn}) combined with the simple observation
that any closed subinterval $I$ of $(0,\lambda(Q))$ of length $\lambda(Q)/2$
must contain at least one of the two closed intervals $\left[\frac{\lambda(Q)}{4},\frac{\lambda(Q)}{2}\right]$
and $\left[\frac{\lambda(Q)}{2},\frac{3\lambda(Q)}{4}\right]$. For
our second result we also note that, since the above interval $I$
is closed, it must also contain the point $\lambda(Q)/2$ in its interior.
Consequently (\ref{eq:ffnn}) also gives us that \[
\mathbf{J}(f,Q,1/2)\ge\frac{1}{2}\left(\lim_{t\nearrow\lambda(Q)/2}\left(f\chi_{Q}\right)^{*}(t)-\left(f\chi_{Q}\right)^{*}\left(\frac{\lambda(Q)}{2}\right)\right)\,.\]

\section{\label{sec:FnOfOneVariable}non increasing functions of one variable
in BMO. Some simple calculations.}

Suppose that $d=1$, that $D$ is a bounded open interval, and that
$f:D\to\mathbb{R}$ is non increasing and right continuous. Our main
aim in this section is to prove that a slight variant of the John-Strömberg
inequality (namely (\ref{eq:odms})) holds for these very special
choices of $D$ and $f$.

The proof of (\ref{eq:odms}) in this special case is of course much
simpler than any known proofs of the John-Nirenberg or John-Strömberg
inequalities for the general case. But the results of other sections
will enable us to deduce the general case from this special case,
albeit with not particularly good constants, and with restrictions
on the range of the parameter $s$ appearing in the John-Strömberg
functional.

We obtain (\ref{eq:odms}) as a consequence of the following two lemmata.
The first of these bounds the functional $\left\Vert f\right\Vert _{BMO\left(I,\mathcal{Q}(I)\right)}^{(\mathbf{J},s)}$
by another functional which has been found to be useful in various
contexts and is more or less connected to the functional $\sup_{t>0}f^{**}(t)-f^{*}(t)$
which was introduced in \cite{bds}. Other results about these and
similar functionals can be found, for example, in \cite{BagbyKurtz1,BagbyKurtz2,lerner98,YoramPasha02}
and in a large number of subsequent papers.

For simplicity, we only consider (and in fact only need to consider)
the interval $I=(0,1)$ at this stage. 
\begin{lem}
\label{lem:bromp}Suppose that $s\in(0,1/2)$ and $\rho=\frac{s}{1-s}$.
Suppose that $f:(0,1)\to\mathbb{R}$ is a non increasing right continuous
function. Then \[
\sup_{t\in(0,1/2]}\left(f(\rho t)-f(t)\right)\le2\left\Vert f\right\Vert _{BMO\left((0,1),\mathcal{Q}((0,1))\right)}^{(\mathbf{J},s)}\,.\]
\end{lem}
\begin{rem}
Lerner \cite[Theorem  3.1, p. 52]{lerner98} has obtained a much more
general result with a much more elaborate proof, which essentially
implies this lemma. 
\end{rem}
\noindent \textit{Proof.} The properties of $f$ permit us to use
the formula (\ref{eq:preffnn}) of Proposition \ref{pro:prop}. Let
$(a,b)$ be an arbitrary open subinterval of $(0,1)$. Via an obvious
change of variables (translation, e.g.~apply (\ref{eq:preaffine})
with $r=1$ and $x_{0}=a$) the formula (\ref{eq:preffnn}) tells
us that \begin{equation}
\mathbf{J}(f,(a,b),s)=\frac{1}{2}\inf\left\{ f(a+u)-f\left(a+u+(1-s)(b-a)\right):0<u<s(b-a)\right\} \,.\label{eq:upq}\end{equation}

Let $[c,d]$ be an arbitrary closed subinterval of $(a,b)$ of length
$(1-s)(b-a)$. Then $d$ must satisfy \[
d>a+(1-s)(b-a)\]
and $c$ must satisfy \[
c<b-\left(1-s\right)(b-a)\,.\]
From these estimates it follows that \[
f(d)\le f(a+(1-s)(b-a))\mbox{ and }f(b-(1-s)(b-a))\le f(c)\,.\]
These estimates imply that \[
f\left(b-\left(1-s\right)(b-a)\right)-f\left(a+(1-s)(b-a)\right)\le f(c)-f(d)\,.\]
Taking the infimum over all subintervals $[c,d]$ of $\left(a,b\right)$
which have length $(1-s)(b-a)$ and applying (\ref{eq:upq}), we see
that \begin{equation}
f\left(b-\left(1-s\right)(b-a)\right)-f\left(a+(1-s)(b-a)\right)\le2\mathbf{J}(f,(a,b),s)\le2\left\Vert f\right\Vert _{BMO\left(I,\mathcal{Q}(I)\right)}^{(\mathbf{J},s)}\label{eq:tds}\end{equation}
 whenever $0\le a<b\le1$.

In particular, for an arbitrary $t\in(0,1/2]$, let us choose $a=0$
and $b=\frac{t}{1-s}$. Since $s\in(0,1/2)$ we have $b\in(t,2t)\subset(t,1)$.
For these choices of $a$ and $b$, the left hand side of (\ref{eq:tds})
equals $f(\rho t)-f(t)$. So the proof of the lemma is complete. $\qed$

\bigskip{}

Our second lemma enables us to bound the size of our function $f$
by an expression depending on the functional $\sup_{t\in(0,1/2]}\left(f(\rho t)-f(t\right)$
and consequently to obtain an inequality which is quite close to the
one that we need.
\begin{lem}
\label{lem:bonk}The inequality \begin{equation}
f(u)-f(v)\le\left(1+\frac{\log\frac{v}{u}}{\log(1/\rho)}\right)\sup_{t\in(0,1/2]}\left(f(\rho t)-f(t)\right)\label{eq:burp}\end{equation}
 holds for every non increasing function $f:(0,1)\to\mathbb{R}$\,,
every $\rho\in(0,1)$, and every $u$ and $v$ satisfying $0<u<v\le1/2$\,.
\end{lem}
As an immediate consequence we obtain 
\begin{cor}
If $f$ and $\rho$ are as in the preceding lemma and if \textup{\[
\sup_{t\in(0,1/2]}\left(f(\rho t)-f(t)\right)\le c\,,\]
}then\begin{equation}
\left|\left\{ t\in(0,1):f(t)-f(1/2)\ge\alpha\right\} \right|\le\frac{1}{2\rho}\exp\left(-\frac{\alpha\log(1/\rho)}{c}\right)\mbox{ for each }\alpha\ge0\,.\label{eq:ejn}\end{equation}

\end{cor}
\noindent \textit{Proof of the lemma and its corollary.} Let $N$
be the unique positive integer for which $\rho^{N}v\le u<\rho^{N-1}v$.
Then $(1/\rho)^{N-1}<\frac{v}{u}$ and so $N<1+\frac{\log\frac{v}{u}}{\log(1/\rho)}$\,.
Hence \begin{eqnarray*}
f(u)-f(v) & \le & f\left(\rho^{N}v\right)-f\left(\rho^{0}v\right)=\sum_{n=1}^{N}\left(f\left(\rho^{n}v\right)-f\left(\rho^{n-1}v\right)\right)\\
 & \le & N\sup_{t\in(0,1/2]}\left(f(\rho t)-f(t)\right)\,.\end{eqnarray*}
 This, combined with our estimate for $N$, establishes (\ref{eq:burp}).

Now let us prove (\ref{eq:ejn}) under the stated hypothesis. Setting
$v=1/2$ in (\ref{eq:burp}) gives us that \begin{equation}
f(u)-f(1/2)\le c\left(1-\frac{\log2u}{\log(1/\rho)}\right)\mbox{ for all }u\in(0,1/2)\,.\label{eq:tsjdt}\end{equation}
For each $\alpha\ge0$, the set $\left\{ t\in(0,1):f(t)-f(1/2)\ge\alpha\right\} $
is of course an interval contained in $(0,1/2]$. It follows from
(\ref{eq:tsjdt}) that the length of this interval cannot exceed $\frac{1}{2\rho}\exp\left(-\frac{\alpha\log(1/\rho)}{c}\right)$.
$\qed$

The preceding two lemmata and corollary have the following immediate
consequence.

Let $f:I\to\mathbb{R}$ be a non increasing right continuous function
on the interval $I=(0,1)$. Then, for each $s\in(0,1/2]$,

\begin{equation}
\left|\left\{ t\in I:f(t)-f\left(\frac{1}{2}\right)\ge\alpha\right\} \right|\le\frac{1-s}{2s}\cdot\exp\left(-\frac{\alpha\log\left(\frac{1}{s}-1\right)}{2\left\Vert f\right\Vert _{BMO\left(I,\mathcal{Q}(I)\right)}^{(\mathbf{J},s)}}\right)\mbox{ for all }\alpha\ge0\,.\label{eq:jtddms}\end{equation}
Note that here we can also permit $s$ to take the limiting value
$s=1/2$ and we can permit $\left\Vert f\right\Vert _{BMO\left(I,\mathcal{Q}(I)\right)}^{(\mathbf{J},s)}$
to be infinite (provided we agree to interpret both $1/\infty$ and
$0/\infty$ as $0$). In such cases the right hand side of (\ref{eq:jtddms})
is greater than or equal to $1/2$ which means that (\ref{eq:jtddms})
is also true, trivially so, in these {}``limiting'' cases.

It will now be a very simple matter to deduce a more general version
of (\ref{eq:jtddms}) for the case where $I$ is an arbitrary open
interval $(a,b)$. Suppose that $f:(a,b)\to\mathbb{R}$ is right continuous
and non increasing. Define $g:(0,1)\to\mathbb{R}$ by $g(t)=f(a+(b-a)t)$.
Then (cf.~(\ref{eq:preaffine}) and (\ref{eq:affinestuff})) we have
$\left\Vert f\right\Vert _{BMO\left(I,\mathcal{Q}(I)\right)}^{(\mathbf{J},s)}=\left\Vert g\right\Vert _{BMO\left((0,1),\mathcal{Q}(0,1)\right)}^{(\mathbf{J},s)}$.
Furthermore, the set \[
E_{1}=\left\{ t\in(a,b):f(t)-f\left(\frac{a+b}{2}\right)\ge\alpha\right\} \]
coincides with the set $(b-a)E_{2}+a$ where \[
E_{2}=\left\{ x\in(0,1):g(x)-g(1/2)\ge\alpha\right\} \,.\]
Therefore, applying (\ref{eq:jtddms}) to the function $g$ and multiplying
both sides of the resulting inequality by $b-a=\left|I\right|$ gives
us that \begin{equation}
\left|\left\{ t\in I:f(t)-f\left(c_{I}\right)\ge\alpha\right\} \right|\le\frac{1-s}{2s}\cdot\left|I\right|\cdot\exp\left(-\frac{\alpha\log\left(\frac{1}{s}-1\right)}{2\left\Vert f\right\Vert _{BMO\left(I,\mathcal{Q}(I)\right)}^{(\mathbf{J},s)}}\right)\label{eq:odms}\end{equation}
 for every $\alpha\ge0$ and every $s\in(0,1/2]$ and for every open
interval $I$, where $c_{I}$ denotes the midpoint of $I$.

This is the inequality that we need to apply in the proof of our main
result, Theorem \ref{thm:mainjs} of Section \ref{sec:Putting}.

One immediate consequence of (\ref{eq:odms}) together with the inequality
(\ref{eq:talpt}) recalled in Remark \ref{rem:Cheb}, is that \begin{equation}
\left|\left\{ t\in I:f(t)-f\left(c_{I}\right)\ge\alpha\right\} \right|\le\frac{1-s}{2s}\cdot\left|I\right|\cdot\exp\left(-\frac{\alpha s\log\left(\frac{1}{s}-1\right)}{2\left\Vert f\right\Vert _{BMO\left(I,\mathcal{Q}(I)\right)}}\right)\label{eq:brerp}\end{equation}

\begin{rem}
\label{rem:transit}The inequality (\ref{eq:odms}) is not exactly
the John-Strömberg inequality, and (\ref{eq:brerp}) is not exactly
the John-Nirenberg inequality for the non increasing function $f$
on the interval $I$, since the John-Str\"omberg and John-Nirenberg
inequalities are for the measure of the sets $\left\{ t\in I:\left|f(t)-f\left(c_{I}\right)\right|>\alpha\right\} $
and $\left\{ t\in I:\left|f(t)-f_{I}\right|>\alpha\right\} $ respectively,
where $f_{I}=\frac{1}{\left|I\right|}\int_{I}f(t)dt$. But a simple
argument (cf.~the proof of Lemma \ref{lem:temjemt}) enables us to
replace $f(x)-f(c_{I})$ by $\left|f(x)-f(c_{I})\right|$ in the left
hand side of (\ref{eq:odms}), provided that we also multipy the right
hand side by $2$. The possibility of replacing $f(c_{I})$ (which
is of course the median of $f$ on $I$) by $f_{I}$ has already been
discussed in Remark \ref{rem:difver}. This replacement can be made,
again at the price of increasing the constant on the right hand side.
(More details about doing this are given in Lemma \ref{lem:av2med}
of Appendix \ref{sub:dfojne}.) We recall that a version of the John-Nirenberg
inequality is given in \cite{koren1990} in which it is shown, for
$d=1$, that the optimal value of the constant $b$ in (\ref{eq:jne})
is $2/e$. See also \cite{korenovskii} p.\ 77. It is not difficult
to check that the inequality obtained from (\ref{eq:brerp}) by applying
these simple steps does not have this optimal value for $b$. 
\end{rem}

\section{\label{sec:simplifyJohnStromberg}A reduction of the proof of the
John-Strömberg Theorem to a special case.}

Having, in the previous section, prepared the auxiliary results that
we need about special functions of one variable, we now turn to consider
functions of several variables.

In our (very slightly) different notation, Lemma 3.1 on p.~517 of
\cite{stromberg} states that

\begin{equation}
\begin{array}{c}
{\displaystyle \lambda\left(\left\{ x\in Q:\left|f(x)-m_{f}(Q)\right|>\alpha\right\} \right)\le C\lambda(Q)\exp\left(-\frac{c\alpha}{\left\Vert f\right\Vert _{BMO_{0,s}}}\right)}\end{array}\label{eq:threeone}\end{equation}
for all $\alpha\ge0$ and $s\in(0,1/2]$. 

Here $C$ and $c$ are positive constants depending only on $d$ and
$m_{f}(Q)$ is a median of $f$ on the arbitrary cube $Q$ in $\mathbb{R}^{d}$.

(Here we are again using the notation specified in (\ref{eq:pone}).
Note that there is a small misprint in \cite{stromberg}, namely the
factor $\lambda(Q)$ (or $\left|Q\right|$) has been omitted there.)

Our main goal in this section is to show that, in order to prove the
inequality (\ref{eq:threeone}) for the specified values of $\alpha$
and $s$, and some other inequalities like it, it suffices to obtain
such inequalities, but with different values of the constants $c$
and $C$, in the special case where $f$ is a non negative function
taking only integer values. This fact will be precisely formulated
as Theorem \ref{thm:temjemt-1}. (The question of whether such inequalities
actually do hold in that special case will be deferred to Section
\ref{sec:Putting}. We will be able to answer it there, with the help
of results from other sections.)

Here we can just as easily work in the rather more general context
of the space $BMO(D,\mathcal{E})$ of Definition \ref{def:bmode}.
Indeed, doing so will be convenient, since we will later want to apply
the result of this section in such a general context, which will include,
for example, the particular case of special rectangles (\ref{eq:wik0})
as well as the case of usual cubes (\ref{eq:jnc}). Thus, throughout
this section $D$ will denote some arbitrary but fixed measurable
subset of $\mathbb{R}^{d}$ and $\mathcal{E}$ will denote some arbitrary
but fixed collection of admissible subsets of $D$.

Our first (easy) step is to reduce everything to the case of non negative
functions.
\begin{lem}
\label{lem:temjemt}Let $E$ be a fixed admissible set in $\mathcal{E}$,
let $s$ be a fixed number in $(0,1/2]$, and let $c$ and $C$ be
positive constants. Suppose that every non negative measurable function
$f:D\to[0,\infty)$ satisfies the inequality

\begin{equation}
\lambda\left(\left\{ x\in E:f(x)>\alpha\right\} \right)\le C\lambda(E)\exp\left(-\frac{c\alpha}{\left\Vert f\right\Vert _{BMO(D,\mathcal{E})}^{(\mathbf{J},s)}}\right)\mbox{ for all }\alpha\ge0\,.\label{eq:pozthreeone}\end{equation}
 Then every measurable function $f:D\to\mathbb{R}$ satisfies

\begin{equation}
\lambda\left(\left\{ x\in E:\left|f(x)-m\right|>\alpha\right\} \right)\le2C\lambda(E)\exp\left(-\frac{c\alpha}{\left\Vert f\right\Vert _{BMO(D,\mathcal{E})}^{(\mathbf{J},s)}}\right)\mbox{ for all }\alpha\ge0\label{eq:newthreeone}\end{equation}
 whenever $m$ is a median of $f$ on $E$.

In fact this same conclusion also holds under weaker hypotheses, namely
if it is only known that (\ref{eq:pozthreeone}) holds for those non
negative functions $f:D\to[0,\infty)$ having the additional property
that \begin{equation}
\lambda\left(\left\{ x\in E:f(x)>0\right\} \right)\le\frac{1}{2}\lambda(E)\,.\label{eq:addprop}\end{equation}
\end{lem}
\begin{rem}
Our {}``natural'' applications of Lemma \ref{lem:temjemt}, will
be in the case where the collection $\mathcal{E}$ includes the set
$D$ itself, and we choose $E=D$. 
\end{rem}
\noindent \textit{Proof}. We have to prove (\ref{eq:newthreeone})
for an arbitrary measurable function $f:D\to\mathbb{R}$ with median
$m$ on $E$\,. Let $g=f-m$. Obviously $\mathbf{J}(g,E,s)=\mathbf{J}(f,E,s)$
and $\left\Vert g\right\Vert _{BMO(D,\mathcal{E})}^{(\mathbf{J},s)}=\left\Vert f\right\Vert _{BMO(D,\mathcal{E})}^{(\mathbf{J},s)}$
and so it will suffice to prove that\begin{equation}
\lambda\left(\left\{ x\in E:\left|g(x)\right|>\alpha\right\} \right)\le2C\lambda(E)\exp\left(-\frac{c\alpha}{\left\Vert g\right\Vert _{BMO(D,\mathcal{E})}^{(\mathbf{J},s)}}\right)\mbox{ for all }\alpha\ge0\,.\label{eq:gthreeone}\end{equation}
 The left hand side of (\ref{eq:gthreeone}) equals \begin{eqnarray*}
 &  & \lambda\left(\left\{ x\in E:g(x)>\alpha\right\} \right)+\lambda\left(\left\{ x\in E:g(x)<-\alpha\right\} \right)\\
 & = & \lambda\left(\left\{ x\in E:g_{+}(x)>\alpha\right\} \right)+\lambda\left(\left\{ x\in E:g_{-}(x)>\alpha\right\} \right)\end{eqnarray*}
 where, as usual, $g_{+}=\max\left\{ g,0\right\} $ and $g_{-}=g_{+}-g=\max\left\{ -g,0\}\right\} $.
We can apply Lemma \ref{lem:JLIP} with $\varphi(t)=\max\left\{ t,0\right\} $
to obtain that \[
\left\Vert g_{+}\right\Vert _{BMO(D,\mathcal{E})}^{(\mathbf{J},s)}\le\left\Vert g\right\Vert _{BMO(D,\mathcal{E})}^{(\mathbf{J},s)}\mbox{ and }\left\Vert g_{-}\right\Vert _{BMO(D,\mathcal{E})}^{(\mathbf{J},s)}\le\left\Vert -g\right\Vert _{BMO(D,\mathcal{E})}^{(\mathbf{J},s)}.\]
Obviously $\left\Vert -g\right\Vert _{BMO(D,\mathcal{E})}^{(\mathbf{J},s)}=\left\Vert g\right\Vert _{BMO(D,\mathcal{E})}^{(\mathbf{J},s)}$\,.
Thus, if we apply (\ref{eq:pozthreeone}) to each of the non negative
functions $g_{+}$ and $g_{-}$ and sum the results, we obtain (\ref{eq:gthreeone})
and therefore that $f$ indeed satisfies (\ref{eq:newthreeone}).

It remains to justify the claim in the last sentence of the statement
of the lemma. Since $m$ is a median of $f$, it follows that $0$
is a median of $g$ and so the two functions $g_{+}$ and $g_{-}$
to which we have applied (\ref{eq:pozthreeone}) satisfy \[
\lambda\left(\left\{ x\in E:g_{+}(x)>0\right\} \right)\le\frac{1}{2}\lambda(E)\mbox{ and }\lambda\left(\left\{ x\in E:g_{-}(x)>0\right\} \right)\le\frac{1}{2}\lambda(E)\,.\]
$\qed$

Our next step is to reduce the proof of (\ref{eq:pozthreeone}) to
the case of appropriate integer valued functions.
\begin{lem}
\label{lem:intval}Let $E$, $\mathcal{E}$, $s$, $c$ and $C$ be
as in the statement of Lemma \ref{lem:temjemt}. Suppose that (\ref{eq:pozthreeone})
holds for every non negative measurable function $f:D\to[0,\infty)$
which takes only integer values and satisfies $\left\Vert f\right\Vert _{BMO(D,\mathcal{E})}^{(\mathbf{J},s)}\le1/2$
and (\ref{eq:addprop}). Then (\ref{eq:pozthreeone}) also holds for
every measurable $f:D\to[0,\infty)$, which satisfies (\ref{eq:addprop}),
but with the constants $C$ and $c$ replaced by $c_{1}=c/4$ and
$C_{1}=\max\left\{ C,e^{c}\right\} $.\end{lem}
\begin{rem}
It will be clear from the following proof that we can also obtain
the following additional result: Suppose that in the above lemma we
impose the stronger condition that (\ref{eq:pozthreeone}) holds also
for every non negative measurable function $f:D\to[0,\infty)$ which
takes only integer values and satisfies $\left\Vert f\right\Vert _{BMO(D,\mathcal{E})}^{(\mathbf{J},s)}\le1/2$
but does not necessarily satisfy (\ref{eq:addprop}). Then we obtain
the stronger conclusion that (\ref{eq:pozthreeone}) holds for every
measurable $f:D\to[0,\infty)$, but with the constants $C$ and $c$
replaced by $c_{1}=c/4$ and $C_{1}=\max\left\{ C,e^{c}\right\} $.
I.e., in this case we can obtain (\ref{eq:pozthreeone}) also for
functions $f:D\to[0,\infty)$ which do not satisfy (\ref{eq:addprop}).
\end{rem}
\noindent \textit{Proof.} We shall use the function $\varphi:[0,\infty)\to[0,\infty)$
which is defined by \[
\varphi(t)=\left\{ \begin{array}{ccc}
0 & , & 0\le t\le1/2_{\phantom{Q}}\\
n & , & n-1/2<t\le n+1/2\mbox{ for each }n\in\mathbb{N^{\phantom{T}}}\,.\end{array}\right.\]
 We will need the following three obvious or easily verified properties
of $\varphi$:

\[
\varphi([0,\infty))=\mathbb{\mathbb{N}}\cup\left\{ 0\right\} \,,\]

\begin{equation}
\varphi(t)-\varphi(s)\in\left\{ 0,1\right\} \mbox{ whenever }0\le s\le t\le s+1/2\,\label{eq:phat}\end{equation}
 and \begin{equation}
\varphi(t)\ge t-1/2\mbox{ for all }t\ge0\,.\label{eq:gerp}\end{equation}

Suppose that $f:D\to[0,\infty)$ is an arbitrary measurable function
which satisfies (\ref{eq:addprop}) and also \begin{equation}
0<\left\Vert f\right\Vert _{BMO(D,\mathcal{E})}^{(\mathbf{J},s)}<\frac{1}{3}\,.\label{eq:ner}\end{equation}
 The composed function $\varphi\circ f$ also satisfies (\ref{eq:addprop})
since \begin{eqnarray*}
\lambda\left(\left\{ x\in E:\varphi\circ f(x)>0\right\} \right) & = & \lambda\left(\left\{ x\in E:\varphi\circ f(x)\ge1\right\} \right)\\
 & = & \lambda\left(\left\{ x\in E:f(x)>1/2\right\} \right)\le\frac{1}{2}\lambda(E)\,.\end{eqnarray*}
 We will next show that, furthermore, $\varphi\circ f$ satisfies
\begin{equation}
\left\Vert \varphi\circ f\right\Vert _{BMO(D,\mathcal{E})}^{(\mathbf{J},s)}\le1/2\,.\label{eq:wih}\end{equation}
 Let $W$ be an arbitrary set in $\mathcal{E}$\,. Then (\ref{eq:ner})
implies that $\mathbf{J}(f,W,s)<1/3$\,. Therefore (cf.~(\ref{eq:altdef}))
there exists some $\alpha\in[0,1/3)$ and some $\gamma\in\mathbb{R}$
such that \begin{equation}
\lambda\left(\left\{ x\in W:\left|f(x)-\gamma\right|\le\alpha\right\} \right)>(1-s)\lambda(W)\,.\label{eq:ovmt}\end{equation}
 Let us choose \[
\gamma_{1}=\frac{1}{2}\left(\varphi(\gamma-\alpha)+\varphi(\gamma+\alpha)\right)\mbox{ and }\alpha_{1}=\frac{1}{2}\left(\varphi(\gamma+\alpha)-\varphi(\gamma-\alpha)\right)\,.\]
Since $0\le\alpha<1/3$ we obtain from (\ref{eq:phat}) that $\alpha_{1}$
is either $0$ or $1/2$. Since $\varphi$ is non decreasing, we also
obtain that \begin{eqnarray*}
\left\{ x\in W:\left|f(x)-\gamma\right|\le\alpha\right\}  & = & \left\{ x\in W:\gamma-\alpha\le f(x)\le\gamma+\alpha\right\} \\
 & \subset & \left\{ x\in W:\varphi(\gamma-\alpha)\le\varphi\circ f(x)\le\varphi(\gamma+\alpha)\right\} \\
 & = & \left\{ x\in W:\gamma_{1}-\alpha_{1}\le\varphi\circ f(x)\le\gamma_{1}+\alpha_{1}\right\} \\
 & = & \left\{ x\in W:\left|\varphi\circ f(x)-\gamma_{1}\right|\le\alpha_{1}\right\} \,.\end{eqnarray*}
 Thus we deduce, using (\ref{eq:ovmt}) and (\ref{eq:altdef}) once
more, that \[
\mathbf{J}\left(\varphi\circ f,W,s\right)\le\alpha_{1}\le1/2\]
 and this establishes (\ref{eq:wih}).

Since the function $\varphi\circ f$ is also non negative and integer
valued, we have, according to the hypotheses of the lemma, that\begin{equation}
\lambda\left(\left\{ x\in E:\varphi\circ f(x)>\alpha\right\} \right)\le C\lambda(E)\exp\left(-\frac{c\alpha}{\left\Vert \varphi\circ f\right\Vert _{BMO(D,\mathcal{E})}^{(\mathbf{J},s)}}\right)\mbox{ for all }\alpha\ge0\,.\label{eq:dddiy}\end{equation}

The inequality (\ref{eq:gerp}) implies that \[
\left\{ x\in E:f(x)>\alpha\right\} \subset\left\{ x\in E:\varphi\circ f(x)>\alpha-1/2\right\} \,,\]
 and consequently, for all $\alpha>1/2$, it follows, using (\ref{eq:dddiy})
and then (\ref{eq:wih}) and then (\ref{eq:ner}), that\begin{eqnarray*}
\lambda\left(\left\{ x\in E:f(x)>\alpha\right\} \right) & \le & C\lambda(E)\exp\left(-\frac{c(\alpha-1/2)}{\left\Vert \varphi\circ f\right\Vert _{BMO(D,\mathcal{E})}^{(\mathbf{J},s)}}\right)\\
 & \le & C\lambda(E)\exp\left(-2c(\alpha-1/2)\right)\,.\end{eqnarray*}
If we now restrict $\alpha$ to the range $\alpha\ge1$, we also have
$\alpha-1/2\ge\alpha/2$ and so \begin{eqnarray*}
\lambda\left(\left\{ x\in E:f(x)>\alpha\right\} \right) & \le & C\lambda(E)\exp\left(-c\alpha\right)\le C_{1}\lambda(E)\exp\left(-c\alpha\right)\,,\end{eqnarray*}
recalling, as stated in the lemma, that $C_{1}=\max\left\{ C,e^{c}\right\} $.

Now let us consider the range of values $0\le\alpha<1$. Of course
\[
\lambda\left(\left\{ x\in E:f(x)>\alpha\right\} \right)\le\lambda(E)\]
for these (and all other) values of $\alpha$. Furthermore, for each
$\alpha\in[0,1)$ we have $C_{1}e^{-c\alpha}\ge1$ and therefore \[
\lambda\left(\left\{ x\in E:f(x)>\alpha\right\} \right)\le C_{1}\lambda(E)\exp\left(-c\alpha\right)\,.\]
 We have thus now shown that, subject to the hypotheses of the lemma,
the inequality \begin{equation}
\lambda\left(\left\{ x\in E:f(x)>\beta\right\} \right)\le C_{1}\lambda(E)\exp\left(-c\beta\right)\label{eq:kromp}\end{equation}
holds for all $\beta\ge0$ and for all those measurable functions
$f:D\to[0,\infty)$ which satisfy \begin{equation}
0<\left\Vert f\right\Vert _{BMO(D,\mathcal{E})}^{(\mathbf{J},s)}<1/3\label{eq:slish}\end{equation}
and (\ref{eq:addprop}). Now we can easily obtain the required inequality
(\ref{eq:pozthreeone}) without having to impose (\ref{eq:slish}).
Given an arbitrary measurable function $f:D\to[0,\infty)$ satisfying
(\ref{eq:addprop}) and any $\alpha>0$, we let $\widetilde{f}=f/4\left\Vert f\right\Vert _{BMO(D,\mathcal{E})}^{(\mathbf{J},s)}$
and choose $\beta=\alpha/4\left\Vert f\right\Vert _{BMO(D,\mathcal{E})}^{(\mathbf{J},s)}$.
Then by homogeneity (cf.~(\ref{eq:homog})), we have $\left\Vert \widetilde{f}\right\Vert _{BMO(D,\mathcal{E})}^{(\mathbf{J},s)}=1/4$
and of course $\widetilde{f}$ also satisfies (\ref{eq:addprop}).
So we can apply (\ref{eq:kromp}) to $\widetilde{f}$ and obtain (\ref{eq:pozthreeone}),
completing the proof of the lemma. $\qed$

We can summarize the results of this section by the following theorem,
whose proof follows immediately from the previous two lemmata.
\begin{thm}
\label{thm:temjemt-1}Let $E$ be a fixed admissible set in $\mathcal{E}$,
let $s$ be a fixed number in $(0,1/2]$, and let $c$ and $C$ be
positive constants. Let $\Phi$ be the set of all non negative measurable
functions $f:D\to[0,\infty)$ which take only integer values, and
which satisfy $\left\Vert f\right\Vert _{BMO(D,\mathcal{E})}^{(\mathbf{J},s)}\le1/2$
and $\lambda\left(\left\{ x\in E:f(x)>0\right\} \right)\le\frac{1}{2}\lambda(E)$.
Suppose that every $f\in\Phi$ satisfies \begin{equation}
\lambda\left(\left\{ x\in E:f(x)>\alpha\right\} \right)\le C\lambda(E)\exp\left(-\frac{c\alpha}{\left\Vert f\right\Vert _{BMO(D,\mathcal{E})}^{(\mathbf{J},s)}}\right)\mbox{ for all }\alpha\ge0\,.\label{eq:brinz}\end{equation}
 Then every measurable function $f:D\to\mathbb{R}$ satisfies\begin{equation}
\lambda\left(\left\{ x\in E:\left|f(x)-m\right|>\alpha\right\} \right)\le2\max\left\{ C,e^{c}\right\} \lambda(E)\exp\left(-\frac{c\alpha}{4\left\Vert f\right\Vert _{BMO(D,\mathcal{E})}^{(\mathbf{J},s)}}\right)\label{eq:inxz}\end{equation}
 for all $\alpha\ge0$, whenever $m$ is a median of $f$ on $E$.
\end{thm}

\section{The {}``geometrical'' component of our proof}

\subsection{A {}``balancing act'' between two subsets of a cube. The {}``bi-density''
constant.}

We begin by stating a simple result which is a sort of {}``prototype''
of the main result that we seek in this section. It will also be a
tool for proving that main result.
\begin{lem}
\label{lem:prototype}Let $Q$ be a cube in $\mathbb{R}^{d}$ and
let $E$ be a measurable subset of $Q$ such that $0<\lambda(E)<\lambda(Q)$.
Then there exists a cube $W$ contained in $Q$ such that \begin{equation}
\lambda(W\setminus E)=\lambda(W\cap E)=\frac{1}{2}\lambda(W)\,.\label{eq:zprop}\end{equation}

\end{lem}
We can just as easily prove a slightly more general result. See Lemma
\ref{lem:prehlfm} below. Setting $s=1/2$ in that lemma will give
the result just stated here.

If the cube $Q$ in Lemma \ref{lem:prototype} is dyadic, and we want
to only consider subcubes $W$ which are also dyadic, then (as is
made explicit below in Remark \ref{rem:bestdyadic}) we cannot hope
in general to obtain one of them which satisfies (\ref{eq:zprop}).
Instead, as the next lemma states, we can obtain a dyadic subcube
$W$ satisfying a rather weaker property.
\begin{lem}
\label{lem:prototype-1}Let $Q$ be a dyadic cube in $\mathbb{R}^{d}$
and let $E$ be a measurable subset of $Q$ such that $0<\lambda(E)<\lambda(Q)$.
Then there exists a dyadic cube $W$ contained in $Q$ such that \begin{equation}
\min\left\{ \lambda(W\setminus E),\lambda(W\cap E)\right\} \ge2^{-d}(1-2^{-d})\lambda(W)\,.\label{eq:wprop}\end{equation}

\end{lem}
\noindent \textit{Proof.} Clearly any measurable set $W$ (whether
or not it is a dyadic cube) satisfies the condition \begin{equation}
2^{-d}\lambda(W)\le\lambda(W\setminus E)\le(1-2^{-d})\lambda(W)\label{eq:mxm8}\end{equation}
 if and only if it satisfies \begin{equation}
2^{-d}\lambda(W)\le\lambda(W\cap E)\le(1-2^{-d})\lambda(W)\,.\label{eq:btr}\end{equation}
 Any cube $W$ satisfying these conditions of course also satisfies
(\ref{eq:wprop}). Thus the only case which we have to consider is
when (\ref{eq:mxm8}) and (\ref{eq:btr}) are not satisfied for any
dyadic subcube $W$ of $Q$. Then, in particular, $Q$ itself must
satisfy either \begin{equation}
\lambda(Q\cap E)<2^{-d}\lambda(Q)\label{eq:aak}\end{equation}
 or $\lambda(Q\cap E)>(1-2^{-d})\lambda(Q)$. We can and will suppose
that the former of these two conditions holds. If the latter condition
holds then the proof is exactly the same, except for an interchange
of the roles of the sets $E$ and $F$ where $F=Q\setminus E$.

Let $x$ be a density point of $E$ in the interior of $Q$. By an
appropriate version of the Lebesgue differentiation theorem, (see
e.g.\ \cite{WheedenZygmund} Theorem 7.16 pp.~108--109) there exists
a dyadic cube $U$ such that $x\in U\subset Q$ and \begin{equation}
\lambda(U\cap E)\ge(1-2^{-d})\lambda(U)\,.\label{eq:aakz}\end{equation}
 Consider the finite sequence of dyadic cubes $\left\{ U_{n}\right\} _{n=1}^{N}$
for which $U_{1}=U$ and $U_{N}=Q$ and which satisfy $U_{n}\subset U_{n+1}$
and $\lambda(U_{n+1})=2^{d}\lambda(U_{n})$ for each $n=1,2,...,N-1$.
In view of (\ref{eq:aak}) and (\ref{eq:aakz}) the inequality \begin{equation}
\lambda(U_{n}\cap E)\ge(1-2^{-d})\lambda(U_{n})\label{eq:imk}\end{equation}
 holds for $n=1$ but not for $n=N$. So there must exist some integer
$n$ with $1\le n\le N-1$ for which (\ref{eq:imk}) holds and also
\begin{equation}
\lambda(U_{n+1}\cap E)<(1-2^{-d})\lambda(U_{n+1})\,.\label{eq:ronz}\end{equation}
 Since we have excluded the possibility that $U_{n+1}$ satisfies
(\ref{eq:btr}), it follows from (\ref{eq:ronz}) that $\lambda(U_{n+1}\cap E)<2^{-d}\lambda(U_{n+1})$,
or, equivalently, that $\lambda(U_{n+1}\setminus E)>(1-2^{-d})\lambda(U_{n+1})$.
Using (\ref{eq:imk}), we see that \[
\lambda(U_{n+1}\cap E)\ge\lambda(U_{n}\cap E)\ge2^{-d}(1-2^{-d})\lambda(U_{n+1})\]
 and so the dyadic cube $W=U_{n+1}$ satisfies (\ref{eq:wprop}).
$\qed$ 
\begin{rem}
\label{rem:bestdyadic}We have not bothered to check whether it is
possible to obtain a stronger conclusion in Lemma \ref{lem:prototype-1}
where the constant $2^{-d}(1-2^{-d})$ in (\ref{eq:wprop}) is replaced
by some larger constant. However any such improvement would not be
very significant, since the simple example where $E$ is a dyadic
subcube of $Q$ shows that the constant in (\ref{eq:wprop}) cannot
exceed $2^{-d}$. 
\end{rem}

Referring back to the terminology introduced in (\ref{eq:jnc}), (\ref{eq:dyadic}),
(\ref{eq:balls}), (\ref{eq:fedja}) and (\ref{eq:wik0}) we can see
that Lemma \ref{lem:prototype} gives us information about the collection
$\mathcal{Q}(D)$ of subcubes of $D$, and Lemma \ref{lem:prototype-1}
gives us analogous information about the collection $\mathcal{D}(D)$
of dyadic subcubes of $D$. The following notion will put the results
of these two lemmata in a more general context. 
\begin{defn}
Let $\mathcal{E}$ be a collection of admissible subsets of $\mathbb{R}^{d}$.
We say that a number $\delta$ is a \textit{bi-density constant} for
$\mathcal{E}$ if, for each $Q\in\mathcal{E}$ and for each measurable
set $E$ for which $0<\lambda(Q\cap E)<\lambda(Q)$, there exists
some set $W\in\mathcal{E}$ with $W\subset Q$ such that \[
\min\left\{ \lambda(W\cap E),\lambda(W\setminus E)\right\} \ge\delta\lambda(W)\,.\]

\end{defn}
Thus Lemma \ref{lem:prototype} tells us that $\delta=1/2$ is a bi-density
constant for $\mathcal{Q}(D)$. The next lemma will show (when we
substitute $s=1/2$) that this is also the case for the collections
$\mathcal{K}(D)$ and $\mathcal{W}(D)$ (defined above in (\ref{eq:fedja})
and (\ref{eq:wik0})). The value $1/2$ is in fact optimal since,
clearly, any bi-density constant for any collection $\mathcal{E}$
always has to satisfy $\delta\le1/2$. Lemma \ref{lem:prototype-1}
and Remark \ref{rem:bestdyadic} tell us that every bi-density constant
$\delta$ for $\mathcal{D}(D)$ must satisfy $2^{-d}(1-2^{-d})\le\delta\le2^{-d}$.
\begin{lem}
\label{lem:prehlfm}Let $K$ be a convex subset of $\mathbb{R}^{d}$
with non empty interior and let $E$ be a measurable subset of $K$
such that $0<\lambda(E)<\lambda(K)$. Then, given an arbitrary number
$s\in(0,1)$, there exists a cube $W\subset K$ for which \begin{equation}
\frac{\lambda(W\cap E)}{\lambda(W)}=s\,.\label{eq:tpdh}\end{equation}

\end{lem}
\noindent \textit{Proof.} Since the boundary of a convex subset of
$\mathbb{R}^{d}$ always has measure zero (see Appendix \ref{sub:convexboundary}),
we may assume without loss of generality that $E$ is contained in
the interior $K^{\circ}$ of $K$.

By the Lebesgue differentiation theorem, almost every point of the
set $E$ is a point of density. So there exists a subcube $W_{a}$
of $K$ centred at one of these density points for which \[
\frac{1}{\lambda(W_{a})}\int_{W_{a}}\chi_{E}d\lambda>s\,.\]
 Similarly there exists another subcube $W_{b}$ of $K$ centred at
some density point of $F:=K^{\circ}\backslash E$ for which \[
\frac{1}{\lambda(W_{a})}\int_{W_{a}}\chi_{E}d\lambda=1-\frac{1}{\lambda(W_{b})}\int_{W_{b}}\chi_{F}d\lambda<s\,.\]
 Now we define a family $\left\{ W(t)\right\} _{t\in[0,1]}$ of subcubes
of $K$ such that $W(0)=W_{a}$ and $W(1)=W_{b}$ and such that $W(t)$
varies continuously as a function of $t$ on $[0,1]$. It is easy
to see how to do this, but we will nevertheless give an explicit recipe.

Let $c(a)$ and $c(b)$ be the centres of $W_{a}$ and $W_{b}$ respectively,
and let $r(a)$ and $r(b)$ be the side lengths of $W_{a}$ and $W_{b}$
respectively. For convenience, and without of loss of generality,
we may suppose that $r(a)\le r(b)$.

If $r(a)=r(b)$ then we simply define $W(t)$ to be the cube of sidelength
$r(a)$ centred at the point $(1-t)c(a)+tc(b)$ for each $t\in[0,1]$.

If $r(a)<r(b)$ then we first define $W(t)$ on $[0,1/2]$ by letting
it be the cube of sidelength $r(a)$ centred at the point $(1-2t)c(a)+2tc(b)$.
Then, for $t\in[1/2,1]$ we let $W(t)$ be the cube centred at $c(b)$
with sidelength $r(a)(2-2t)+r(b)(2t-1)$.

In both of the above cases, the fact that $K$ is convex guarantees
that $W(t)\subset K$ for all $t\in[0,1]$.

It is easy to see that $t\mapsto\frac{1}{\lambda(W(t))}\int_{W(t)}\chi_{E}d\lambda$
is a continuous function of $t$ on $[0,1]$. Therefore it must take
the value $s$ at some point $t_{*}\in(0,1)$. The cube $W=W\left(t_{*}\right)$
satisfies (\ref{eq:tpdh}). $\qed$

\subsection{Some preparations for a {}``balancing act'' between three subsets
of a cube.}

In the previous subsection we considered cubes $Q$ which are the
unions of \textit{two} disjoint sets $E$ and $Q\setminus E$ which
both have positive measure, and we have obtained a subcube $W$ of
$Q$ whose intersections with $E$ and with $Q\setminus E$ are both
{}``significant'' proportions of $W$. Our main goal in the sequel
will be to {}``upgrade'' these kinds of results to a situation where
the cube $Q$ (or a more general set) is the union of \textit{three}
disjoint sets, which we may denote by $E_{+}$ and $E_{-}$ and $G$.
We will show (in Theorem \ref{thm:maingt}) that, under certain conditions,
there is a kind of {}``tri-density constant''. Let us try to express
this a little more explicitly: If $G$ is a {}``relatively small''
part of $Q$ then we will show that there is a subcube $W$ of $Q$
whose intersections with $E_{+}$ and $E_{-}$ are {}``significant''
proportions of $W$\,. We will formulate this result in a more general
context where $Q$ and $W$ are not necessarily cubes, but are members
of some suitable collection $\mathcal{E}$ of admissible sets. In
Section \ref{sec:crux} we will see the implications of this property
of three subsets for the study of various versions of $BMO$.

In order to express our main result for a more general choice of collections
$\mathcal{E}$ of admissible sets, we need to define some more notions.
Our point of departure for doing this comes from considering two important
examples:

(i) Every cube $Q$ in $\mathbb{R}^{d}$ is of course the union of
$2^{d}$ non overlapping subcubes, each having volume $2^{-d}\lambda(Q)$.
Then each of these subcubes can of course in turn be subdivided into
$2^{d}$ non overlapping subcubes of volume $2^{-2d}\lambda(Q)$,
.... and this process can be continued indefinitely.

(ii) Every special rectangle $Q$ in $\mathbb{R}^{d}$ is of course
the union of $2$ non overlapping special rectangles, each having
volume $2^{-1}\lambda(Q)$. Then each of these special rectangles
can of course in turn subdivided into $2$ non overlapping special
rectangles of volume $2^{-2d}\lambda(Q)$, .... and this process too
can be continued indefinitely.

Here is a notion which incorporates these two examples, and ultimately
other examples. 
\begin{defn}
\label{def:multilevel}Let $E$ be an admissible set in $\mathbb{R}^{d}$.
Let $M\ge2$ be a positive integer. We will say that the doubly indexed
sequence $\left\{ E_{j,k}\right\} _{j\ge0,1\le k\le M^{j}}$ of admissible
sets is a \textit{multilevel decomposition of $E$ with multiplicity
$M$} if it satisfies the following conditions:

(i) $E_{0,1}=E$, and, more generally, \[
E=\bigcup_{k=1}^{M^{j}}E_{j,k}\]
for each fixed $j\ge0$\,.

(ii) For each fixed $j\ge1$ the sets $E_{j,k}$ satisfy $\lambda\left(E_{j,k}\cap E_{j,k'}\right)=0$
whenever $k\ne k'$.

(iii) $\lambda(E_{j,k})=M^{-j}\lambda(E)$ for each $j\ge0$ and $k\in\left\{ 1,2,....,M^{j}\right\} $\,.

(iv) For each fixed $j\ge1$, and $k\in\left\{ 1,2,...,M^{j}\right\} $,
the set $E_{j,k}$ is the union of $M$ sets from among the $M^{j+1}$
sets $E_{j+1,m}$. More explicitly, \[
E_{j,k}=\bigcup_{m=M(k-1)+1}^{Mk}E_{j+1,m}\,.\]

(v) The diameters of the sets $E_{j,k}$ tend to zero uniformly as
$j$ tends to infinity, i.e.,\begin{equation}
\lim_{j\to\infty}\left(\max_{1\le k\le M^{j}}\mathrm{diam}E_{j,k}\right)=0\,.\label{eq:diamtz}\end{equation}

\end{defn}
The preceding definition leads us immediately to this next one. 
\begin{defn}
\label{def:Mmultidecomp}Let $\mathcal{E}$ be a collection of admissible
subsets of $\mathbb{R}^{d}$ and let $M\ge2$ be an integer. We will
say that $\mathcal{E}$ is \textit{$M$-multidecomposable} if every
set $E\in\mathcal{E}$ has a multilevel decomposition $\left\{ E_{j,k}\right\} _{j\ge0,1\le k\le M^{j}}$
of multiplicity $M$ where all of the sets $E_{j,k}$ are also in
$\mathcal{E}$. 
\end{defn}
So, of course, for any open subset $D$ of $\mathbb{R}^{d}$, the
collections $\mathcal{Q}(D)$ and $\mathcal{D}(D)$ are both $2^{d}$-multidecomposable,
and the collection $\mathcal{W}(D)$ is $2$-multidecomposable. It
is probably easy to show that the collection $\mathcal{K}(D)$ is
also $2$-decomposable.

\subsection{Our main {}``geometrical'' result. The promised {}``balancing
act'' between three subsets of a (generalized) cube. }
\begin{thm}
\label{thm:maingt}Let $\mathcal{E}$ be a $M$-multidecomposable
collection of admissible subsets of $\mathbb{R}^{d}$ for some $M\ge2$.
Let $\tau$ be a positive number. Let $\delta$ be a bi-density constant
for $\mathcal{E}$.

Suppose that $Q$ is a set in $\mathcal{E}$ and there exist three
pairwise disjoint measurable sets $E_{+}$, $E_{-}$ and $G$ which
satisfy

\[
Q=E_{+}\cup E_{-}\cup G\]
 and \begin{equation}
\min\left\{ \lambda(E_{+}),\lambda(E_{-})\right\} >\tau\lambda(G)\,.\label{eq:tlm}\end{equation}
 Then there exists a subset $W$ of $Q$ such that $W\in\mathcal{E}$
and \begin{equation}
\min\left\{ \lambda(E_{+}\cap W),\lambda(E_{-}\cap W)\right\} \ge s\lambda(W)\label{eq:pdev}\end{equation}
 where \[
s=\left\{ \begin{array}{ccc}
{\displaystyle \min\left\{ \frac{\tau-\tau^{2}}{M(1+\tau)},\delta\right\} } & , & 0<\tau\le\sqrt{2}-1\\
\\{\displaystyle \min\left\{ \frac{3-2\sqrt{2}}{M},\delta\right\} } & , & \sqrt{2}-1\le\tau\,.\end{array}\right.\]

\end{thm}
It will be convenient to explicitly state some immediate consequences
of Theorem \ref{thm:maingt} for some special choices of $Q$ and
$\mathcal{E}$\,, 
\begin{cor}
\label{cor:cmaingt}Suppose that $Q$ is, respectively (i), a cube,
or (ii) a dyadic cube or (iii) a special rectangle in $\mathbb{R}^{d}$
and that, respectively,

(i) $\mathcal{E=Q}(Q)$, or (ii) $\mathcal{E=D}(Q)$, or (iii) $\mathcal{E=W}(Q)$.
Suppose that $Q$ is the disjoint union of the three sets $E_{-}$,
$E_{+}$ and $G$ and that \[
\min\left\{ \lambda(E_{+}),\lambda(E_{-})\right\} >\left(\sqrt{2}-1\right)\lambda(G)\,.\]
 Then there exists a set $W\subset Q$ which is, respectively, (i)
a cube, or (ii) a dyadic cube or (iii) a special rectangle, and which
satisfies \[
\min\left\{ \lambda(E_{+}\cap W),\lambda(E_{-}\cap W)\right\} \ge s\lambda(W)\]
 where, respectively, (i) $s=2^{-d}(3-2\sqrt{2})$\,, or (ii) $s=2^{-d}(3-2\sqrt{2})$
(again), or (iii) $s=(3-2\sqrt{2})/2$.
\end{cor}
\noindent \textit{Proof of Corollary \ref{cor:cmaingt}.} We simply
apply Theorem \ref{thm:maingt}, substituting the known values for
$\delta$ and $M$ in the formula $s=\min\left\{ \frac{3-2\sqrt{2}}{M},\delta\right\} $
in each of the three cases. $\qed$

Theorem \ref{thm:maingt} and Corollary \textit{\ref{cor:cmaingt}}
motivate us to introduce another notion. This notion will enable a
convenient formulation of the main question raised by this paper,
and also a convenient proof of the consequences that an affirmative
answer to that question would have. 
\begin{defn}
\label{def:jsconst}Let $\mathcal{E}$ be a collection of admissible
subsets of $\mathbb{R}^{d}$. Let $\tau$ and $s$ be positive numbers
with the following property:

Let $Q$ be a set in $\mathcal{E}$ and let $E_{+}$, $E_{-}$ and
$G$ be arbitrary pairwise disjoint admissible sets whose union is
$Q$. Suppose that\begin{equation}
\min\left\{ \lambda(E_{+}),\lambda(E_{-})\right\} >\tau\lambda(G)\,.\label{eq:isms}\end{equation}
 Then there exists a set $W\subset Q$ which is also in $\mathcal{E}$
and for which\begin{equation}
\min\left\{ \lambda(E_{+}\cap W),\lambda(E_{-}\cap W)\right\} \ge s\lambda(W)\,.\label{eq:blerk}\end{equation}
 Then we will say that \textit{$\left(\tau,s\right)$ is a John-Strömberg
pair for $\mathcal{E}$}. 
\end{defn}

\begin{rem}
The preceding definition is formulated for all possible positive values
of $\tau$ and $s$. However, for our particular applications we are
interested only in cases where $\tau<1/2$. This is why, in the formulation
of Question A, we apply this latter restriction to $\tau$\,.\end{rem}
\begin{example}
\label{exa:emaingt}We can reformulate Corollary \ref{cor:cmaingt}
as follows: We take $\tau=\sqrt{2}-1$. Then, for any open set $D$
of $\mathbb{R}^{d}$, we have that

(i) $\left(\sqrt{2}-1,2^{-d}(3-2\sqrt{2})\right)$ is a John-Strömberg
pair for $\mathcal{Q}(D)$ and also for $\mathcal{D}(D)$.

(ii) $\left(\sqrt{2}-1,\frac{3-2\sqrt{2}}{2}\right)$ is a John-Strömberg
pair for $\mathcal{W}(D)$. 
\end{example}

\begin{rem}
\label{rem:unitcube}In the special case where $\mathcal{E}=\mathcal{Q}(D)$
for some open subset $D$ of $\mathbb{R}^{d}$ then, obviously, $\left(\tau,s\right)$
is a John-Str\"omberg pair if and only if the condition appearing
in the second paragraph of Definition \ref{def:jsconst} holds for
just one particular cube $Q$ in $\mathbb{R}^{d}$, for example for
the unit cube $Q=[0,1]^{d}$. The cube $Q$ does not have to be contained
in $D$. (In fact the particular choice of $D$ is irrelevant in this
case. $\left(\tau,s\right)$ is a John-Str\"omberg pair for $\mathcal{Q}(D)$
if and only if it is a John-Str\"omberg pair for $\mathcal{Q}(\mathbb{R}^{d})$.)
\end{rem}

\begin{rem}
\label{rem:tauIsLittle}We can now express Question A concisely in
the language of Definition \ref{def:jsconst}. Question A simply asks
whether there exist two absolute constants $s>0$ and $\tau\in(0,1/2)$
such that $(\tau,s)$ is a John-Strömberg pair for $\mathcal{Q}(\mathbb{R}^{d})$
for \textit{every} $d\in\mathbb{N}$. 
\end{rem}

\begin{rem}
\label{rem:JSimpliesBIDENSITY}In Definition \ref{def:jsconst} the
set $G$ may be chosen to be empty. Therefore any number $s$, which
happens to form a John-Strömberg pair $(\tau,s)$ for $\mathcal{E}$
with some positive number $\tau$, will also automatically be a bi-density
constant for $\mathcal{E}$. The particular value of $\tau$ is immaterial
here. 
\end{rem}

\begin{rem}
\label{rem:sbth}Although the only restriction that we explicitly
impose on $s$ is that it has to be positive, we see from (\ref{eq:blerk})
that some cube $W$, which has positive measure, must satisfy \[
\lambda(W)\ge\lambda(E_{+}\cap W)+\lambda(E_{-}\cap W)\ge2\min\left\{ \lambda(E_{+}\cap W),\lambda(E_{-}\cap W)\right\} \ge2s\lambda(W)\,.\]
Thus the constant $s$ in the above definition can never be greater
than $1/2$.
\end{rem}

\begin{rem}
\label{rem:dim}The first of the two inequalities which appear in
Definition \ref{def:jsconst} (and also in Question A) is strict and
the second is not. But would it really change anything if both or
neither of them were strict? It would seem that not. We refer to Remark
\ref{rem:equality} in connection with this issue.
\end{rem}

\noindent \textit{Proof of Theorem \ref{thm:maingt}.} In view of
the regularity of $\lambda$ there exist compact subsets $H_{+}$
and $H_{-}$ of $E_{+}$ and $E_{-}$ respectively such that \[
\min\left\{ \lambda(H_{+}),\lambda(H_{-})\right\} >\tau\lambda(G)\,.\]
 If there exists $W\in\mathcal{E}$ such that $W\subset Q$ and \[
\min\left\{ \lambda(H_{+}\cap W),\lambda(H_{-}\cap W)\right\} \ge s\lambda(W)\]
then obviously $W$ also satisfies (\ref{eq:pdev}). This means that
we may assume without loss of generality that $E_{+}$ and $E_{-}$
are themselves compact sets. Since they are also disjoint, it follows
that\[
\rho:=\mathrm{dist}\left(E_{+},E_{-}\right)>0\,.\]

We will say that the set $W$ is a \textit{good set} if $W\in\mathcal{E}$
and $W\subset Q$ and \[
\min\left\{ \lambda(E_{+}\cap W),\lambda(E_{-}\cap W)\right\} >\tau\lambda(G\cap W)\,.\]

(We mention that a slight variant of this definition will play a role
later, in Subsection \ref{sub:diffspec}, in particular in the proof
of Theorem \ref{thm:newred}, which will have some similarities with
some of the arguments here.)

Let $\left\{ Q_{j,k}\right\} _{j\ge0,1\le k\le M^{j}}$ be a multilevel
decomposition of $Q$ of multiplicity $M$ where all of the sets $Q_{j,k}$
are in $\mathcal{E}$. A sequence $\left\{ Q^{(n)}\right\} _{0\le n<\ell}$
of sets in $\mathcal{E}$, where $\ell$ can be finite or infinite,
will be called a \textit{chain} if $Q^{(0)}=Q$ and, if for each $n$
such that $1\le n<\ell$ we have $Q^{(n)}=Q_{n,k_{n}}$ for some integer
$k_{n}\in\left\{ 1,2,...,M^{n}\right\} $. When $\ell=\infty$, we
have (by (\ref{eq:diamtz})) that $\lim_{n\to\infty}\mathrm{diam}Q^{(n)}=0$.
So there must exist some $n_{1}$ such that $\mathrm{diam}Q^{(n)}<\rho$
for all $n\ge n_{1}$. Thus, for $n\ge n_{1}$ the set $Q^{(n)}$
cannot intersect with both of $E_{+}$ and $E_{-}$ and so cannot
be a good set.

Let construct a particular chain $\left\{ Q^{(n)}\right\} _{0\le n<\ell}$
in the following way. We of course have to start with $Q^{(0)}=Q$.
If among all the sets $Q_{1,k}$ for $k\in\left\{ 1,2,...,M\right\} $
there is no set which is good, then we set $\ell=1$ and our construction
is complete. Otherwise we choose $Q^{(1)}$ to be a good set from
the above list. Next we check whether, among those of the sets $Q_{2,k}$
for $k\in\left\{ 1,2,...,M^{2}\right\} $ which are contained in $Q^{(1)}\,$,
there is one which is a good set. If so we choose $Q^{(2)}$ to be
such a set. If not, we set $\ell=2$ and our construction is complete.
The continuation of this process is now clear. At the $n$th stage
we seek a good set, which we will call $Q^{(n)}$, from among those
of the sets $Q_{n,k}$ which are contained in $Q^{(n-1)}$. If no
such set exists, then we set $\ell=n$ and the construction is complete.
In view of the arguments presented in the previous paragraph, we must
have $\ell<\infty$, i.e., the construction necessarily has to terminate
after finitely many steps.

Thus we have obtained a good set $Q^{(\ell-1)}$ which is one of the
sets $Q_{\ell-1,k}$ for some integer $k$ and it will enable us to
complete the proof of the theorem. Let us denote the $M$ sets of
the form $Q_{\ell,m}$ which are contained in $Q^{(\ell-1)}$ by $W_{1}$,
$W_{2}$,...., $W_{M}$. By our construction, none of these sets are
good sets. For the rest of the proof we will assume that \[
\lambda(Q^{(\ell-1)}\cap E_{-})\le\lambda(Q^{(\ell-1)}\cap E_{+})\,.\]
 If the reverse inequality holds then we will simply use an exact
analogue of the proof that we are about to give, where we will simply
interchange the roles of $E_{+}$ and $E_{-}$.

We have to consider three cases:

Case (i). Suppose that $\lambda(Q^{(\ell-1)}\cap G)=0$. Then, since
$\delta$ is a bi-density constant for $\mathcal{E}$, and since $\lambda(Q^{(\ell-1)}\cap E_{+})$
and $\lambda(Q^{(\ell-1)}\setminus E_{+})=\lambda(Q^{(\ell-1)}\cap E_{-})$
are both positive, there exists $W\in\mathcal{E}$ such that $W\subset Q^{(\ell-1)}$
and \begin{equation}
\min\left\{ \lambda(W\cap E_{+}),\lambda(W\cap E_{-})\right\} \ge\delta\lambda(W)\ge s\lambda(W)\,,\label{eq:tcpt}\end{equation}
 completing the proof of the theorem.

Case (ii). Suppose that $\lambda(W_{m}\cap E_{-})=\lambda(W_{m})$
for some $m\in\left\{ 1,2,...,M\right\} $\,.

Then \begin{eqnarray*}
\lambda(Q^{(\ell-1)}\cap E_{+}) & \ge & \lambda(Q^{(\ell-1)}\cap E_{-})\ge\lambda(W_{m}\cap E_{-})=\lambda(W_{m})=M^{-1}\lambda(Q^{(\ell-1)})\\
 & \ge & s\lambda(Q^{(\ell-1)})\,.\end{eqnarray*}

So we see that in this case the set $W=Q^{(\ell-1)}$ has the properties
required to complete the proof of the theorem.

Case (iii). This is the remaining case where cases (i) and (ii) are
excluded. I.e., we have $\lambda(Q^{(\ell-1)}\cap G)>0$ and \begin{equation}
\lambda(W_{m}\cap E_{-})<\lambda(W_{m})\mbox{ for all }m\in\left\{ 1,2,...,M\right\} \,.\label{eq:omp-1}\end{equation}

The inequality \[
\lambda(Q^{(\ell-1)}\cap E_{-})>\tau\lambda(Q^{(\ell-1)}\cap G)>0\]
 (which holds because $Q^{(\ell-1)}$ is good) can be rewritten (in
view of condition (ii) of Definition \ref{def:multilevel}) as\begin{equation}
\sum_{m=1}^{M}\lambda(W_{m}\cap E_{-})>\sum_{m=1}^{M}\tau\lambda(W_{m}\cap G)>0\,.\label{eq:qtyp}\end{equation}
Let $N$ be the (possibly empty) set of all integers $m$ in $\left\{ 1,2,...,M\right\} $
which satisfy \[
\lambda(W_{m}\cap G)=0\,.\]
If $\lambda(W_{m}\cap E_{-})>0$ for some $m\in N$, then, since (\ref{eq:omp-1})
also holds, we can again, analogously to what was done in Case (i),
invoke the bi-density condition to obtain some $W\in\mathcal{E}$
with $W\subset W_{m}$ which satisfies (\ref{eq:tcpt}) and so completes
the proof. This means that we can now assume that $\lambda(W_{m}\cap E_{-})=0$
for all $m\in N$. Therefore (\ref{eq:qtyp}) can be rewritten as\[
\sum_{m\in\{1,2,..,M\}\backslash N}\lambda(W_{m}\cap E_{-})>\sum_{m\in\{1,2,..,M\}\backslash N}\tau\lambda(W_{m}\cap G)>0\,.\]
 It follows that there exists at least one $m\in\left\{ 1,2,...,M\right\} \setminus N$
which satisfies \[
\lambda(W_{m}\cap E_{-})>\tau\lambda(W_{m}\cap G)>0\,.\]
 Recall that, by our construction, $W_{m}$ is not good. Therefore
we must have \[
\lambda(W_{m}\cap E_{+})\le\tau\lambda(W_{m}\cap G)\,.\]
We use this and the preceding inequality to obtain that\begin{eqnarray*}
\lambda(W_{m}\cap E_{-}) & > & \tau\lambda(W_{m}\cap G)=\tau\left[\lambda(W_{m})-\lambda(W_{m}\cap E_{-})-\lambda(W_{m}\cap E_{+})\right]\\
 & \ge & \tau\left[\lambda(W_{m})-\lambda(W_{m}\cap E_{-})-\tau\lambda(W_{m}\cap G)\right]\\
 & \ge & \tau\left[\lambda(W_{m})-\lambda(W_{m}\cap E_{-})-\tau\lambda(W_{m})\right]\\
 & = & \left(\tau-\tau^{2}\right)\lambda(W_{m})-\tau\lambda(W_{m}\cap E_{-})\,.\end{eqnarray*}
 This implies that \[
\lambda\left(W_{m}\cap E_{-}\right)>\frac{\tau-\tau^{2}}{1+\tau}\lambda(W_{m})=M^{-1}\frac{\tau-\tau^{2}}{1+\tau}\lambda(Q^{(\ell-1)})\,.\]
Let us consider the case where $\tau$ is in the range $0<\tau\le\sqrt{2}-1$.
Since \[
\lambda(Q^{(\ell-1)}\cap E_{+})\ge\lambda(Q^{(\ell-1)}\cap E_{-})\ge\lambda(W_{m}\cap E_{-})\]
and $M^{-1}\frac{\tau-\tau^{2}}{1+\tau}=s$ for this range of values
of $\tau$, we see that in this case $W=Q^{(\ell-1)}$ satisfies (\ref{eq:pdev}).
In the remaining case, where $\tau>\sqrt{2}-1$, we simply observe
that the given condition (\ref{eq:tlm}) also holds when $\tau$ is
replaced by the smaller number $\sqrt{2}-1$ and so we can apply the
same argument for this smaller value of $\tau$ to obtain the required
conclusion. This completes the proof of the theorem. $\qed$

\section{\label{sec:crux}Applying our {}``geometrical'' result. The non
increasing rearrangement of a BMO function.}

Suppose that a function $f$ of $d$-variables is in $BMO\left(D,\mathcal{Q}(D)\right)$
for some cube $D$ in $\mathbb{R}^{d}$. It is known \cite{bds},
\cite{garsiarodemich74} that this implies that the function of one
variable $f^{*}$, i.e., the non increasing rearrangement of $f$,
is in $BMO\left(I,\mathcal{Q}(I)\right)$ and that \begin{equation}
\left\Vert f^{*}\right\Vert _{BMO\left(I,\mathcal{Q}(I)\right)}\le C\left\Vert f\right\Vert _{BMO\left(D,\mathcal{Q}(D)\right)}\label{eq:ffstar}\end{equation}
for some constant $C$ depending only on $d$. (Of course here we
are in fact considering only the restriction of $f^{*}$ to the interval
$\left(0,\lambda(D)\right)$.) Apparently the optimal value of $C$
for which (\ref{eq:ffstar}) holds for all such $f$ is not yet known.
Nor is it known yet whether or not $C$ can be chosen to in fact be
independent of $d$. It is known \cite{koren1990} that $C=1$ when
$d=1$. An analogous result holds for dyadic intervals \cite{klemes}.

In this section we wish to obtain inequalities analogous to (\ref{eq:ffstar}),
in terms of the functional of John and Strömberg, namely inequalities
of the form \begin{equation}
\left\Vert f^{*}\right\Vert _{BMO\left(I,\mathcal{Q}(I)\right)}^{(\mathbf{J},\sigma)}\le C\left\Vert f\right\Vert _{BMO\left(D,\mathcal{E}\right)}^{(\mathbf{J},s)}\label{eq:vmjet}\end{equation}
where, as above, $I$ is the interval $I=(0,\lambda(D))$ and $s$
and $\sigma$ are suitably chosen numbers in $(0,1/2]$. We will be
particularly interested in the cases where $D$ is a cube or a special
rectangle and $\mathcal{E}$ is $\mathcal{Q}(D)$ or $\mathcal{W}(D)$.
But our results for these will be consequences of an analogous result
for more general choices of $D$ and $\mathcal{E}$.

The constant $C$ in our versions of (\ref{eq:vmjet}) will be $C=1$
and this is apparently the best possible constant for any and every
choice of the parameters $\sigma$ and $s$ in $(0,1/2]$. We will
be able to take our parameter $\sigma$ to be any number satisfying
\[
\frac{1}{2}\ge\sigma>\frac{2\sqrt{2}-2}{2\sqrt{2}-1}\approx0.453082\,.\]
(As mentioned earlier, values of $\sigma>1/2$ are not particularly
interesting, since $\left\Vert f^{*}\right\Vert _{BMO\left(I,\mathcal{Q}(I)\right)}^{(\mathbf{J},\sigma)}$
can equal $0$ also when $f$ is not a constant function.) Our parameter
$s$ will depend on our choice of $\mathcal{E}$. In fact it will
be exactly the parameter given by the formula $s=\min\left\{ \frac{3-2\sqrt{2}}{M},\delta\right\} $
in Theorem \ref{thm:maingt}. Thus, exactly as in Corollary \ref{cor:cmaingt},
we will have $s=2^{-d}(3-2\sqrt{2})$ if $\mathcal{E}$ is either
$\mathcal{Q}(D)$ or $\mathcal{D}(D)$ and $D$ is, respectively a
cube, or a dyadic cube, and, furthermore, we will have $s=(3-2\sqrt{2})/2$
if $\mathcal{E=W}(D)$ and $D$ is a special rectangle.

We should now stress that we will not obtain our versions of the inequality
(\ref{eq:vmjet}), which we have spent the last few paragraphs describing,
in as much generality as the reader may have been led to expect. We
will only obtain them for the very special class of those measurable
functions $f$ on $D$ which take only non negative integer values,
and which satisfy $\left\Vert f\right\Vert _{BMO\left(D,\mathcal{E}\right)}^{(\mathbf{J},s)}\le1/2$.
But, in view of the results of Section \ref{sec:simplifyJohnStromberg},
notably Theorem \ref{thm:temjemt-1}, this will be sufficient for
our subsequent purposes.

It is perhaps surprising, and perhaps even amusing, that, at a certain
stage, it will turn out to be sufficient to consider an even more
restricted subclass of the very special class of functions just referred
to, namely those functions which only assume at most three consecutive
integer values, which may just as well be $0$, $1$ and $2$.

All the results which we have just described are immediate consequences
of the following theorem (Theorem \ref{thm:UltraNewMain}). In turn
Theorem \ref{thm:UltraNewMain} will follow from a more general theorem
(Theorem \ref{thm:Nultra}) which will be formulated after this one.
\begin{rem}
In the formulation of both of the following two theorems we consider
a collection $\mathcal{E}$ of admissible sets and a particular set
$Q\in\mathcal{E}$. Instead of $\mathcal{E}$, we use the subcollection
$\mathcal{E}(Q)$ of all sets in $\mathcal{E}$ which are contained
in $Q$. We are essentially forced to do this because our function
$f$ is defined only on $Q$. Obviously, if $f$ happens to be defined
on all sets in $\mathcal{E}$ the conclusions of both theorems will
remain true if we replace $\mathcal{E}(Q)$ by the larger collection
$\mathcal{E}$.\end{rem}
\begin{thm}
\label{thm:UltraNewMain}Let $\mathcal{E}$ be an $M$-multidecomposable
collection for some $M\ge2$. Let $\delta$ be a bi-density constant
for $\mathcal{E}$. Let $Q$ be a set in $\mathcal{E}$ and let $\mathcal{E}(Q)$
be the collection of all sets in $\mathcal{E}$ which are contained
in $Q$. Suppose that the function $f:Q\to\mathbb{N}\cup\left\{ 0\right\} $
is measurable and satisfies $\left\Vert f\right\Vert _{BMO(Q,\mathcal{E}(Q))}^{(\mathbf{J},s)}\le1/2$
for $s=\min\left\{ \delta,\frac{3-2\sqrt{2}}{M}\right\} $. Let $\sigma$
be a number in the range \begin{equation}
\frac{2\sqrt{2}-2}{2\sqrt{2}-1}<\sigma\le\frac{1}{2}\,.\label{eq:ams}\end{equation}
Then the function $f^{*}:(0,\lambda(Q))\to[0,\infty)$, i.e., the
non increasing rearrangement of $f$ restricted to the interval $I:=(0,\lambda(Q))$
satisfies \[
\left\Vert f^{*}\right\Vert _{BMO(I,\mathcal{Q}(I))}^{(\mathbf{J},\sigma)}\le\left\Vert f\right\Vert _{BMO(Q,\mathcal{E}(Q))}^{(\mathbf{J},s)}\,.\]

\end{thm}
We want to obtain this theorem as a consequence of the following somewhat
more abstract and general theorem, which is formulated in terms of
John-Strömberg pairs $(\tau,s)$ for the collection $\mathcal{E}$.
By introducing this extra level of abstraction we also make it possible
to formulate the consequences of an affirmative answer to Question
A in a (hopefully) clearer and more organized way.
\begin{thm}
\label{thm:Nultra}Let $\mathcal{E}$ be a collection of admissible
subsets of $\mathbb{R}^{d}$. Let $Q$ be a set in $\mathcal{E}$
and let $\mathcal{E}(Q)$ be the collection of all sets in $\mathcal{E}$
which are contained in $Q$. Let $\tau\in(0,1/2)$ and $s\in(0,1/2)$
be such that $(\tau,s)$ is a John-Strömberg pair for $\mathcal{E}$.
Suppose that the function $f:Q\to\mathbb{N}\cup\left\{ 0\right\} $
is measurable and satisfies $\left\Vert f\right\Vert _{BMO(Q,\mathcal{E}(Q))}^{(\mathbf{J},s)}\le1/2$\,.
Let $\sigma$ be a number in the range \begin{equation}
\frac{2\tau}{1+2\tau}<\sigma\le\frac{1}{2}\,.\label{eq:Nams}\end{equation}
 Then the function $f^{*}:(0,\lambda(Q))\to[0,\infty)$, i.e., the
non increasing rearrangement of $f$ restricted to the interval $I:=(0,\lambda(Q))$,
satisfies \begin{equation}
\left\Vert f^{*}\right\Vert _{BMO(I,\mathcal{Q}(I))}^{(\mathbf{J},\sigma)}\le\left\Vert f\right\Vert _{BMO(Q,\mathcal{E}(Q))}^{(\mathbf{J},s)}\,.\label{eq:Nyphwu}\end{equation}
\end{thm}
\begin{rem}
\label{rem:conversethm}We find it interesting that the {}``geometric''
condition of the existence of a John-Strömberg pair turns out to be
in some sense {}``almost equivalent'' to the {}``analytic'' condition
expressed by the inequality (\ref{eq:Nyphwu}). This is revealed by
combining Theorem \ref{thm:Nultra} with an auxiliary result (Theorem
\ref{thm:Ninverse}) which we will defer to the end of this section,
since it will not be needed for obtaining our main result in Section
\ref{sec:Putting}. Theorem \ref{thm:Ninverse} will be a sort of
converse to Theorem \ref{thm:Nultra}. (The fact that we use $\mathcal{E}(Q)$
rather than $\mathcal{E}$ in the previous two theorems makes their
results more closely comparable with the result of Theorem \ref{thm:Ninverse}.) 
\end{rem}
\noindent \textit{Proofs of Theorems \ref{thm:UltraNewMain} and\ref{thm:Nultra}.}

Let us first show that Theorem \ref{thm:UltraNewMain} is a consequence
of \ref{thm:Nultra}. Since the conclusions of both theorems are the
same, this simply amounts to showing that the conditions imposed on
$\mathcal{E}$ and $s$ and $\sigma$ in Theorem \ref{thm:UltraNewMain}
suffice to guarantee that $\mathcal{E}$ and $s$ and $\sigma$ satisfy
the hypotheses of Theorem \ref{thm:Nultra} for some suitable choice
of $\tau$. In fact we will choose $\tau=\sqrt{2}-1$ so that $\frac{2\tau}{1+2\tau}=\frac{2\sqrt{2}-2}{2\sqrt{2}-1}$.
So when, in Theorem \textit{\ref{thm:UltraNewMain}}, we require $\sigma$
to satisfy (\ref{eq:ams}), this ensures that $\sigma$ will be in
the range specified in (\ref{eq:Nams}). It remains only to check
that $(\tau,s)$ is a John-Strömberg pair for $\mathcal{E}$ when
$\tau=\sqrt{2}-1$ and $s=\min\left\{ \delta,\frac{3-2\sqrt{2}}{M}\right\} $.
But this is exactly what is stated by Theorem \ref{thm:maingt} for
these choices of $\tau$ and $s\,$.

Thus we can now turn to the proof of Theorem \ref{thm:Nultra}. Since
we only have to deal with sets of $\mathcal{E}$ which are contained
in $Q$ we may suppose from here onwards that $\mathcal{E}=\mathcal{E}(Q)$.

The fact that all values taken by $f$ are in $\mathbb{N\cup}\left\{ 0\right\} $
readily implies that the same is true for all values of $f^{*}$.
To explain this more precisely, since $\lambda(Q)<\infty$, we can
invoke (\ref{eq:df4}) to obtain that \begin{equation}
\lambda\left(\left\{ x\in Q:f(x)=m\right\} \right)=\left|\left\{ t\in\left(0,\lambda(Q)\right):f^{*}(t)=m\right\} \right|\,\mbox{for all }m\in\mathbb{N}\cup\left\{ 0\right\} \,.\label{eq:intval}\end{equation}
Then the fact that \[
\sum_{m=0}^{\infty}\lambda\left(\left\{ x\in Q:f(x)=m\right\} \right)=\lambda(Q)\]
implies that the subset of $\left(0,\lambda(Q)\right)$ where $f^{*}$
takes non integer values has measure zero.

Let us first dispose of three easier special cases where we can readily
see that (\ref{eq:Nyphwu}) holds.

The first of these cases is when $f$ takes only one value (on sets
of positive measure). Then of course $\left\Vert f\right\Vert _{BMO(Q,\mathcal{E})}^{(\mathbf{J},s)}=\left\Vert f^{*}\right\Vert _{BMO(I,\mathcal{Q}(I))}^{(\mathbf{J},\sigma)}=0$
no matter how we choose $\sigma$ and $s$. So we obtain (\ref{eq:Nyphwu}).

The second case is when $f$ takes only two values (on sets of positive
measure). In this case $f=a\chi_{A}+b\chi_{Q\setminus A}$ for some
$A\subset Q$ such that $0<\lambda(A)<\lambda(Q)$ and where $b<a$.
Since (cf.~Remark \ref{rem:JSimpliesBIDENSITY}) $s$ is also a bi-density
constant for $\mathcal{E}$, there exists some set $W_{0}\in\mathcal{E}$
contained in $Q$ such that \[
\min\left\{ \lambda(W_{0}\cap A),\lambda(W_{0}\setminus A)\right\} \ge s\lambda(W_{0})\,.\]
This means that the restriction of $\left(f\chi_{W_{0}}\right)^{*}$
to the interval $\left(0,\lambda(W_{0})\right)$ is given by the formula
$\left(f\chi_{W_{0}}\right)^{*}=a\chi_{(0,r)}+b\chi_{[r,\lambda(W_{0}))}$
for some number $r\in[s\lambda(W_{0}),(1-s)\lambda(W_{0})]$\,. For
all choices of the number $u$ which satisfy $u\in\left(0,s\lambda(W_{0})\right)$
we have $u\in(0,r)$ and \[
u+(1-s)\lambda(W_{0})\in[r,\lambda(W_{0}))\,.\]
Consequently, $\left(f\chi_{W_{0}}\right)^{*}(u)=a$ and $\left(f\chi_{W_{0}}\right)^{*}\left(u+(1-s)\lambda(W_{0})\right)=b$.
This implies, by Proposition \ref{pro:prop}, that $\mathbf{J}(f,W_{0},s)=\frac{a-b}{2}$
for this particular set $W_{0}$. Proposition \ref{pro:prop} also
tells us that $\mathbf{J}(f,W,s)\le\frac{a-b}{2}$ for all other sets
$W$ in $\mathcal{E}$, so we conclude that $\left\Vert f\right\Vert _{BMO(Q,\mathcal{E})}^{(\mathbf{J},s)}=\frac{a-b}{2}$.
Since $1/2$ is a bi-density constant for $\mathcal{Q}(I)$ and $\sigma\in(0,1/2]$
it follows that $\sigma$ is also a bi-density constant for $\mathcal{Q}(I)$.
So we can show, by applying reasoning to $f^{*}$ and $\sigma$, which
is exactly analogous to the reasoning just applied to $f$ and $s$,
that $\left\Vert f^{*}\right\Vert _{BMO(I,\mathcal{Q}(I))}^{(\mathbf{J},\sigma)}=\frac{a-b}{2}$.
Thus we see that (\ref{eq:Nyphwu}) holds in this case also.

The third and last of these easier cases is when $\left\Vert f\right\Vert _{BMO(Q,\mathcal{E})}^{(\mathbf{J},s)}<1/2$.
We will deal with this case by showing that here $f$ has to be a
constant. The fact that $f$ and $f^{*}$ are integer valued, together
with the formula (\ref{eq:ffnn}), tells us that the two functionals
$\left\Vert f^{*}\right\Vert _{BMO(I,\mathcal{Q}(I))}^{(\mathbf{J},\sigma)}$
and $\left\Vert f\right\Vert _{BMO(Q,\mathcal{E})}^{(\mathbf{J},s)}$
are each infima of appropriate subsets of the set $\left\{ (n-1)/2:n\in\mathbb{N}\right\} $
of non negative half-integers. Thus these functionals themselves can
only take non negative half integer values. More explicitly, if $\left\Vert f\right\Vert _{BMO(Q,\mathcal{E})}^{(\mathbf{J},s)}<1/2$,
then $\left\Vert f\right\Vert _{BMO(Q,\mathcal{E})}^{(\mathbf{J},s)}=0$.
If $f$ is constant a.e., then so is $f^{*}$ and, as in the first
case, we see, trivially that (\ref{eq:Nyphwu}) holds. We shall now
show that this is the only possibility. Suppose, on the contrary,
that $f$ is not a constant a.e. Then let $b$ be the smallest non
negative integer (there must exist at least two such integers) for
which the set $B=\left\{ x\in Q:f(x)=b\right\} $ has positive measure.
Then the set \[
A=\left\{ x\in Q:f(x)>b\right\} =\left\{ x\in Q:f(x)\ge b+1\right\} \]
must also have positive measure and we must have $\lambda(B)+\lambda(A)=\lambda(Q)$.

Let us now define the functions $\varphi:\mathbb{R}\to\mathbb{R}$
and $\psi:\mathbb{R}\to\mathbb{R}$ by $\psi(t)=\max\left\{ t,b\right\} $
and then $\varphi(t)=\min\left\{ b+1,\psi(t)\right\} $. Obviously
$\psi$ and therefore also $\varphi$ are both 1-Lipschitz functions
and therefore, by Lemma \ref{lem:JLIP}, we have that \begin{equation}
\left\Vert \varphi\circ f\right\Vert _{BMO(Q,\mathcal{E})}^{(\mathbf{J},s)}\le\left\Vert f\right\Vert _{BMO(Q,\mathcal{E})}^{(\mathbf{J},s)}=0\,.\label{eq:mma}\end{equation}
 We also have $\varphi\circ f=b\chi_{B}+(b+1)\chi_{A}$ almost everywhere.
Therefore, the same calculation that we did in case (ii) for a function
taking \textit{only two different values} a.e., gives here that \[
\left\Vert \varphi\circ f\right\Vert _{BMO(Q,\mathcal{E})}^{(\mathbf{J},s)}=\frac{\left|(b+1)-b\right|}{2}=\frac{1}{2}\,.\]
This contradicts (\ref{eq:mma}) and shows that $f$ indeed must be
a constant, completing our treatment of this case.

Having disposed of these three cases, we can, from this point onwards,
assume that $f$ takes three or more different values on sets of positive
measure, and we can also assume that $\left\Vert f\right\Vert _{BMO(Q,\mathcal{E})}^{(\mathbf{J},s)}=1/2$.
Let us suppose that (\ref{eq:Nyphwu}) does not hold, i.e., that \begin{equation}
\left\Vert f^{*}\right\Vert _{BMO(I,\mathcal{Q}(I))}^{(\mathbf{J},\sigma)}>1/2\,.\label{eq:wa}\end{equation}
 We will complete the proof of Theorem \ref{thm:Nultra} by showing
that (\ref{eq:wa}) leads to a contradiction. It follows from (\ref{eq:wa})
that there exists some subinterval $I^{\prime}=(a,b)$ of $I$ such
that $\mathbf{J}(f^{*},I^{\prime},\sigma)>1/2$. This means that the
number $\alpha=1/2$ is not in the {}``competition'' for the infimum
in the formula analogous to (\ref{eq:altdef}) for $\mathbf{J}(f^{*},I^{\prime},\sigma)$\,.
So, for every $c\in\mathbb{R}$\,, we have \[
\left|\left\{ t\in(a,b):\left|f^{*}(t)-c\right|\le1/2\right\} \right|\le(1-\sigma)(b-a)\,.\]
 If $k$ is any integer and if we choose $c=k+1/2$, then, since $f$
and $f^{*}$ only take integer values, we have that \begin{eqnarray*}
\left\{ t\in(a,b):\left|f^{*}(t)-c\right|\le1/2\right\}  & = & \left\{ t\in(a,b):c-1/2\le f^{*}(t)\le c+1/2\right\} \\
 & = & \left\{ t\in(a,b):f^{*}(t)\in\{k,k+1\}\right\} \,.\end{eqnarray*}
 So, in fact, \begin{equation}
\left|\left\{ t\in(a,b):f^{*}(t)\in\{k,k+1\}\right\} \right|\le(1-\sigma)(b-a)\mbox{ for each }k\in\mathbb{Z}\,.\label{eq:bomp}\end{equation}
 (It may of course happen that this set is empty for some or even
most values of $k$.)

Since \[
(a,b)=\bigcup_{k\in\mathbb{N}}\left\{ t\in(a,b):f^{*}(t)<k\right\} \]
 and since the set $\left\{ t\in(a,b):f^{*}(t)<0\right\} $ is empty,
there exists a unique integer $k\in\mathbb{N}$ such that\begin{equation}
\begin{array}{rcl}
\left|\left\{ t\in(a,b):f^{*}(t)<k-1\right\} \right| & \le & {\displaystyle \frac{\sigma}{2}}(b-a)<\left|\left\{ t\in(a,b):f^{*}(t)<k\right\} \right|\\
\\ & = & \left|\left\{ t\in(a,b):f^{*}(t)\le k-1\right\} \right|\,\end{array}\label{eq:jomp}\end{equation}

The interval $(a,b)$ is the union of the three disjoint sets $\left\{ t\in(a,b):f^{*}(t)<k-1\right\} $,
$\left\{ t\in(a,b):f^{*}(t)\in\{k-1,k\}\right\} $ and $\left\{ t\in(a,b):f^{*}(t)\ge k+1\right\} $.
By (\ref{eq:bomp}) and (\ref{eq:jomp}), the measure of the union
of the first two of these does not exceed $\left(1-\frac{\sigma}{2}\right)(b-a)$.
Therefore we conclude that \begin{equation}
\left|\left\{ t\in(a,b):f^{*}(t)\ge k+1\right\} \right|\ge\frac{\sigma}{2}(b-a)\,.\label{eq:klomp}\end{equation}

This implies that there exists some integer $m\ge k+1$ for which
the set \[
\left\{ t\in(a,b):f^{*}(t)=m\right\} \]
 has positive measure. Let us also show that the set $\left\{ t\in(a,b):f^{*}(t)=k-1\right\} $
has positive measure. Its measure satisfies \begin{eqnarray*}
 &  & \left|\left\{ t\in(a,b):f^{*}(t)=k-1\right\} \right|\\
 & = & \left|\left\{ t\in(a,b):f^{*}(t)\le k-1\right\} \right|-\left|\left\{ t\in(a,b):f^{*}(t)\le k-2\right\} \right|\\
 & = & \left|\left\{ t\in(a,b):f^{*}(t)<k\right\} \right|-\left|\left\{ t\in(a,b):f^{*}(t)<k-1\right\} \right|\end{eqnarray*}
 which, by (\ref{eq:jomp}), is indeed positive.

We now know that on the interval $(a,b)$ the function $f^{*}$ assumes
at least one value strictly larger than $k$ and at least one value
strictly less than $k$. Since $f$ is non increasing, this means
that \begin{equation}
\left\{ t\in(a,b):f^{*}(t)=k\right\} =\left\{ t\in(0,\lambda(Q)):f^{*}(t)=k\right\} \,.\label{eq:fromp}\end{equation}
Now let us define the three sets \[
E_{-}=\left\{ x\in Q:f(x)\le k-1\right\} ,\, G=\left\{ x\in Q:f(x)=k\right\} \mbox{ and }E_{+}=\left\{ x\in Q:f(x)\ge k+1\right\} .\]
By properties of the non increasing rearrangement, or, more specifically,
in view of (\ref{eq:intval}), we obtain that the $\lambda$ measures
of these three sets are equal, respectively, to \[
\left|\left\{ t\in(0,\lambda(Q)):f^{*}(t)\le k-1\right\} \right|\mbox{ and }\left|\left\{ t\in(0,\lambda(Q)):f^{*}(t)=k\right\} \right|\]
 and $\left|\left\{ t\in(0,\lambda(Q)):f^{*}(t)\ge k+1\right\} \right|$.
Consequently, using (\ref{eq:fromp}) and then (\ref{eq:bomp}), we
see that $\lambda(G)\le(1-\sigma)(b-a)$. Then (\ref{eq:jomp}) and
(\ref{eq:klomp}) give us that $\lambda(E_{-})\ge\frac{\sigma}{2}(b-a)$
and $\lambda(E_{+})\ge\frac{\sigma}{2}(b-a)$. Therefore we have \begin{equation}
\min\left\{ \lambda(E_{-}),\lambda(E_{+})\right\} \ge\frac{\sigma}{2(1-\sigma)}\lambda(G)\,.\label{eq:utmdtv}\end{equation}
 It is a routine matter to check that the condition (\ref{eq:Nams})
is equivalent to \begin{equation}
\tau<\frac{\sigma}{2(1-\sigma)}\le\frac{1}{2}\,.\label{eq:sywc}\end{equation}

Let us pause for a moment to point out that we have finally reached
the \textit{only} step of the proof which needs a non trivial {}``geometric''
input, namely the fact that $(\tau,s)$ is a John-Strömberg pair for
$\mathcal{E}$.

To know that we have this fact in the particular case that appears
in the formulation of Theorem \ref{thm:UltraNewMain} we need to apply
our {}``geometric'' Theorem \ref{thm:maingt}. To know that we have
this fact for other particular collections $\mathcal{E}$, or with
better values of the constant $s$, we would need an affirmative answer
to Question A, or to some other question. In the present theorem we
have simply {}``axiomatized the 'geometric' problem away'' by invoking
a convenient definition.

Now let us resume our formal proof: The estimates (\ref{eq:sywc})
and (\ref{eq:utmdtv}) give us the inequality (\ref{eq:isms}) which
appears in Definition \ref{def:jsconst}. Since we have required that
$(\tau,s)$ is a John-Strömberg pair for $\mathcal{E}$\,, this guarantees
that there exists a set $W\in\mathcal{E}$ for which $W\subset Q$
and \begin{equation}
\min\left\{ \lambda(E_{+}\cap W),\lambda(E_{-}\cap W)\right\} \ge s\lambda(W)\,.\label{eq:smdtv}\end{equation}

Let us now, analogously to what we did in the easy case (iii) above,
define the functions $\varphi:\mathbb{R}\to\mathbb{R}$ and $\psi:\mathbb{R}\to\mathbb{R}$
by $\psi(t)=\max\left\{ t,k-1\right\} $ and then $\varphi(t)=\min\left\{ k+1,\psi(t)\right\} $.
Here again it is obvious that $\psi$ and therefore also $\varphi$
are both 1-Lipschitz functions. Therefore, again by Lemma \ref{lem:JLIP}
and by our hypotheses on $f$, we have that \begin{equation}
\left\Vert \varphi\circ f\right\Vert _{BMO(Q,\mathcal{E})}^{(\mathbf{J},s)}\le\left\Vert f\right\Vert _{BMO(Q,\mathcal{E})}^{(\mathbf{J},s)}\le1/2\,.\label{eq:ityb}\end{equation}
 Note also that the function $\varphi\circ f$ takes precisely three
values, namely $k-1$, $k$ and $k+1$.

Our final step will be to show that the set $W$ which satisfies (\ref{eq:smdtv})
must also satisfy \begin{equation}
\mathbf{J}(\varphi\circ f,W,s)\ge1\,.\label{eq:rc}\end{equation}
 This will contradict (\ref{eq:ityb}) and so show that the assumption
(\ref{eq:wa}) must be false, and thus will suffice to complete the
proof of Theorem r\ref{thm:Nultra}.

To simplify the notation, let us set $g=\varphi\circ f$. As already
remarked above, this non negative function takes only the three values
$k-1$, $k$, and $k+1$. More precisely, when we restrict $g$ to
the cube $W$, it takes these three values, respectively, on the sets
$E_{-}\cap W$, $G\cap W$ and $E_{+}\cap W$, whose union is $W$.

Thus the restriction of $(g\chi_{W})^{*}$ to the interval $\left(0,\lambda(W)\right)$
is given by \[
(g\chi_{W})^{*}(t)=\left\{ \begin{array}{ccc}
k+1 & , & 0<t<a\\
k & , & a\le t<b\\
k-1 & , & b\le t<\lambda(W)\end{array}\right.\]
 where $a=\lambda(E_{+}\cap W)$ , $b-a=\lambda(G\cap W)$ and $\lambda(W)-b=\lambda(E_{-}\cap W)$.

Let $I=\left[u,u+(1-s)\lambda(W)\right]$ be an arbitrary closed interval
of length $(1-s)\lambda(W)$ which is contained in $\left(0,\lambda(W)\right)$.
Since the left interval $(0,a)$ has length not less than $s\lambda(W)$
we conclude that the left endpoint of $I$ must line in $(0,a)$.
Similarly, since the right interval $[b,\lambda(W))$ also has length
not less than $s\lambda(W)$, the right endpoint of must lie in $[b,\lambda(W))$.
We deduce that \[
\left(g\chi_{W}\right)^{*}(u)-\left(g\chi_{W}\right)^{*}\left(u+(1-s)\lambda(Q)\right)=(k+1)-(k-1)=2\,.\]
Therefore, by Proposition \ref{pro:prop}, we have that $\mathbf{J}(g,W,s)=1$
which establishes (\ref{eq:rc}). As already explained above, this
suffices to complete the proof of Theorem \ref{thm:UltraNewMain}.
$\qed$

We conclude this section by stating and proving the auxiliary result
alluded to above in Remark \ref{rem:conversethm}, which is a sort
of converse to Theorem \ref{thm:Nultra} but which will not be needed
for other purposes here. Note that here the connection (\ref{eq:Nams})
between $\tau$ and $\sigma$ of Theorem \ref{thm:Nultra} has to
be replaced by the different connection (\ref{eq:NNams}). Also the
collection $\mathcal{E}$ appearing in (\ref{eq:Nyphwu}) has to be
replaced by the closely related collection $\mathcal{E}(Q)$.
\begin{thm}
\label{thm:Ninverse}Let $\mathcal{E}$ be a collection of admissible
subsets of $\mathbb{R}^{d}$. Let $s$ and $\sigma$ be two given
numbers in $(0,1/2)$ and let $\tau$ be any number in $(0,1)$ satisfying\begin{equation}
0<\sigma<\frac{\tau}{1+2\tau}\,.\label{eq:NNams}\end{equation}
Suppose that, for each $Q\in\mathcal{E}$, the inequality \begin{equation}
\left\Vert f^{*}\right\Vert _{BMO(I,\mathcal{Q}(I))}^{(\mathbf{J},\sigma)}\le\left\Vert f\right\Vert _{BMO(Q,\mathcal{E}(Q))}^{(\mathbf{J},s)}\label{eq:NNyphwu}\end{equation}
holds for the interval $I=\left(0,\lambda(Q)\right)$ and for every
measurable function $f:Q\to\mathbb{R}$ which assumes only the three
values $0$, $1$ and $2$. Here $\mathcal{E}(Q)$ denotes the collection
of all those sets in $\mathcal{E}$ which are contained in $Q$.

Then $(\tau,s)$ is a John-Strömberg pair for $\mathcal{E}$\,.
\end{thm}
\noindent \textit{Proof.} Let $Q$ be an arbitrary set in $\mathcal{E}$.
Suppose that $Q$ is the disjoint union of measurable sets $E_{+}$,
$E_{-}$ and $G$. Suppose that \begin{equation}
\min\left\{ \lambda(E_{+}),\lambda(E_{-})\right\} >\tau\lambda(G)\,.\label{eq:ctau}\end{equation}
 We have to show that there exists some $W\in\mathcal{E}$ such that
$W\subset Q$ and \begin{equation}
\min\left\{ \lambda(E_{+}\cap W),\lambda(E_{-}\cap W)\right\} \ge s\lambda(W)\,.\label{eq:hrp}\end{equation}
 We may also suppose, without loss of generality, that \begin{equation}
\lambda(E_{-})\le\lambda(E_{+})\,,\label{eq:mlp}\end{equation}
since, if not, we can simply interchange the roles of $E_{+}$ and
$E_{-}$. Let $f:Q\to[0,\infty)$ be the measurable function $f=2\chi_{E_{-}}+\chi_{G}$.
Then \[
f^{*}=2\chi_{(0,\lambda(E_{-}))}+\chi_{[\lambda(E_{-}),\lambda(E_{-})+\lambda(G))}\,.\]
In view of (\ref{eq:mlp}), the interval $I_{0}:=(0,2\lambda(E_{-})+\lambda(G))$
is contained in $\left(0,\lambda(Q)\right)$. In view of (\ref{eq:ctau}),
\[
\left|I_{0}\right|=2\lambda(E_{-})+\lambda(G)>(2\tau+1)\lambda(G)\]
and therefore $\lambda(G)<\frac{1}{2\tau+1}\left|I_{0}\right|$\,.
Since the inequality (\ref{eq:NNams}) implies that $\frac{1}{2\tau+1}<1-2\sigma$,
we deduce that \begin{equation}
\lambda(G)<(1-2\sigma)\left|I_{0}\right|\,.\label{eq:omq}\end{equation}
On the interval $I_{0}$, the function $f^{*}$ takes the value $1$
on a subinterval $I_{G}$ of length $\lambda(G)$ located centrally
in $I_{0}$ and it takes the values $2$ and $0$, respectively, on
two intervals $I_{-}$ and $I_{+}$, both of length $\lambda(E_{-})$,
located respectively on the left and right sides of $I_{G}$. In view
of (\ref{eq:omq}), each of these two intervals has length greater
than $\sigma\left|I_{0}\right|$. Now let \[
I_{1}=\left[u,u+(1-\sigma)\left|I_{0}\right|\right]\]
be an arbitrary closed interval of length $(1-\sigma)\left|I_{0}\right|$
contained in $I_{0}$. This means that $0<u<\sigma\left|I_{0}\right|$.
The fact that $u<\sigma\left|I_{0}\right|$ ensures that the left
endpoint of $I_{0}$ must lie in the interior of $I_{-}$. The fact
that $u>0$ ensures that the right endpoint of $I_{1}$ must lie in
the interior of $I_{+}$. It follows that \[
f^{*}(u)-f^{*}(u+(1-\sigma)\left|I_{0}\right|)=2\]
for each such interval $I_{1}$. Consequently, using (\ref{eq:ffnn}),
we deduce that $\left\Vert f^{*}\right\Vert _{BMO(I,\mathcal{Q}(I))}^{(\mathbf{J},\sigma)}=1$.
The hypotheses of the current theorem include the assumption that
our function $f$ satisfies (\ref{eq:NNyphwu}). Note also that the
formula (\ref{eq:ffnn}) and the fact that the ranges of $f$ and
$f^{*}$ are both contained in $[0,2]$ guarantee that neither $\left\Vert f\right\Vert _{BMO(Q,\mathcal{E}(Q))}^{(\mathbf{J},s)}$
nor $\left\Vert f^{*}\right\Vert _{BMO(I,\mathcal{Q}(I))}^{(\mathbf{J},\sigma)}$
can be greater than $1$. We deduce that \begin{equation}
\left\Vert f\right\Vert _{BMO(Q,\mathcal{E}(Q))}^{(\mathbf{J},s)}=1\,.\label{eq:mqtpat}\end{equation}
We now use (\ref{eq:mqtpat}) to deduce the existence of a set $W\in\mathcal{E}(Q)$
such that \begin{equation}
\mathbf{J}\left(f,W,s\right)=1\,.\label{eq:otpat}\end{equation}
 (When calculating $\left\Vert f\right\Vert _{BMO(Q,\mathcal{E}(Q))}^{(\mathbf{J},s)}$
we take the supremum of a set of numbers which is a subset of $\left\{ 0,\frac{1}{2},1\right\} $.)

To complete the proof we have to show that $W$ is a set with the
required property (\ref{eq:hrp}). Let us suppose that this is not
the case and show that this leads to a contradiction. Our supposition
means that at least one of the two sets $(E_{+}\cup G)\cap W$ and
$(E_{-}\cup G)\cap W$ must have $\lambda$ measure strictly greater
than $(1-s)\lambda(W)$.

Let $h:(0,\lambda(W))\to[0,\infty)$ be the restriction to $\left(0,\lambda(W)\right)$
of the non increasing rearrangement of $f\chi_{W}$. Then the statement
of the last sentence of the preceding paragraph means that there exists
an subinterval $I_{0}$ of the interval $(0,\lambda(W))$ which has
length strictly greater than $(1-s)\lambda(W)$ on which $h$ takes
only two values, and that these two values are either $0$ and $1$
or $1$ and $2$. This enables us to find a closed interval $\left[u,u+(1-s)\lambda(W)\right]$
contained in the interior of $I_{0}$ for which \[
h(u)-h\left(u+(1-s)\lambda(W)\right)\le1\]
and therefore $\mathbf{J}\left(f,W,s\right)\le1/2$. This contradicts
(\ref{eq:otpat}) and therefore completes the proof of the theorem.
$\qed$

\section{\label{sec:Putting}Putting all the pieces together.}

Now at last we can combine our results from previous sections to obtain
our versions of the John-Nirenberg and John-Strömberg inequalities.
The following theorem does this. It also explicitly and immediately
shows (keeping in mind Remark \ref{rem:tauIsLittle}) the consequence
of an affirmative answer to Question A, thus proving what we claimed
at the very beginning of this paper.

The hypotheses imposed on $\mathcal{E}$, $Q$, $\tau$ and $s$ in
this theorem are exactly those which were imposed in Theorem \ref{thm:Nultra}\,.
But here the parameter $\sigma$ does not need to be explicitly mentioned.
\begin{thm}
\label{thm:mainjs}Let $\mathcal{E}$ be a collection of admissible
subset of $\mathbb{R}^{d}$. Let $Q$ be a set in $\mathcal{E}$ which
contains all other sets of $\mathcal{E}$. Let $\tau\in(0,1/2)$ and
$s\in(0,1/2)$ be such that $(\tau,s)$ is a John-Strömberg pair for
$\mathcal{E}$\,. Then, for every constant $r$ in the range $1\le r\le1/2\tau$,
the inequalities\begin{equation}
\lambda\left(\left\{ x\in Q:\left|f(x)-m\right|\ge\alpha\right\} \right)\le\max\left\{ r,2\sqrt{r}\right\} \cdot\lambda(Q)\cdot\exp\left(-\frac{\alpha\log r}{8\left\Vert f\right\Vert _{BMO(Q,\mathcal{E})}^{(\mathbf{J},s)}}\right)\label{eq:Neimed}\end{equation}
 and\begin{equation}
\lambda\left(\left\{ x\in Q:\left|f(x)-m\right|\ge\alpha\right\} \right)\le\max\left\{ r,2\sqrt{r}\right\} \cdot\lambda(Q)\cdot\exp\left(-\frac{\alpha s\log r}{8\left\Vert f\right\Vert _{BMO(Q,\mathcal{E})}}\right)\label{eq:Nbeimed}\end{equation}
 hold for every $\alpha\ge0$, every measurable $f:Q\to\mathbb{R}$\,,
and every median $m$ of $f$ on $Q$.\end{thm}
\begin{rem}
\label{rem:choosingr}The inequality (\ref{eq:newurp}) which appears
at the very beginning of this paper is of course simply (\ref{eq:Nbeimed})
with $r=1/2\tau$. This is in some sense the most {}``pertinent''
value of $r$ to choose since it gives the best control of the left
hand side of (\ref{eq:Nbeimed}) when we consider large values of
$\alpha$ 
\end{rem}
\noindent \textit{Proof of Theorem \ref{thm:mainjs}.} Let us prepare
some ingredients which will later enable us to apply Theorem \ref{thm:temjemt-1}.
In our application the sets $D$ and $E$ which appear in the statement
of Theorem \ref{thm:temjemt-1} will both be taken to equal the set
$Q$ specified in the formulation here of Theorem \ref{thm:mainjs}.
Let $f:Q\to\mathbb{N}\cup\left\{ 0\right\} $ be an arbitrary function
in the class $\Phi$ which is defined in Theorem \ref{thm:temjemt-1}.
Let $\sigma$ be a number satisfying the condition (\ref{eq:Nams})
which is imposed in Theorem \ref{thm:Nultra}. Then, since $f$ and
$\sigma$ and $s$ satisfy the hypotheses of Theorem \ref{thm:Nultra},
we can deduce from Theorem \ref{thm:Nultra} that \begin{equation}
\left\Vert f^{*}\right\Vert _{BMO(I,\mathcal{Q}(I))}^{(\mathbf{J},\sigma)}\le\left\Vert f\right\Vert _{BMO(Q,\mathcal{E})}^{(\mathbf{J},s)}\,,\label{eq:flimp}\end{equation}
 where $I$ is the interval $\left(0,\lambda(Q)\right)$. The fact
that $f^{*}$ is right continuous and non increasing on $I$ ensures
that, for every $\alpha\ge0$, it satisfies the inequality corresponding
to (\ref{eq:odms}) which here takes the form \begin{equation}
\left|\left\{ t\in I:f^{*}(t)-f^{*}\left(\lambda(Q)/2\right)\ge\alpha\right\} \right|\le\frac{1-\sigma}{2\sigma}\cdot\left|I\right|\cdot\exp\left(-\frac{\alpha\log\left(\frac{1}{\sigma}-1\right)}{2\left\Vert f^{*}\right\Vert _{BMO\left(I,\mathcal{Q}(I)\right)}^{(\mathbf{J},\sigma)}}\right)\,.\label{eq:jerp}\end{equation}

As an element of $\Phi,$ the function $f$ must satisfy \begin{equation}
\lambda\left(\left\{ x\in Q:f(x)>0\right\} \right)\le\frac{1}{2}\lambda(Q)\,.\label{eq:vergp}\end{equation}
 This implies that $f^{*}(t)=0$ for all $t>\lambda(Q)/2$ and therefore
also that $f^{*}\left(\lambda(Q)/2\right)=0$. (Use properties (i)
and (ii) of Section \ref{sec:SignedRearrangements}.)

Now, for every $\alpha>0$ we can apply (\ref{eq:df2}) to the left
hand side of (\ref{eq:jerp}) and use (\ref{eq:flimp}) to bound the
right hand side of (\ref{eq:jerp}) from above. This gives us that
\begin{equation}
\lambda\left(\left\{ x\in Q:f(x)\ge\alpha\right\} \right)\le\frac{1-\sigma}{2\sigma}\cdot\lambda(Q)\cdot\exp\left(-\frac{\alpha\log\left(\frac{1}{\sigma}-1\right)}{2\left\Vert f\right\Vert _{BMO(Q,\mathcal{E})}^{(\mathbf{J},s)}}\right)\,.\label{eq:oink}\end{equation}

We deduce, using (\ref{eq:oink}) when $\alpha>0$, and using (\ref{eq:vergp}),
together with the fact that $\frac{1-\sigma}{2\sigma}\ge\frac{1}{2}$,
when $\alpha=0$, that \[
\lambda\left(\left\{ x\in Q:f(x)>\alpha\right\} \right)\le\frac{1-\sigma}{2\sigma}\cdot\lambda(Q)\cdot\exp\left(-\frac{\alpha\log\left(\frac{1}{\sigma}-1\right)}{2\left\Vert f\right\Vert _{BMO(Q,\mathcal{E})}^{(\mathbf{J},s)}}\right)\]
 for all $\alpha\ge0$. This last inequality corresponds to the inequality
(\ref{eq:brinz}) of Theorem \ref{thm:temjemt-1}, and the fact that
it holds for every $f\in\Phi$ is exactly what we need to justify
applying Theorem \ref{thm:temjemt-1}. Since here the constants $C$
and $c$ of (\ref{eq:brinz}) are respectively \[
\frac{1-\sigma}{2\sigma}\mbox{ and }\frac{\log\left(\frac{1}{\sigma}-1\right)}{2},\]
the formula (\ref{eq:inxz}) furnished by Theorem \ref{thm:temjemt-1}
takes the form

\begin{equation}
\lambda\left(\left\{ x\in Q:\left|f(x)-m\right|\ge\alpha\right\} \right)\le2\max\left\{ \frac{1-\sigma}{2\sigma},\sqrt{\frac{1-\sigma}{\sigma}}\right\} \cdot\lambda(Q)\cdot\exp\left(-\frac{\alpha\log\left(\frac{1}{\sigma}-1\right)}{8\left\Vert f\right\Vert _{BMO(Q,\mathcal{E})}^{(\mathbf{J},s)}}\right)\,.\label{eq:eimed}\end{equation}
 Let us put $r=\frac{1}{\sigma}-1=\frac{1-\sigma}{\sigma}$. Since
$\sigma$ satisfies (\ref{eq:Nams}) and since, as already observed
during the proof of Theorem \ref{thm:Nultra}, (\ref{eq:Nams}) is
equivalent to (\ref{eq:sywc}), it follows that $\tau<1/2r\le1/2$
and so $1\le r<1/2\tau$. Thus (\ref{eq:eimed}) gives us (\ref{eq:Neimed})
for all $r\in[1,1/2\tau)$, and therefore also, by continuity, at
the endpoint $r=1/2\tau$. Finally (\ref{eq:Nbeimed}) follows immediately,
in view of (\ref{eq:talpt}). $\qed$ 
\begin{rem}
\label{rem:specialrectangles} Let us consider Theorem \ref{thm:mainjs}
in the particular case where $Q$ is a special rectangle and $\mathcal{E}=\mathcal{W}(Q)$.
Then we know (cf.~Corollary \ref{cor:cmaingt} or Example \ref{exa:emaingt})
that we can take $\tau=\sqrt{2}-1$ and $s=(3-2\sqrt{2})/2\approx0.0857864$.
So the parameter $r$ can range between $1$ and $1/(2\sqrt{2}-2)\approx1.20711$.
Consequently, the right hand side of (\ref{eq:Nbeimed}) can be, for
example, \[
2/\sqrt{2\sqrt{2}-2}\cdot\lambda(Q)\cdot\exp\left(-\frac{\alpha(3-2\sqrt{2})\log\left(1/(2\sqrt{2}-2)\right)}{16\left\Vert f\right\Vert _{BMO}^{\prime}}\right)\]
where $\left\Vert f\right\Vert _{BMO}^{\prime}$ is the seminorm defined
in (\ref{eq:WikDef}). This expression is approximately equal to $2.197\cdot\lambda(Q)\cdot\exp\left(-\frac{0.002\alpha}{\left\Vert f\right\Vert _{BMO}^{\prime}}\right)$\,.
Wik obtains a smaller expression, approximately equal to $2\lambda(Q)\exp\left(-\frac{0.043\alpha}{\left\Vert f\right\Vert _{BMO}^{\prime}}\right)$\,. 
\end{rem}

\section{\label{sec:attempts}Towards an answer to Question A}

The reader who has already traversed all the previous sections of
this paper to get to here is presumably convinced by now that it is
worth trying to answer Question A. So let us spend this section trying
to offer some help towards that goal. We will formulate three new
questions. The first and second of them seem to be somewhat easier
to answer than Question A, and we shall see that they are each essentially
equivalent to Question A. The third question is of a somewhat different
nature, and not equivalent. But we shall indicate why we consider
it also to be well worth considering.

For each $x\in\mathbb{R}^{d}$ and $r>0$ we recall the standard notation
$Q(x,r)$ for the set \begin{equation}
Q\left(x,r\right)=\left\{ y\in\mathbb{R}^{d}:\left\Vert y-x\right\Vert _{\ell_{d}^{\infty}}\le r\right\} \,.\label{eq:defqxr}\end{equation}
This is of course the closed cube in $\mathbb{R}^{d}$ centred at
the point $x$ and having side length $2r$\,. We will need to use
two simple properties of such cubes, which we state in the following
lemma:
\begin{lem}
\label{lem:SimplePropsCubes}Let $\left\{ x_{n}\right\} _{n\in\mathbb{N}}$
be a sequence of points in $\mathbb{R}^{d}$ and let $\left\{ r_{n}\right\} _{n\in\mathbb{N}}$
be a sequence of positive numbers which converge, respectively to
the point $x_{*}\in\mathbb{R}^{d}$ and the positive number $r_{*}$. 

(i) If all of the cubes $Q(x_{n},r_{n})$ are contained in some closed
set $Q$ then $Q(x_{*},r_{*})$ is also contained in $Q$.

(ii) Let $A$ be a fixed measurable subset of $\mathbb{R}^{d}$. Then
\begin{equation}
\lim_{n\to\infty}\lambda(A\cap Q(x_{n},r_{n}))=\lambda(A\cap Q(x_{*},r_{*}))\,.\label{eq:ttz}\end{equation}

\end{lem}
The proof of this lemma is an easy exercise. In this preliminary version
of our paper we include it as Appendix \ref{sub:ProveSimpleProps}.

\bigskip{}

\subsection{Reducing the description of John-Str\"omberg pairs to the case of
sets which are finite unions of cubes.\phantom{AAAAAAAAAAAAAAAAAAAAAAAAAAAAAAAAAAAAAAAAAAAAAAAAAAAAAAAAAAAAAAAAAAAAAAAAAAAAAAAAAAAAAAAAAAAAAAAAAAAAAAAAAAAAAAAAAAAAAAAAAAAAAAAAAAAAAAAAAAAAAAAAAAAAAAAAAAAAAA}}

\medskip{}

\phantom{.}

The following theorem shows that, in order to obtain an affirmative
answer to Question A, or merely to determine any pairs $(\tau,s)$
which are John-Str\"omberg pairs for the collection of all cubes
in $\mathbb{R}^{d}$, it suffices to consider only those sets $E_{+}$
and $E_{-}$ which are of a comparatively simple form. This also suggests
that one might choose to reformulate Question A and the question of
finding John-Str\"omberg pairs for cubes in $\mathbb{R}^{d}$ in
terms which are more in the realm of combinatorics, counting lattice
points etc.
\begin{thm}
\label{thm:onlyneedcubes}Suppose that the numbers $\tau\in(0,1/2)$
and $s>0$ have the following property:

Whenever $F_{+}$ and $F_{-}$ are disjoint subsets of the closed
unit cube $Q=[0,1]^{d}$ in $\mathbb{R}^{d}$ such that \textbf{\[
\min\left\{ \lambda(F_{+}),\lambda(F_{-})\right\} >\tau\lambda(Q\setminus F_{+}\setminus F_{-})\,,\]
}and also each of the sets $F_{+}$ and $F_{-}$ is the union of finitely
many dyadic cubes, then there exists some cube $W$ contained in $Q$
for which\textbf{\[
\min\left\{ \lambda(W\cap F_{+}),\lambda(W\cap F_{-})\right\} \ge s\lambda(W)\,.\]
}

Then $(\tau,s)$ is a John-Str\"omberg pair for the collection of
cubes in $\mathbb{R}^{d}$.
\end{thm}
\noindent \textit{Proof.} We will use the notation {}``$\flat$''
to stand for a subscript that is either {}``$+$'' or {}``$-$''.
More precisely, whenever we write a formula where $\flat$ appears
as a subscript in one or more places, this expresses the fact that
the formula holds in both of the cases: 

(i) whenever that subscript is replaced throughout by the subscript
$+$, and also

(ii) whenever it is replaced throughout by the subscript $-$.

We start by writing down two simple formulae which will be useful
later:

Suppose that $V_{+}$ and $V_{-}$ are two disjoint measurable subsets
of $Q$ and that $U_{+}$ and $U_{-}$ are measurable subsets, respectively,
of $V_{+}$ and of $V_{-}$. Then, obviously, \begin{equation}
\lambda(V_{\flat})=\lambda(U_{\flat})+\lambda\left(V_{\flat}\setminus U_{\flat}\right)\label{eq:os1}\end{equation}
 and so, consequently, \begin{eqnarray}
\lambda\left(Q\setminus U_{+}\setminus U_{-}\right) & = & \lambda(Q)-\lambda(U_{+})-\lambda(U_{-})\nonumber \\
 & = & \lambda(Q)-\lambda(V_{+})-\lambda(V_{-})+\lambda(V_{+}\setminus U_{+})+\lambda(V_{-}\setminus U_{-})\nonumber \\
 & = & \lambda\left(Q\setminus V_{+}\setminus V_{-}\right)+\lambda(V_{+}\setminus U_{+})+\lambda(V_{-}\setminus U_{-})\,.\label{eq:os2}\end{eqnarray}

After this preparation, let us suppose that $E_{+}$ and $E_{-}$
are disjoint arbitrary measurable subsets of the cube $Q$ which satisfy\textbf{\textit{\[
\min\left\{ \lambda(E_{+}),\lambda(E_{-})\right\} >\tau\lambda(Q\setminus E_{+}\setminus E_{-})\,.\]
}}In order to prove the theorem, we have to find a cube $W$ contained
in $Q$ for which \textbf{\textit{\begin{equation}
\min\left\{ \lambda(W\cap E_{+}),\lambda(W\cap E_{-})\right\} \ge s\lambda(W)\,.\label{eq:wwnn}\end{equation}
}}The obvious fact that we need only consider the case where $Q$
is the unit cube $Q=[0,1]^{d}$ was already pointed out in Remark
\ref{rem:unitcube}.We can of course assume without loss of generality
that $E_{+}$ and $E_{-}$ are both contained in $Q^{\circ}$, the
interior of $Q$. 

Since Lebesgue measure is inner regular, there exist compact sets
$H_{+}$ and $H_{-}$ contained respectively in $E_{+}$ and $E_{-}$
such that $\lambda$$(E_{+}\setminus H_{+})$ and $\lambda$$(E_{-}\setminus H_{-})$
are both sufficiently small to guarantee (of course via formulae like
(\ref{eq:os1}) and (\ref{eq:os2})) that\textbf{\textit{\[
\min\left\{ \lambda(H_{+}),\lambda(H_{-})\right\} >\tau\lambda(Q\setminus H_{+}\setminus H_{-})\,.\]
}}

We let $\varepsilon_{0}$ be a positive number which is chosen sufficiently
small so that it satisfies \begin{equation}
2\varepsilon_{0}+4\tau\varepsilon_{0}<\min\left\{ \lambda(H_{+}),\lambda(H_{-})\right\} -\tau\lambda(Q\setminus H_{+}\setminus H_{-})\,.\label{eq:cdz}\end{equation}
Let $\delta=\mathrm{dist}\left(H_{+},H_{-}\right)$.This is of course
a positive number, since $H_{+}$ and $H_{-}$ are disjoint and compact.
We shall use $\delta$ to obtain two disjoint open sets $\Omega_{+}$
and $\Omega_{-}$ contained in $Q^{\circ}$, such that $H_{\flat}\subset\Omega_{\flat}$,
and \begin{equation}
\mathrm{dist}(\Omega_{+},\Omega_{-})\ge\frac{\delta}{2}\,.\label{eq:ompmd}\end{equation}
Initially we can choose $\Omega_{\flat}$ to be the set $\Omega_{\flat}=Q^{\circ}\cap\bigcup_{x\in H_{\flat}}(x+B)$,
where $B$ is the open ball of radius $\delta/4$ centred at the origin,
and this indeed will guarantee that (\ref{eq:ompmd}) holds. But then,
furthermore, since Lebesgue measure is outer regular, we can, by replacing
$\Omega_{\flat}$ if necessary by its intersection with some other
open set containing $H_{\flat}$, assume also that \begin{equation}
\lambda(\Omega_{\flat}\setminus H_{\flat})<\varepsilon_{0}\label{eq:pinq}\end{equation}
and so (cf.~(\ref{eq:os1}) \begin{equation}
\lambda(H_{\flat})\le\lambda\left(\Omega_{\flat}\right)<\lambda(H_{\flat})+\varepsilon_{0}\,.\label{eq:tema}\end{equation}
Now we use the fact (see e.g.~\cite{WheedenZygmund} Theorem (1.11)
p.~8) that every open set $\Omega$ in $\mathbb{R}^{d}$ is the union
of some sequence $\left\{ D_{n}\right\} _{n\in\mathbb{N}}$ of non
overlapping closed dyadic cubes. In particular we shall write $\Omega_{\flat}=\bigcup_{n\in\mathbb{N}}D_{\flat,n}$
for both choices of $\flat$. Obviously all of the dyadic cubes $D_{\flat,n}$
have to be contained in $Q$.

Let $N$ be a positive integer which is sufficiently large to ensure
that \begin{equation}
\sum_{n=N+1}^{\infty}\lambda(D_{\flat,n})<\varepsilon_{0}\label{eq:tevma}\end{equation}
for both choices of $\flat$. Of course there exists such an $N$
since \[
\sum_{n=1}^{\infty}\lambda(D_{\flat,n})=\lambda(\Omega_{\flat})\le\lambda(Q^{\circ})<\infty\,.\]
Then let $F_{\flat}$ be the compact set $F_{\flat}=\bigcup_{n=1}^{N}D_{\flat,n}$. 

We want to be able to apply the hypotheses of the theorem to the two
disjoint sets $F_{+}$ and $F_{-}$. So we need estimates from below
for the measures of each of these sets. But first we note that $F_{\flat}\subset\Omega_{\flat}$
and, by (\ref{eq:tevma}), \begin{equation}
\lambda\left(\Omega_{\flat}\setminus F_{\flat}\right)=\lambda\left(\bigcup_{n>N}D_{\flat,n}\right)<\varepsilon_{0}\,.\label{eq:inq}\end{equation}
It follows (cf.~(\ref{eq:os1})) that \[
\lambda(F_{\flat})\le\lambda(\Omega_{\flat})<\lambda(F_{\flat})+\varepsilon_{0}\,.\]
This last estimate, together with (\ref{eq:tema}), gives us that
\begin{equation}
\left|\lambda(F_{\flat})-\lambda(H_{\flat})\right|<2\varepsilon_{0}\,.\label{eq:twinq}\end{equation}

Now we substitute $V_{\flat}=\Omega_{\flat}$ and $U_{\flat}=H_{\flat}$
in (\ref{eq:os2}) and then apply (\ref{eq:pinq}) to obtain that
\[
\left|\lambda(Q\setminus H_{+}\setminus H_{-})-\lambda(Q\setminus\Omega_{+}\setminus\Omega_{-})\right|<2\varepsilon_{0}\,.\]
Then, analogously, we substitute $V_{\flat}=\Omega_{\flat}$ and $U_{\flat}=F_{\flat}$
in (\ref{eq:os2}) and then apply (\ref{eq:inq}) to obtain that \[
\left|\lambda(Q\setminus F_{+}\setminus F_{-})-\lambda(Q\setminus\Omega_{+}\setminus\Omega_{-})\right|<2\varepsilon_{0}\,.\]
The preceding two inequalities imply that \begin{equation}
\left|\lambda(Q\setminus F_{+}\setminus F_{-})-\lambda(Q\setminus H_{+}\setminus H_{-})\right|<4\varepsilon_{0}\,.\label{eq:brompz}\end{equation}

We can now deduce our required estimates from below for $\lambda(F_{\flat})$.
In the following calculation we shall use (\ref{eq:twinq}) in the
first line, and then (\ref{eq:cdz}) in the third line, and then (\ref{eq:brompz})
in the fourth line.

\begin{eqnarray*}
\lambda(F_{\flat}) & \ge & \lambda(H_{\flat})-2\varepsilon_{0}\\
 & \ge & \min\left\{ \lambda(H_{+}),\lambda(H_{-})\right\} -2\varepsilon_{0}\\
 & > & \tau\lambda\left(Q\setminus H_{+}\setminus H_{-}\right)+(2\varepsilon_{0}+4\tau\varepsilon_{0})-2\varepsilon_{0}\\
 & > & \tau\lambda\left(Q\setminus F_{+}\setminus F_{-}\right)+(2\varepsilon_{0}+4\tau\varepsilon_{0})-2\varepsilon_{0}-4\tau\varepsilon_{0}\\
 & = & \tau\lambda\left(Q\setminus F_{+}\setminus F_{-}\right)\,.\end{eqnarray*}

Since the sets $F_{+}$ and $F_{-}$ are disjoint and are each finite
unions of dyadic intervals, the estimate which we have just obtained
for $\lambda(F_{+})$ and $\lambda(F_{-})$ is exactly the one which
we require to apply the hypothesis of our theorem. That hypothesis
implies that there exists a cube, which we will denote by $W(\varepsilon_{0})$
contained in $Q$, for which\textbf{\textit{\begin{equation}
\min\left\{ \lambda(W(\varepsilon_{0})\cap F_{+}),\lambda(W(\varepsilon_{0})\cap F_{-})\right\} \ge s\lambda(W(\varepsilon_{0}))>0\,.\label{eq:op131}\end{equation}
}}

We deduce from (\ref{eq:ompmd}) that $\mathrm{dist}(F_{+},F_{-})\ge\frac{\delta}{2}$.
This in turn implies that the cube $W(\varepsilon_{0})$ must satisfy
\begin{equation}
\mathrm{diam}W(\varepsilon_{0})\ge\frac{\delta}{2}\label{eq:diamb}\end{equation}
since it contains points of $F_{+}$ and also points of $F_{-}$.

It follows from (\ref{eq:pinq}) that \[
\lambda\left((W(\varepsilon_{0})\cap\Omega_{\flat})\setminus(W(\varepsilon_{0})\cap H_{\flat})\right)<\varepsilon_{0}\]
 and, from (\ref{eq:inq}), that \[
\lambda\left((W(\varepsilon_{0})\cap\Omega_{\flat})\setminus(W(\varepsilon_{0})\cap F_{\flat})\right)<\varepsilon_{0}\,.\]
Therefore, for both choices of $\flat$,\begin{eqnarray*}
 &  & \left|\lambda(W(\varepsilon_{0})\cap F_{\varepsilon})-\lambda(W(\varepsilon_{0})\cap H_{\flat})\right|\\
 & \le & \left|\lambda(W(\varepsilon_{0})\cap F_{\varepsilon})-\lambda(W(\varepsilon_{0})\cap\Omega_{\flat})\right|+\left|\lambda(W(\varepsilon_{0})\cap\Omega_{\varepsilon})-\lambda(W(\varepsilon_{0})\cap H_{\flat})\right|\\
 & \le & 2\varepsilon_{0}\,.\end{eqnarray*}
Combining this inequality with (\ref{eq:op131}) gives us that

\begin{equation}
\min\left\{ \lambda(W(\varepsilon_{0})\cap H_{+}),\lambda(W(\varepsilon_{0})\cap H_{-})\right\} \ge s\lambda(W(\varepsilon_{0}))-2\varepsilon_{0}\,.\label{eq:miy}\end{equation}

Note that we fixed the sets $H_{+}$ and $H_{-}$ and therefore also
the number $\delta$, before we chose the number $\varepsilon_{0}$.
This means that, for each number $\varepsilon$ satisfying $0<\varepsilon\le\varepsilon_{0}$,
we can use exactly the same arguments as were used from (\ref{eq:cdz})
to (\ref{eq:miy}) to obtain a subcube $W\left(\varepsilon\right)$
of $Q$ which satisfies (\ref{eq:diamb}) and (\ref{eq:miy}), but
with $\varepsilon$ in place of $\varepsilon_{0}$. In fact we shall
do this for each $\varepsilon_{n}$ in a sequence $\left\{ \varepsilon_{n}\right\} _{n\in\mathbb{N}}$
of numbers in the interval $(0,\varepsilon_{0}]$ which converges
to $0$. Thus we obtain a sequence $\left\{ W(\varepsilon_{n})\right\} _{n\in\mathbb{N}}$
of subcubes of $Q$ which, for each $n$, satisfy the counterparts
of (\ref{eq:diamb}) and (\ref{eq:miy}), namely\begin{equation}
\mathrm{diam}W(\varepsilon_{n})\ge\frac{\delta}{2}\label{eq:slmi}\end{equation}
and \begin{equation}
\min\left\{ \lambda(W(\varepsilon_{n})\cap H_{+}),\lambda(W(\varepsilon_{n})\cap H_{-})\right\} \ge s\lambda(W(\varepsilon_{n}))-2\varepsilon_{n}\,.\label{eq:injb}\end{equation}

For each $n\in\mathbb{N}$ let $r_{n}>0$ and $x_{n}\in Q$ be the
half side length and centre, respectively, of $W(\varepsilon_{n})$.
I.e., as in (\ref{eq:defqxr}), we have $W(\varepsilon_{n})=Q(x_{n},r_{n})$.
By passing, if necessary, to a subsequence, we can assume that the
sequences $\left\{ r_{n}\right\} _{n\in\mathbb{N}}$ and $\left\{ x_{n}\right\} _{n\in\mathbb{N}}$
are both convergent, to limits $r_{*}\ge0$ and $x_{*}\in Q$ respectively.
The condition (\ref{eq:slmi}) ensures that in fact $r_{*}>0$. Using
part (i) of Lemma \ref{lem:SimplePropsCubes}, we see that the cube
$W:=Q(x_{*},r_{*})$ is contained in $Q$. Then, applying part (ii)
of the same lemma with $A$ chosen to be the set $H_{\flat}$, we
obtain that \[
\lim_{n\to\infty}\lambda(W(\varepsilon_{n})\cap H_{\flat})=\lambda(W\cap H_{\flat})\]
for both choices of the subscript $\flat$. Obviously we also have
\[
\lim_{n\to\infty}\lambda(W(\varepsilon_{n}))=\lim_{n\to\infty}r_{n}^{d}=r_{*}^{d}=\lambda(W)\,.\]
Thus a passage to the limit as $n$ tends to $\infty$ in (\ref{eq:injb}),
shows that the cube $W$ satisfies $\min\left\{ \lambda(W\cap H_{+}),\lambda(W\cap H_{-})\right\} \ge s\lambda(W)$\,.
Since $H_{\flat}\subset E_{\flat}$ this last inequality immediately
implies that $W$ satisfies (\ref{eq:wwnn}). This, together with
the fact that $W\subset Q$ and the fact that the sets $E_{+}$ and
$E_{-}$ were chosen arbitrarily, suffices to show that $(\tau,s)$
is a John-Str\"omberg pair for cubes in $\mathbb{R}^{d}$. $\qed$
\begin{rem}
\label{rem:equality}We have not bothered to write out the details,
but it seems very likely that techniques similar to those in the proof
of the preceding theorem might show that whenever $(\tau,s)$ is a
John-Str\"omberg pair for the collection of all cubes in $\mathbb{R}^{d}$
then it also has the following slightly stronger property: 

Whenever $E_{+}$ and $E_{-}$ are disjoint admissible subsets of
the cube $Q$ in $\mathbb{R}^{d}$ such that \textbf{\textit{\begin{equation}
\min\left\{ \lambda(E_{+}),\lambda(E_{-})\right\} \ge\tau\lambda(Q\setminus E_{+}\setminus E_{-})\,,\label{eq:gyiq}\end{equation}
}}then there exists some cube $W$ contained in $Q$ for which\textbf{\textit{\[
\min\left\{ \lambda(W\cap E_{+}),\lambda(W\cap E_{-})\right\} \ge s\lambda(W)\,.\]
}}In other words, it seems very likely that, even if equality holds
in (\ref{eq:gyiq}), this is still enough to imply the existence of
a cube $W$ with the same properties as before. 
\end{rem}

\subsection{\label{sub:diffspec}A reduction of Question A to a different special
case.\phantom{AAAAAAAAAAAAAAAAAAAAAAAAAAAAAAAAAAAAAAAAAAAAAAAAAAAAAAAAAAAAAAAAAAAAAAAAAAAAAAAAAAAAAAAAAAAAAAAAAAAAAAAAAAAAAAAAAAAAAAAAAAAAAAAAAAAAAAAAAAAAAAAAAAAAAAAAAAAAAA}}

\smallskip{}

\phantom{.}

\smallskip{}

Here is a variant of Question A. We will call it Question A$^{\prime}$.

\medskip{}

\textit{Do there exist two absolute constants $\tau'\in(0,1/2)$ and
$s>0$ which have the following property?}

\textbf{\textit{For every positive integer $d$ and for every closed
cube $Q$ in $\mathbb{R}^{d}$, whenever $F_{+}$ and $F_{-}$ are
two disjoint compact subsets of $Q$ which each have positive measure,
which are each the union of finitely many closed rectangles, and whose
$d$-dimensional Lebesgue measures satisfy \begin{equation}
\min\left\{ \lambda(F_{+}),\lambda(F_{-})\right\} \ge\tau'\lambda(Q\setminus F_{+}\setminus F_{-})\,,\label{eq:ehyp}\end{equation}
then there exists some cube $W$ contained in $Q$ for which\begin{equation}
\min\left\{ \lambda(W\cap F_{+}),\lambda(W\cap F_{-})\right\} \ge s\lambda(W)\,.\label{eq:econ}\end{equation}
}}

\medskip{}
 Here the terminology \textit{closed rectangle} means a subset of
$\mathbb{R}^{d}$ which is the cartesian product of $d$ bounded closed
intervals. Our reason for using rectangles rather than cubes here
will become apparent later (in Remark \ref{rem:splurg}). Note that
if we did not include the requirement that both of $F_{+}$ and $F_{-}$
have positive measure, then, in view of the possibility of taking
$F_{+}=Q$, the answer to the above question would be negative. 

It follows immediately from Theorem \ref{thm:onlyneedcubes}, that
an affirmative answer to Question A$^{\prime}$ for some constants
$\tau'$ and $s$ would imply an affirmative answer to Question A
for every choice of the constant $\tau$ satisfying $\tau\in(\tau',1/2)$
and for the same constant $s$. Thus we may transfer our attention
from Question A to Question A$^{\prime}$. By the same simple reasoning
as in Remark \ref{rem:sbth}, if Question A$^{\prime}$ has an affirmative
answer, then the constant relevant constant $s$ has to satisfy $s\le1/2$.

We shall show that we can reduce Question A$^{\prime}$ to yet another
question which is, in principle, easier to answer. However its formulation
is more technical and requires the following three definitions. (The
second of them is reminiscent of, but slightly different from a definition
used in the course of the proof of Theorem \ref{thm:maingt}.)
\begin{defn}
\label{def:toptimal}Let $Q$ be a closed cube in $\mathbb{R}^{d}$
and let $F_{+}$ and $F_{-}$ be two disjoint compact subsets of $Q$
which (as in Question A$^{\prime}$) are each unions of finitely many
closed rectangles. Let $\tau'$ be a number in $(0,1/2)$. Let $V$
be a closed subcube of $Q$. 

(i) We will say that $V$ is an \textit{exceptional subcube} if $\lambda(V\cap F_{+})$
and $\lambda(V\cap F_{-})$ are both strictly positive and $\lambda(V\setminus F_{+}\setminus F_{-})=0$. 

(ii) We will say that $V$ is a \textit{good $\tau'$-subcube} if
$\lambda(V\cap F_{+})$ and $\lambda(V\cap F_{-})$ are both strictly
positive and \[
\min\left\{ \lambda(V\cap F_{+}),\lambda(V\cap F_{-})\right\} \ge\tau'\lambda(V\setminus F_{+}\setminus F_{-})>0\,.\]

(iii) We will say that $V$ is $\tau'$\textit{-minimal} if it is
a good $\tau'$-subcube but every strictly smaller closed subcube
of $V$ is neither exceptional nor a good $\tau'$-subcube. \end{defn}
\begin{rem}
\label{rem:splurg}Of course, even though our chosen terminologies
do not explicitly express it, all three of the above notions depend
crucially on the choice of the cube $Q$ and of its subsets $F_{+}$
and $F_{-}$. When occasionally necessary, we can replace these terminologies
by the more explicit \textit{exceptional $(Q,F_{+},F_{-})$-subcube}
and \textit{good $(\tau',Q,F_{+},F_{-})$-subcube }and\textit{ $(\tau',Q,F_{+},F_{-})$-minimal
subcube. }In connection with this we will need the following three
simple observations. We make them in the context where $V$, $V_{*}$
and $Q$ are closed cubes which satisfy $V\subset V_{*}\subset Q$,
and $F_{+}$ and $F_{-}$ are two disjoint compact subsets of $Q$
which are unions of finitely many closed rectangles. Then of course
$V_{*}\cap F_{+}$ and $V_{*}\cap F_{-}$ are disjoint compact subsets
of $V_{*}$ which are finite unions of closed rectangles. (We could
not make an analogous claim if we considered finite unions of cubes
instead of rectangles.) 

Our observations are that

\noindent \begin{flushleft}
\noindent(i) $V$ is an exceptional $(Q,F_{+},F_{-})$-subcube if
and only if \\
it is an exceptional $(V_{*},V_{*}\cap F_{+},V_{*}\cap F_{-})$-subcube,
\par\end{flushleft}

\noindent \begin{flushleft}
\noindent(ii) $V$ is a good $(\tau',Q,F_{+},F_{-})$-subcube if
and only if \\
it is a good $(\tau',V_{*},V_{*}\cap F_{+},V_{*}\cap F_{-})$-subcube, 
\par\end{flushleft}

\noindent \begin{flushleft}
and so, in view of (i) and (ii),
\par\end{flushleft}

\noindent \begin{flushleft}
\noindent(iii) $V$ is $(\tau',Q,F_{+},F_{-})$-minimal if and only
if \\
it is good $(\tau',V_{*},V_{*}\cap F_{+},V_{*}\cap F_{-})$-minimal.
\par\end{flushleft}
\end{rem}

When we try to obtain a positive answer for Question A$^{\prime}$
we of course have to start with a given cube $Q$ and subsets $F_{+}$
and $F_{-}$ of $Q$ such that, in the language of Definition \ref{def:toptimal},
$Q$ itself is a good $\tau'$-subcube. Our simplification, which
we will state formally in a moment (as Theorem \ref{thm:newred}),
is that we only have to consider the special case where $Q$ is also
$\tau'$-minimal. Thus it would be of interest to study the properties
of $\tau'$-minimal cubes. One simple and easily established property
of such cubes will be given below in Lemma \ref{lem:aprelim}. Later
(in Subsection \ref{sub:AnotherQ}) the reader will be invited to
consider whether $\tau'$-minimal cubes have a certain other simple
but much less evident property. If they do, this would lead to an
affirmative answer to Questions A$^{\prime}$ and A.

Here then is the reduction of Question A$^{\prime}$ alluded to above. 
\begin{thm}
\label{thm:newred}Let $\tau'$ and $s$ be constants satisfying $\tau'\in(0,1/2)$
and $s\in(0,1/2]$. Let $Q$ be a closed cube in $\mathbb{R}^{d}$. 

Suppose that, whenever $F_{+}$ and $F_{-}$ are two disjoint compact
subsets of $Q$ which are each unions of finitely many closed rectangles,
and are such that $Q$ is a $\tau'$-minimal subcube (of itself),
then there exists a cube $W$ contained in $Q$ such that\textbf{\begin{equation}
\min\left\{ \lambda(W\cap F_{+}),\lambda(W\cap F_{-})\right\} \ge s\lambda(W)\,.\label{eq:omqd}\end{equation}
}Then Question A$^{\prime}$ has an affirmative answer for the same
constants $\tau'$ and $s$\,.
\end{thm}

It will be convenient to present some parts of the proof of this theorem
separately in the following lemma.
\begin{lem}
\label{lem:aprelim}Let $Q$, $F_{+}$, $F_{-}$ and $\tau'$ be as
in Definition \ref{def:toptimal}. 

(i) If there exists an exceptional subcube of $Q$, then there exists
a subcube $W$ of $Q$ which satisfies \begin{equation}
\lambda(W\cap F_{+})=\lambda(W\cap F_{-})=\min\left\{ \lambda(W\cap F_{+}),\lambda(W\cap F_{-})\right\} =\frac{1}{2}\lambda(W)\,.\label{eq:urp}\end{equation}

(ii) If $V$ is a $\tau'$-minimal subcube of $Q$, then it satisfies
\begin{equation}
\min\left\{ \lambda(V\cap F_{+}),\lambda(V\cap F_{-})\right\} =\tau'\lambda(V\setminus F_{+}\setminus F_{-})>0\label{eq:fap}\end{equation}
and every subcube $W$ of $V$ which is strictly smaller than $V$
satisfies at least one of the two conditions\begin{equation}
\min\left\{ \lambda(W\cap F_{+}),\lambda(W\cap F_{-})\right\} =0\label{eq:asap}\end{equation}
and \begin{equation}
\min\left\{ \lambda(W\cap F_{+}),\lambda(W\cap F_{-})\right\} <\tau'\lambda(W\setminus F_{+}\setminus F_{-})\,.\label{eq:bsap}\end{equation}

\end{lem}
\noindent\textit{Proof.} We first deal with part (i). Suppose that
$V$ is an exceptional subcube of $Q$. We want to use Lemma \ref{lem:prototype}.
The roles of the cube denoted by $Q$ and of the subset $E$ of $Q$
in the statement of that lemma will now be played here, respectively,
by the cube $V$ and by its subset $V\cap F_{+}$. To apply the lemma
we need to know that $\lambda(V\cap F_{+})>0$ which is part of the
definition of exceptional subcubes, and we also need to know that
$\lambda(V\cap F_{+})<\lambda(V)$. This second inequality follows
from the given condition $\lambda(V\cap F_{-})>0$ (again part of
the definition) and the inclusion $V\cap F_{-}\subset V\setminus F_{+}=V\setminus(V\cap F_{+})$
which together give that $0<\lambda(V\cap F_{-})\le\lambda(V)-\lambda(V\cap F_{+})$.
Thus Lemma \ref{lem:prototype} can be applied to provide us with
a subcube $W$ of $V$ for which \begin{equation}
\lambda\left(W\setminus(V\cap F_{+})\right)=\lambda\left(W\cap(V\cap F_{+})\right)=\frac{1}{2}\lambda(W)\,.\label{eq:sbema}\end{equation}
Obviously \begin{equation}
\lambda\left(W\cap(V\cap F_{+})\right)=\lambda\left(W\cap F_{+}\right)\,,\label{eq:mdpcm}\end{equation}
and, since $V$ is essentially the union of $V\cap F_{+}$ and $V\cap F_{-}$,
it will also be easy to deduce that \begin{equation}
\lambda\left(W\setminus(V\cap F_{+})\right)=\lambda\left(W\cap F_{-}\right)\,.\label{eq:meb}\end{equation}
More explicitly, since $F_{+}\cap F_{-}=\emptyset$ and $W\subset V$,
we have \begin{eqnarray*}
W\cap F_{-} & \subset & W\setminus F_{+}\subset W\setminus(V\cap F_{+})\subset W\setminus(W\cap F_{+})\\
 & \subset & W\setminus F_{+}=\left(\left(W\setminus F_{+}\right)\setminus F_{-}\right)\cup\left(\left(W\setminus F_{+}\right)\cap F_{-}\right)\\
 & \subset & \left(V\setminus F_{+}\setminus F_{-}\right)\cup\left(W\cap F_{-}\right)\,.\end{eqnarray*}
These inclusions and then the fact that $V$ is an exceptional subcube,
imply that \begin{eqnarray*}
\lambda(W\cap F_{-}) & \le & \lambda\left(W\setminus(V\cap F_{+})\right)\\
 & \le & \lambda\left(V\setminus F_{+}\setminus F_{-}\right)+\lambda\left(W\cap F_{-}\right)=0+\lambda\left(W\cap F_{-}\right)\end{eqnarray*}
which establishes (\ref{eq:meb}). The required formula (\ref{eq:urp})
now follows immediately from (\ref{eq:sbema}), (\ref{eq:mdpcm})
and (\ref{eq:meb}). 

We now deal with part (ii) of the lemma. Suppose that $V$ is a $\tau'$-minimal
subcube of $Q$. The fact that every strictly smaller subcube $W$
of $V$ satisfies at least one of the two conditions (\ref{eq:asap})
and (\ref{eq:bsap}) is simply a rewriting of definitions. More explicitly,
suppose that some such subcube $W$ fails to satisfy both (\ref{eq:asap})
and (\ref{eq:bsap}). Then, in the case where $\tau'\lambda(W\setminus F_{+}\setminus F_{-})=0$\,,
this implies that $W$ is an exceptional subcube. In the case where
$\tau'\lambda(W\setminus F_{+}\setminus F_{-})>0$ this implies that
$W$ is a good $\tau'$-subcube. The given condition on $V$ excludes
both of these possibilities, so at least one of (\ref{eq:asap}) and
(\ref{eq:bsap}) must be satisfied. 

Finally, suppose that $V$ does not satisfy (\ref{eq:fap}). Then,
since $V$ is a good $\tau'$-subcube, it must satisfy \begin{equation}
\min\left\{ \lambda(V\cap F_{+}),\lambda(V\cap F_{-})\right\} >\tau'\lambda(V\setminus F_{+}\setminus F_{-})>0\,.\label{eq:vfpm}\end{equation}
Now let $W$ be a subcube of $V$ such that $\lambda(V\setminus W)<\varepsilon$.
If $\varepsilon$ is sufficiently small, it follows from (\ref{eq:vfpm})
(see Remark \ref{rem:brichj}) that \begin{equation}
\min\left\{ \lambda(W\cap F_{+}),\lambda(W\cap F_{-})\right\} >\tau'\lambda(W\setminus F_{+}\setminus F_{-})>0\label{eq:bromx}\end{equation}
which means that $W$ is also a good $\tau'$-subcube, contradicting
the $\tau'$-minimality of $V$. This shows that (\ref{eq:fap}) holds
and so completes the proof of the lemma. $\qed$
\begin{rem}
\label{rem:brichj}More explicitly, the inequalities which enable
us to deduce (\ref{eq:bromx}) from (\ref{eq:vfpm}) for sufficiently
small $\varepsilon$ are, first of all \begin{eqnarray*}
\lambda(W\cap F_{+}) & \le & \lambda(V\cap F_{+})\le\lambda(W\cap F_{+})+\lambda\left((V\setminus W)\cap F_{+}\right)\\
 & \le & \lambda(W\cap F_{+})+\varepsilon\,,\end{eqnarray*}
then the counterpart of this where $F_{+}$ is replaced by $F_{-}$,
and then, finally, \begin{eqnarray*}
\lambda(W\setminus F_{+}\setminus F_{-}) & \le & \lambda(V\setminus F_{+}\setminus F_{-})\le\lambda(V\setminus W\setminus F_{+}\setminus F_{-})+\lambda(W\setminus F_{+}\setminus F_{-})\\
 & \le & \varepsilon+\lambda(W\setminus F_{+}\setminus F_{-})\,.\end{eqnarray*}

\end{rem}
\noindent \textit{The proof of Theorem \ref{thm:newred}.} Our approach
here will have some features in common with the proof of Theorem \ref{thm:maingt}.
Let $Q$ be an arbitrary cube in $\mathbb{R}^{d}$ and let $F_{+}$
and $F_{-}$ be arbitrary disjoint subsets of $Q$ which are both
finite unions of bounded closed rectangles and satisfy (\ref{eq:ehyp}).
We have to show that there exists a subcube $W$ of $Q$ which satisfies
(\ref{eq:econ}). In one case this is very easy to do, namely when
$Q$ has a subcube which is an exceptional cube. We simply invoke
part (i) of Lemma \ref{lem:aprelim} to obtain a cube $W$ satisfying
(\ref{eq:urp}) and therefore (\ref{eq:econ}), since $s\le1/2$.
This leaves us free to assume, for the rest of this proof, that $Q$
does not contain any exceptional subcubes. 

Let us now consider the collection $\mathcal{G}$ of all subcubes
of $Q$ which are good $\tau'$-subcubes. This is non empty since
$Q$ itself is such a cube. Since $F_{+}$ and $F_{-}$ are disjoint
and compact, we have that \[
\rho:=\mathrm{dist}\left(F_{+},F_{-}\right)>0\,.\]
Since every good cube must intersect with both $F_{+}$ and $F_{-}$,
it follows that the diameter of each $V\in\mathcal{G}$ satisfies
$\rho\le\mathrm{diam}\, V\le\mathrm{diam}\, Q$\,. Consequently the
infimum $\rho_{*}:=\inf_{V\in\mathcal{G}}\mathrm{diam}\, V$ is strictly
positive. Furthermore, there exists a sequence of cubes $\left\{ V_{n}\right\} _{n\in\mathbb{N}}$
in $\mathcal{G}$ such that $\lim_{n\to\infty}\mathrm{diam\,}V_{n}=\rho_{*}$.
For each $n$, let $V_{n}=Q\left(x_{n},r_{n}\right)$ (here again
using the standard notation (\ref{eq:defqxr})). By passing to a subsequence
of $\left\{ V_{n}\right\} _{n\in\mathbb{N}}$ if necessary, we can
assume that the sequences $\left\{ x_{n}\right\} _{n\in\mathbb{N}}$
and $\left\{ r_{n}\right\} _{n\in\mathbb{N}}$ converge, respectively,
to a point $x_{*}\in Q$, and to the positive number $r_{*}=\rho_{*}/2\sqrt{d}$\,.
We let $V_{*}=Q(x_{*},r_{*})$. Since each $V_{n}$ is in $\mathcal{G}$
we have \begin{equation}
\min\left\{ \lambda(F_{+}\cap Q(x_{n},r_{n})),\lambda(F_{-}\cap Q(x_{n},r_{n})\right\} \ge\tau'\lambda\left(Q(x_{n},r_{n})\setminus F_{+}\setminus F_{-}\right)>0\,.\label{eq:ntn}\end{equation}
We can pass to the limit in these inequalities, using three applications
of part (ii) of Lemma \ref{lem:SimplePropsCubes}, where we choose
the set $A$ to be, respectively, $F_{+}$, $F_{-}$ and $\mathbb{R}^{d}\setminus F_{+}\setminus F_{-}$.
This gives \begin{equation}
\min\left\{ \lambda(F_{+}\cap Q(x_{*},r_{*})),\lambda(F_{-}\cap Q(x_{*},r_{*})\right\} \ge\tau'\lambda\left(Q(x_{*},r_{*})\setminus F_{+}\setminus F_{-}\right)\ge0\,.\label{eq:tmtdb}\end{equation}

In fact, we have \begin{equation}
\lambda\left(Q(x_{*},r_{*})\setminus F_{+}\setminus F_{-}\right)>0\,.\label{eq:spy}\end{equation}
We will defer the proof of (\ref{eq:spy}) for a moment. It follows
from part (i) of Lemma \ref{lem:SimplePropsCubes} that $V_{*}\subset Q$.
This, together with (\ref{eq:tmtdb}) and (\ref{eq:spy}), implies
that $V_{*}$ is a good $\tau'$-subcube. Since $\mathrm{diam}\, V_{*}=\rho_{*}$
no strictly smaller subcube of $V_{*}$ can be a good $\tau'$-subcube.
Furthermore, using the assumption that we showed that we could make
above, no subcube of $V_{*}$ can be exceptional. This means that
$V_{*}$ is $\tau'$-minimal. Thus, in the terminology of Remark \ref{rem:splurg},
the cube $V_{*}$ is also $\left(\tau',V_{*},V_{*}\cap F_{+},V_{*}\cap F_{-}\right)$-minimal.
Now we are ready to invoke the hypothesis of our theorem, but where
here the role of the cube $Q$ is now played by $V_{*}$ and the roles
of the sets $F_{+}$ and $F_{-}$ are played by $V_{*}\cap F_{+}$
and $V_{*}\cap F_{-}$. The hypothesis ensures that there is a subcube
$W$ of $V_{*}$, and therefore also of $Q$ which satisfies\[
\min\left\{ \lambda(W\cap V_{*}\cap F_{+}),\lambda(W\cap V_{*}\cap F_{-})\right\} \ge s_{0}\lambda(W)\,.\]
Since of course $W\cap V_{*}=W$ and $s_{0}\ge s$ the cube $W$ satisfies
(\ref{eq:econ}) and is therefore the cube required to complete the
proof of the theorem.

It remains only to show that (\ref{eq:spy}). Suppose, on the contrary,
that \begin{equation}
\lambda\left(Q(x_{*},r_{*})\setminus F_{+}\setminus F_{-}\right)=0\,.\label{eq:urz}\end{equation}
Since we have excluded the possibility that the cube $V_{*}=Q(x_{*},r_{*})$
is an exceptional subcube, this means that at least one of $\lambda(V_{*}\cap F_{+})$
and $\lambda(V_{*}\cap F_{-})$ must equal $0$. We can suppose that
$\lambda(V_{*}\cap F_{-})=0$, since the other case, where $\lambda(V_{*}\cap F_{+})=0$,
can be treated exactly analogously. This supposition, together with
(\ref{eq:urz}), implies that \[
0=\lambda(V_{*}\setminus F_{+}\setminus F_{-})=\lambda(V_{*}\setminus F_{+})-\lambda\left((V_{*}\setminus F_{+})\cap F_{-}\right)=\lambda(V_{*}\setminus F_{+})\]
and so $\lambda(V_{*})=\lambda(V_{*}\setminus F_{+})+\lambda(V_{*}\cap F_{+})=\lambda(V_{*}\cap F_{+})$.
This means that every subset of $V_{*}$ with positive measure must
contain points of $F_{+}$. Consider the cube $Q\left(x_{*},r_{*}+\frac{\rho}{2\sqrt{d}}\right)$,
where, as above, $\rho=\mathrm{dist}(F_{+},F_{-})$. Every point in
this cube is at distance strictly less than $\rho$ from some points
in $V_{*}$ and therefore also at distance strictly less than $\rho$
from points in $F_{+}$. Therefore,\begin{equation}
F_{-}\cap Q\left(x_{*},r_{*}+\frac{\rho}{2\sqrt{d}}\right)=\emptyset\,.\label{eq:nta}\end{equation}
On the other hand, if $n$ is chosen large enough, we obtain that
\[
\left\Vert x_{*}-x_{n}\right\Vert _{\ell_{d}^{\infty}}+r_{n}<r_{*}+\frac{\rho}{2\sqrt{d}}\,,\]
which implies that the cube $Q(x_{n},r_{n})$ is contained in $Q\left(x_{*},r_{*}+\frac{\rho}{2\sqrt{d}}\right)$.
In view of (\ref{eq:ntn}) the cube $Q(x_{n},r_{n})$ contains points
of $F_{-}$. This contradicts (\ref{eq:nta}), and thus shows that
(\ref{eq:urz}) cannot hold. This proves (\ref{eq:spy}) and completes
the proof of the theorem. $\qed$

\subsection{\label{sub:AnotherQ}Another question which should be considered.
\phantom{AAAAAAAAAAAAAAAAAAAAAAAAAAAAAAAAAAAAAAAAAAAAAAAAAAAAAAAAAAAAAAAAAAAAAAAAAAAAAAAAAAAAAAAAAAAAAAAAAAAAAAAAAAAAAAAAAAAAAAAAAAAAAAAAAAAAAAAAAA$A_{A_A}$AAAAAAAAAAAAAAAAAAAA}}

As already remarked in the previous subsection, in view of Theorem
\ref{thm:newred}, it could be very helpful if we can discover some
concrete consequences of the apparently very stringent condition on
a cube that it is $\tau'$-minimal. More explicitly we invite the
reader to consider the following question, which we shall call Question
B.

\textbf{\textit{Let $Q$ be a cube in $\mathbb{R}^{d}$ and let $F_{+}$
and $F_{-}$ be disjoint compact subsets of $Q$ which are each finite
unions of closed rectangles. Suppose further that $Q$ is $\tau'$-minimal
for some $\tau'\in(0,1/2)$. Does this imply that $\lambda(F_{+})=\lambda(F_{-})$?}}

This question attracts our attention for reasons expressed by the
following two propositions.
\begin{prop}
\label{pro:cdeo}In the case where $d=1$, the answer to Question
B is affirmative for every choice of the constant $\tau'\in(0,1/2)$.
\end{prop}

\begin{prop}
\label{pro:bimpa}If Question B has an affirmative answer for arbitrary
dimension $d$ and for some constant $\tau'\in(0,1/2)$ which does
not depend on $d$\,, then this implies affirmative answers for Questions
A and A$^{\prime}$.
\end{prop}
In the light of these two propositions we have very strong motivation
for attempting to answer Question B for the case $d=2$. An affirmative
answer may point the way to an affirmative answer for all $d$ and
thus also for Questions A$^{\prime}$ and A. A negative answer for
$d=2$ would apparently lead to a negative answer for all $d>2$.
While this in itself would not imply negative answers to Questions
A or A$^{\prime}$ it could perhaps indicate some path towards such
negative answers.\medskip{}

\noindent \textit{Proof of Proposition \ref{pro:cdeo}.} We begin
with the remark that, since $d=1$, we can and will make use of the
following convenient fact: 

\begin{center}
\textit{Whenever $A$ and $B$ are cubes with a common endpoint and
}\\
\textit{such that $A\subset B$, then $B\setminus A$ also coincides
a.e.~with a cube. }
\par\end{center}

\noindent Unfortunately this fact is not available when $d>1$.

Let $\tau'$ be an arbitrary constant in $\left(0,1/2\right)$. Let
$Q$ be a closed interval $Q=[a,b]$ and suppose that $F_{+}$ and
$F_{-}$ are disjoint subsets of $Q$ which are each unions of finitely
many closed intervals. (In fact our proof will only use the fact that
$F_{+}$ and $F_{-}$ are measurable and disjoint). Suppose furthermore
that $Q$ is $\tau'$- minimal, (i.e., that it is $\left(\tau',Q,F_{+},F_{-}\right)$-minimal
in the notation of Remark \ref{rem:splurg}). It will be convenient
to write $G:=Q\setminus F_{+}\setminus F_{-}$. We may apply part
(ii) of Lemma \ref{lem:aprelim} to obtain that equality must hold
in at least one of the two inequalities \[
\lambda(F_{+})\ge\tau'\lambda(G)\mbox{ and }\lambda(F_{-})\ge\tau'\lambda(G).\]
(Cf.~(\ref{eq:fap}).) We will suppose, without loss of generality,
that $\lambda(F_{-})=\tau'\lambda(G)$. (The other case can be treated
exactly analogously.) 

Suppose, in contradiction to what we seek to prove, that $\lambda(F_{+})\ne\lambda(F_{-})$.
For this to happen we must have \begin{equation}
\lambda(F_{+})>\tau'\lambda(G)\,.\label{eq:xrtp}\end{equation}

Let us define the functions $u:[a,b]\to\mathbb{R}$ and $v:[a,b]\to\mathbb{R}$
by \[
u(t)=\lambda\left([a,t]\cap F_{-}\right)-\tau'\lambda\left([a,t]\cap G\right)\]
and \[
v(t)=\lambda\left([t,b]\cap F_{-}\right)-\tau'\lambda\left([t,b]\cap G\right)\,.\]

Note that $u(t)=v(t)=0$ for $t=a$ and for $t=b$. Furthermore \begin{equation}
u(t)+v(t)=\lambda(F_{-})-\tau'\lambda(G)=0\mbox{ for all }t\in[a,b]\,.\label{eq: baawer}\end{equation}

Suppose that $u(s)=0$ for some $s\in(a,b)$. Then $v(s)=0$. Since
$Q$ is the union of the two non overlapping intervals $I_{1}:=[a,s]$
and $I_{2}:=[s,b]$ it follows from (\ref{eq:xrtp}) that the inequality
\[
\lambda\left(F_{+}\cap I_{j}\right)>\tau'\lambda\left(G\cap I_{j}\right)\]
must hold for at least one of the two values $j=1$ and $j=2$. Since
we also have $\lambda\left(F_{-}\cap I_{j}\right)=\tau'\lambda\left(G\cap I_{j}\right)$
for both these values of $j$, it follows that $I_{j}$ is a good
$\tau'$-subcube of $Q$ for at least one value of $j$. This contradicts
the $\tau'$-minimality of $Q$ and shows that we must have $u(s)\ne0$
for all $s\in(a,b)$. It follows by an analogous argument, or simply
from (\ref{eq: baawer}), that $v(s)\ne0$ for all $s\in(a,b)$. Since
both $u$ and $v$ are continuous functions, they cannot change sign
on $(a,b)$ and so one of them is strictly positive and the other
is strictly negative on the whole interval $(a,b)$. 

Let us consider the first case, where $u(t)>0$ for all $t\in(a,b)$.
In view of (\ref{eq:xrtp}) the continuous function \[
w(t):=\lambda\left([a,t]\cap F_{+}\right)-\tau'\lambda\left([a,t]\cap G\right)\]
is strictly positive for $t=b$ and therefore for some choice (in
fact infinitely many choices) of $s\in(a,b)$ sufficiently close to
$b$ we have $w(s)>0$. This, together with the fact that $u(s)>0$,
implies that the interval $[a,s]$ is a $\tau'$-good subcube of $[a,b]$
for these values of $s$. This contradicts the supposition that $[a,b]$
is $\tau'$-minimal. In the remaining case, where $v(t)>0$ for all
$t\in(a,b)$, an analogous argument implies that $[s,b]$ is a $\tau'$-good
subcube for all values of $s\in(a,b)$ sufficiently close to $a$
and thus also gives a contradiction. Since we have shown that, in
all cases, the assumption that $\lambda(F_{+})\ne\lambda(F_{-})$
leads to a contradiction, our proof is complete. $\qed$ 

\bigskip{}

\noindent \textit{Proof of Proposition \ref{pro:bimpa}.} Suppose
then that Question B has an affirmative answer for all $d\in\mathbb{N}$
and for at least one value of the constant $\tau'\in(0,1/2)$, a value
which does not depend on $d$. We will see that this implies that
the hypothesis of Theorem \ref{thm:newred} is always satisfied, in
fact for a positive constant $s$ which depends only on the given
constant $\tau'$. 

Explicitly, suppose that $Q$ is a cube in $\mathbb{R}^{d}$ and that
$F_{+}$ and $F_{-}$ are two disjoint subsets of $Q$ which are finite
unions of closed rectangles. Suppose furthermore that $Q$ is a $\tau'$-minimal
subcube of itself, i.e., a $\left(\tau',Q,F_{+},F_{-}\right)$-minimal
subcube in the terminology of Remark \ref{rem:splurg}. In order to
be able to apply Theorem \ref{thm:newred} we have to find a subcube
$W$ of $Q$ which satisfies (\ref{eq:omqd}) for a value of $s$
which does not depend on $d$. We will show that we can simply take
$W=Q$.

In view of the supposed affirmative answer to Question B, we have
that $\lambda(F_{+})=\lambda(F_{-})$. This, together with the fact
that $Q$ is a good $\tau'$-subcube (of itself), implies that \begin{eqnarray*}
\lambda(Q) & = & \lambda(F_{+})+\lambda(F_{-})+\lambda(Q\setminus F_{+}\setminus F_{-})\\
 & = & 2\min\left\{ \lambda(Q\cap F_{+}),\lambda(Q\cap F_{-})\right\} +\lambda(Q\setminus F_{+}\setminus F_{-})\\
 & \ge & 2\tau'\lambda(Q\setminus F_{+}\setminus F_{-})+\lambda(Q\setminus F_{+}\setminus F_{-}).\end{eqnarray*}
The previous three lines imply that $\lambda(Q\setminus F_{+}\setminus F_{-})\le\frac{1}{2\tau'+1}\lambda(Q)$
and also that \[
2\min\left\{ \lambda(Q\cap F_{+}),\lambda(Q\cap F_{-})\right\} =\lambda(Q)-\lambda(Q\setminus F_{+}\setminus F_{-})\]
which together give us that \begin{eqnarray*}
\min\left\{ \lambda(Q\cap F_{+}),\lambda(Q\cap F_{-})\right\}  & = & \frac{1}{2}\left(\lambda(Q)-\lambda(Q\setminus F_{+}\setminus F_{-})\right)\\
 & \ge & \frac{1}{2}\left(1-\frac{1}{2\tau'+1}\right)\lambda(Q)\,.\end{eqnarray*}
Thus the cube $Q=W$ indeed satisfies (\ref{eq:omqd}) for \[
s=\frac{1}{2}\left(1-\frac{1}{2\tau'+1}\right)=\frac{1}{2+1/\tau'}\,.\]
This shows that an affirmative answer to Question B ensures the validity
of the condition that is required in Theorem \ref{thm:newred} to
imply a positive answer to Question A$^{\prime}$ for the given constant
$\tau'$ and for $s=\frac{1}{2+1/\tau'}$. Therefore an affirmative
answer to Question B also implies a positive answer to Question A
for any $\tau\in(\tau',1/2)$ and for the same value of $s$. $\qed$

\section{\label{sec:Appendices}Appendices}

The material in this section is quite standard and/or elementary.
This section, or some parts of it, may be removed from future versions
of the paper.

\subsection{\label{sub:mediansandoscillation}Medians and the mean oscillation
of a function}

Throughout this subsection $\left(\Omega,\Sigma,\lambda\right)$ is
an arbitrary measure space, $E$ is a measurable subset of $\Omega$
satisfying $0<\lambda(E)<\infty$, and $f$ is a measurable real function
whose domain of definition contains $E$. If $f$ is integrable we
set $f_{E}=\frac{1}{\lambda(E)}\int_{E}fd\lambda$. We will sometimes
use the notation $f^{-1}(H)$ to mean the set $\left\{ x\in E:f(x)\in H\right\} $,
where $H$ is some subset of $\mathbb{R}$.

In essentially all our applications in this paper of the (standard)
results of this subsection, $\Omega$ is $\mathbb{R}^{d}$ or some
Lebesgue measurable subset of $\mathbb{R}^{d}$, and $\Sigma$ consists
of all Lebesgue measurable subsets of $\Omega$\,, and $\lambda$
is $d$-dimensional Lebesgue measure.

\medskip{}

\begin{lem}
\label{lem:medianexists}There exists at least one median of $f$
on $E$, i.e., a number $m$ which satisfies \begin{equation}
\lambda\left(\left\{ x\in E:f(x)<m\right\} \right)\le\frac{1}{2}\lambda\left(E\right)\mbox{ and }\lambda\left(\left\{ x\in E:f(x)>m\right\} \right)\le\frac{1}{2}\lambda\left(E\right)\,.\label{eq:fynz}\end{equation}

\end{lem}
\noindent \textit{Proof.} Let \[
A=\left\{ \alpha\in\mathbb{R}:\lambda\left(\left\{ x\in E:f(x)>\alpha\right\} \right)>\frac{1}{2}\lambda\left(E\right)\right\} \]
and \[
B=\left\{ \alpha\in\mathbb{R}:\lambda\left(\left\{ x\in E:f(x)<\alpha\right\} \right)>\frac{1}{2}\lambda\left(E\right)\right\} \,.\]
Both of these sets are non empty, by the expanding sequence theorem.
The set $A$ is an interval whose left endpoint is $-\infty$. The
set $B$ is an interval whose right endpoint is $+\infty$. It is
also clear that $A\cap B$ is empty. Thus the set $M=\mathbb{R}\backslash A\backslash B$
is non empty and is a bounded interval (which may also happen to be
a single point) whose endpoints are $\sup A$ and $\inf B$. Every
number $m\in M$ must satisfy (\ref{eq:fynz}). $\qed$
\begin{lem}
\label{lem:medianoptimal}The mean oscillation of $f$ on $E$ satisfies
\begin{equation}
\mathbf{O}(f,E)=\frac{1}{\lambda(E)}\int_{E}\left|f-m\right|d\lambda\label{eq:vfm}\end{equation}
for each median $m$ of $f$ on $E$. 
\end{lem}
\noindent \textit{Proof.} Suppose that $m$ is any median of $f$
on $E$. This implies that\begin{equation}
\lambda\left(\left\{ x\in E:f(x)<m\right\} \right)\le\frac{1}{2}\lambda\left(E\right)\le\lambda\left(\left\{ x\in E:f(x)\ge m\right\} \right)\,.\label{eq:hromp}\end{equation}
We can perform the following calculation for each number $c$ satisfying
$c\le m$. The transition from the third and fourth lines to the fifth
and sixth lines uses (\ref{eq:hromp}) and the fact that $m-c\ge0$.
\begin{eqnarray*}
 &  & \int_{E}\left|f-m\right|d\lambda\\
 & = & \int_{f^{-1}\left((-\infty,m)\right)}(m-f)d\lambda+\int_{f^{-1}\left([m,\infty)\right)}(f-m)d\lambda\\
 & = & \int_{f^{-1}\left((-\infty,m)\right)}(c-f)d\lambda+\int_{f^{-1}\left([m,\infty)\right)}(f-c)d\lambda\\
 &  & +(m-c)\lambda(f^{-1}\left((-\infty,m)\right))+(c-m)\lambda(f^{-1}\left([m,\infty)\right))\\
 & \le & \int_{f^{-1}\left((-\infty,m)\right)}(c-f)d\lambda+\int_{f^{-1}\left([m,\infty)\right)}(f-c)d\lambda\\
 &  & +(m-c)\lambda(f^{-1}\left([m,\infty)\right))+(c-m)\lambda(f^{-1}\left([m,\infty)\right))\\
 & = & \int_{f^{-1}\left((-\infty,m)\right)}(c-f)d\lambda+\int_{f^{-1}\left([m,\infty)\right)}(f-c)d\lambda+0_{_{-}}\\
 & \le & \int_{f^{-1}\left((-\infty,m)\right)}\left|c-f\right|d\lambda+\int_{f^{-1}\left([m,\infty)\right)}\left|f-c\right|d\lambda=\int_{E}\left|f-c\right|d\lambda\,.\end{eqnarray*}
In the remaining case, i.e., for each $c>m$, we can apply the previous
calculation to the function $-f$. Since $-m$ is a median of $-f$
on $E$ and $-c<-m$, we obtain that \[
\int_{E}\left|-f+m\right|d\lambda\le\int_{E}\left|-f+c\right|d\lambda\,.\]
From these two cases we see that \[
\int_{E}\left|f-m\right|d\lambda\le\int_{E}\left|f-c\right|d\lambda\]
for every $c\in\mathbb{R}$\,.We can now obtain (\ref{eq:vfm}) by
taking the infimum over all $c\in\mathbb{R}$. $\qed$
\begin{lem}
\label{lem:gtyl} The inequality \begin{equation}
\left|f_{E}-m\right|\le\mathbf{O}(f,E)\label{eq:gty}\end{equation}
holds for every median $m$ of $f$ on $E$.
\end{lem}
\noindent \textit{Proof.} \[
f_{E}-m=\frac{1}{\lambda(E)}\int_{E}(f_{E}-m)d\lambda=\frac{1}{\lambda(E)}\int_{E}(f-m)d\lambda\le\frac{1}{\lambda(E)}\int_{E}|f-m|d\lambda=\mathbf{O}(f,E)\]
and essentially the same argument shows that $m-f_{E}\le\mathbf{O}(f,E)$.
These two estimates of course give us (\ref{eq:gty}). $\qed$

Let us now prove (\ref{eq:avao}).

Suppose that $m$ is a median of $f$ on $E$. Then, by Lemma \ref{lem:medianoptimal},
\begin{eqnarray*}
\frac{1}{\lambda(E)}\int_{E}\left|f(x)-m\right|d\lambda(x) & \le & \frac{1}{\lambda(E)}\int_{E}\left|f(x)-f_{E}\right|d\lambda(x)\\
 & = & \frac{1}{\lambda(E)}\int_{E}\left|\frac{1}{\lambda(E)}\int_{E}\left(f(x)-f(y)\right)d\lambda(y)\right|d\lambda(x)\\
 & \le & \frac{1}{\lambda(E)^{2}}\iint_{E\times E}\left|f(x)-f(y)\right|d\lambda(x)d\lambda(y)\end{eqnarray*}
establishing the first two inequalities of (\ref{eq:avao}). For the
remaining inequality we observe that

\begin{eqnarray*}
 &  & \frac{1}{\lambda(E)^{2}}\iint_{E\times E}\left|f(x)-f(y)\right|d\lambda(x)d\lambda(y)\\
 & \le & \frac{1}{\lambda(E)^{2}}\iint_{E\times E}\left|f(x)-m\right|d\lambda(x)d\lambda(y)+\frac{1}{\lambda(E)^{2}}\iint_{E\times E}\left|m-f(y)\right|d\lambda(x)d\lambda(y)\\
 & = & \frac{2}{\lambda(E)}\int_{E}\left|f-m\right|d\lambda\,,\end{eqnarray*}
which completes the proof of (\ref{eq:avao}).
\begin{lem}
Let $m_{1}$ and $m_{2}$ be medians of $f$ on $E$ with $m_{1}<m_{2}$.
Then \[
\lambda\left(\left\{ x\in E:m_{1}<f(x)<m_{2}\right\} \right)=0\,.\]

\end{lem}
\noindent \textit{Proof.} The sets $f^{-1}\left((-\infty,m_{1}]\right)$
and $f^{-1}\left([m_{2},\infty)\right)$ must each have $\lambda$
measure greater than or equal to $\frac{1}{2}\lambda(E)$. So, since
they are disjoint subsets of $E$, they must in fact both have measure
equal to $\frac{1}{2}\lambda(E)$. Thus $f^{-1}\left((m_{1},m_{2})\right)$,
the complement of their union in $E$, must have zero $\lambda$ measure.
$\qed$

\subsection{\label{sub:dfojne}Implications between different forms of the John-Nirenberg
inequality.}

As mentioned in Section \ref{sec:Notation-and-terminology}, various
papers present slightly different versions of the inequality (\ref{eq:jne}).
Let us explicitly show the connections between these versions. Although
we have formulated some of these connections for the special case
where the collection $\mathcal{E}$ of admissible subsets of $D$
is taken to consist only of cubes, i.e, for $\mathcal{E}=\mathcal{Q}(D)$,
the same connections clearly apply when $\mathcal{E}$ is chosen to
be any other collection of admissible subsets, even not necessarily
all contained in $D$.

We have to consider three kinds of {}``transition'':

$\bullet$ In some variants of (\ref{eq:jne}) the seminorm $\left\Vert f\right\Vert _{BMO(D,\mathcal{Q}(D))}$
or $\left\Vert f\right\Vert _{BMO(D,\mathcal{Q}(D))}^{(\mathbf{D})}$
appears in place of $\left\Vert f\right\Vert _{BMO(D,\mathcal{Q}(D))}^{(\mathbf{A})}$.
But, since these seminorms are equivalent to each other to within
a factor of $2$, an inequality using one of them obviously implies
analogous inequalities using the others, of course sometimes with
the constant $b$ replaced by $b/2$.

$\bullet$ Sometimes the set $\left\{ x\in D:\left|f(x)-f_{D}\right|>\alpha\right\} $
may be replaced by the possibly larger set $\left\{ x\in D:\left|f(x)-f_{D}\right|\ge\alpha\right\} $.
But, for $\alpha>0$, this of course does not change anything. The
given inequality remains valid since the right hand side of the inequality
is a continuous function of $\alpha$ and \begin{eqnarray*}
\lambda\left(\left\{ x\in D:\left|f(x)-f_{D}\right|>\alpha\right\} \right) & \le & \lambda\left(\left\{ x\in D:\left|f(x)-f_{D}\right|\ge\alpha\right\} \right)\\
 & = & \lim_{n\to\infty}\lambda\left(\left\{ x\in D:\left|f(x)-f_{D}\right|>\alpha-1/n\right\} \right)\,.\end{eqnarray*}
This explanation cannot be applied for $\alpha=0$. But since we know
(by the argument of Remark \ref{rem:bge1} or some obvious variant
of it) that the constant $B$ satisfies $B\ge1$, the case $\alpha=0$
is a triviality. An exactly analogous argument shows that the given
inequality remains valid if we replace the set $\left\{ x\in D:\left|f(x)-m\right|>\alpha\right\} $
by the possibly larger set $\left\{ x\in D:\left|f(x)-m\right|\ge\alpha\right\} $,
when $m$ is a median of $f$ on $D$\,.

$\bullet$ Sometimes the average $f_{D}$ in the set $\left\{ x\in D:\left|f(x)-f_{D}\right|>\alpha\right\} $
or in the set \[
\left\{ x\in D:\left|f(x)-f_{D}\right|\ge\alpha\right\} \]
may be replaced by a median $m$ of $f$ on $D$. We can invoke the
following lemma to describe the implications of such a change. 
\begin{lem}
\label{lem:av2med}Let $D$ be some subset of $\mathbb{R}^{d}$ and
let $\mathcal{E}$ be some collection of admissible subsets $E$ of
$\mathbb{R}^{d}$. Suppose that $N(f)$ denotes one of the three seminorms
(\ref{eq:osn}), (\ref{eq:blonk}) or (\ref{eq:fronk}) and that it
is known that the inequality \begin{equation}
\lambda\left(\left\{ x\in E:\left|f(x)-f_{E}\right|>\alpha\right\} \right)\le B\lambda(E)\exp\left(-\frac{b\alpha}{N(f)}\right)\label{eq:a1}\end{equation}
holds for some fixed $E\in\mathcal{E}$ and for all $\alpha\ge0$
and for certain fixed positive constants $b$ and $B$. Then $B\ge1$
and the inequality

\[
\lambda\left(\left\{ x\in E:\left|f(x)-m\right|\ge\alpha\right\} \right)\le e^{b}B\lambda(E)\exp\left(-\frac{b\alpha}{N(f)}\right)\]
holds for every median $m$ of $f$ on $E$ and each $\alpha\ge0$.

Conversely, suppose that it is known that the inequality \begin{equation}
\lambda\left(\left\{ x\in E:\left|f(x)-m\right|>\alpha\right\} \right)\le B\lambda(E)\exp\left(-\frac{b\alpha}{N(f)}\right)\label{eq:b1}\end{equation}
holds for some fixed $E\in\mathcal{E}$ and for all $\alpha\ge0$
and for certain fixed positive constants $b$ and $B$ and for some
median $m$ of $f$ on $E$. Then $B\ge1$ and the inequality

\begin{equation}
\lambda\left(\left\{ x\in E:\left|f(x)-f_{E}\right|\ge\alpha\right\} \right)\le e^{b}B\lambda(E)\exp\left(-\frac{b\alpha}{N(f)}\right)\label{eq:b2}\end{equation}
holds for every $\alpha\ge0$. 
\end{lem}
\noindent \textit{Proof.} If the inequality (\ref{eq:a1}) or, respectively,
the inequality (\ref{eq:b1}) holds for all $\alpha\ge0$, then (cf.~Remark
\ref{rem:bge1}) the constant $B$ necessarily satisfies $B\ge1$,
and, furthermore , as already explained above, it also follows that
the same inequality (respectively (\ref{eq:a1}) or (\ref{eq:b1}))
still holds when {}``$>\alpha$'' is replaced by {}``$\ge\alpha$''
in the definition of the set on the left hand side.

In view of Lemma \ref{lem:gtyl} and (\ref{eq:avao}) we have $\left|f_{E}-m\right|\le N(f)$
for every median $m$ of $f$ on $E$. So, for each $x\in E$, \[
\left|f(x)-m\right|-N(f)\le\left|f(x)-f_{E}\right|\]
and \[
\left|f(x)-f_{E}\right|-N(f)\le\left|f(x)-m\right|\,.\]
These inequalities imply, respectively, that

\begin{eqnarray}
\left\{ x\in E:\left|f(x)-m\right|\ge\alpha\right\}  & \subset & \left\{ x\in E:\left|f(x)-f_{E}\right|\ge\alpha-N(f)\right\} \,.\label{eq:hlte}\end{eqnarray}
and that

\begin{eqnarray}
\left\{ x\in E:\left|f(x)-f_{E}\right|\ge\alpha\right\}  & \subset & \left\{ x\in E:\left|f(x)-m\right|\ge\alpha-N(f)\right\} \,.\label{eq:fpdh}\end{eqnarray}

If (\ref{eq:a1}) holds for all $\alpha$, then, by (\ref{eq:hlte}),
for all $\alpha\ge N(f)$, \begin{eqnarray*}
\lambda\left(\left\{ x\in E:\left|f(x)-m\right|\ge\alpha\right\} \right) & \le & \lambda\left(\left\{ x\in E:\left|f(x)-f_{E}\right|\ge\alpha-N(f)\right\} \right)\\
 & \le & B\lambda(E)\exp\left(-\frac{b(\alpha-N(f))}{N(f)}\right)\\
 & = & Be^{b}\lambda(E)\exp\left(-\frac{b\alpha}{N(f)}\right)\,.\end{eqnarray*}
If $\alpha\in[0,N(f))$, then\[
\lambda\left(\left\{ x\in E:\left|f(x)-m\right|\ge\alpha\right\} \right)\le\lambda(E)=\lambda(E)e^{b}\cdot e^{-b}\le\lambda(E)e^{b}\cdot\exp\left(-\frac{b\alpha}{N(f)}\right)\]
and so, in all cases, \[
\lambda\left(\left\{ x\in E:\left|f(x)-m\right|\ge\alpha\right\} \right)\le e^{b}B\lambda(E)\exp\left(-\frac{b\alpha}{N(f)}\right)\,,\]
as required.

The proof that (\ref{eq:b1}) implies (\ref{eq:b2}) is exactly analogous,
of course using (\ref{eq:fpdh}) in place of (\ref{eq:hlte}). $\qed$

\subsection{\label{sub:LipschitzAndBMO} Compositions of BMO functions with Lipschitz
functions}

We would not be at all surprised if the result (Proposition \ref{pro:ttbf})
presented in this appendix is already known. But somehow we have not
found a reference for it yet. The following completely obvious lemma
and not very difficult proposition are not needed for obtaining the
main results of this paper. But they may be of independent interest.
They are also a kind of motivation for some of the steps for proving
our main results.
\begin{lem}
\label{lem:lipz}Suppose that the function $\varphi:\mathbb{R}\to\mathbb{R}$
satisfies the Lipschitz condition \begin{equation}
\left|\varphi(s)-\varphi(t)\right|\le\left|s-t\right|\mbox{ for all }s,t\in\mathbb{R}\,.\label{eq:lcwc}\end{equation}
 Then, for each function $f\in BMO(\mathbb{R}^{d})$ the composed
function $\varphi\circ f$ is also in $BMO(\mathbb{R}^{d})$ and satisfies
\[
\left\Vert \varphi\circ f\right\Vert _{BMO(\mathbb{R}^{d})}\le\left\Vert f\right\Vert _{BMO(\mathbb{R}^{d})}\,.\]

\end{lem}
\noindent \textit{Proof.} Let $Q$ be an arbitrary cube in $\mathbb{R}^{d}$
and let $c$ be an arbitrary real constant. It follows immediately
from (\ref{eq:lcwc}) that \[
\int_{Q}\left|\varphi\circ f-\varphi(c)\right|d\lambda\le\int_{Q}\left|f-c\right|d\lambda\]
which completes the proof. $\qed$\bigskip{}

Note that in the following proposition we do not even need to require
the function $f$ to be locally integrable. One of its immediate consequences
is that a measurable function $f:\mathbb{R}^{d}\to\mathbb{R}$ is
in $BMO(\mathbb{R}^{d})$ if and only if the $BMO$ seminorms of all
the bounded functions $f_{k}(x)=k\arctan\left(\frac{f(x)}{k}\right)$
are all dominated by a (finite) constant which is independent of $k\in\mathbb{N}$.
The same conclusion follows when the bounded functions $f_{k}$ are
defined instead by \[
f_{k}(x)=\left\{ \begin{array}{ccc}
-k & , & f(x)<-k\\
f(x) & , & -k\le f(x)\le k\\
k & , & f(x)>k\,.\end{array}\right.\]

\begin{prop}
\label{pro:ttbf}Let $\left\{ \varphi_{k}\right\} _{k\in\mathbb{N}}$
be a sequence of non decreasing functions which each satisfy the Lipschitz
condition (\ref{eq:lcwc}). Suppose also that $\lim_{k\to\infty}\varphi_{k}(t)=t$
for all real $t$ and that, for each bounded interval $(a,b)$ there
exists an integer $N$ such that $\varphi_{k}$ is strictly increasing
on $(a,b)$ for each integer $k\ge N$.

For each integer $k\in\mathbb{N}$, let $f_{k}:\mathbb{R}^{d}\to\mathbb{R}$
be defined by $f_{k}(x)=\varphi_{k}(f(x))$. Then $f\in BMO(\mathbb{R}^{d})$
if and only $f_{k}\in BMO$ for all $k\in\mathbb{N}$ and $\limsup_{k\in\mathbb{N}}\left\Vert f_{k}\right\Vert _{BMO}$
is finite. In fact, for each $f\in BMO(\mathbb{R}^{d})$, the sequence
$\left\{ \left\Vert f_{k}\right\Vert _{BMO}\right\} _{k\in\mathbb{N}}$
converges to a finite limit and \[
\left\Vert f\right\Vert _{BMO}=\lim_{k\to\infty}\left\Vert f_{k}\right\Vert _{BMO}\,.\]

\end{prop}
\noindent \textit{Proof.} In view of Lemma\ref{lem:lipz} we have
$\left\Vert f\right\Vert _{BMO}\ge\limsup_{k\to\infty}\left\Vert f_{k}\right\Vert _{BMO}$
which immediately establishes one of the implications and means that
it remains only to prove that \begin{equation}
\left\Vert f\right\Vert _{BMO}\le\liminf_{k\to\infty}\left\Vert f_{k}\right\Vert _{BMO}\label{eq:mtj}\end{equation}
whenever $\limsup_{k\to\infty}\left\Vert f_{k}\right\Vert _{BMO}$
is finite.

So let us indeed suppose that $\limsup_{k\to\infty}\left\Vert f_{k}\right\Vert _{BMO}$
is finite.

Choose an arbitrary cube $Q$ in $\mathbb{R}^{d}$. Suppose that $c$
is a median of $f$ on $Q$. Let us choose an integer $N$ such that,
for each $k\ge N$, the non decreasing function $\varphi_{k}$ is
strictly increasing on the interval $(c-1,c+1)$. Then $\varphi_{k}(f(x))>\varphi_{k}(c)$
if and only if $f(x)>c$, and, similarly, $\varphi_{k}(f(x))<\varphi_{k}(c)$
if and only if $f(x)<c$. It follows that $c_{k}:=\varphi_{k}(c)$
is a median of $f_{k}$ on $Q$. We also have $\lim_{k\to\infty}c_{k}=c$.
For each fixed positive number $M$ we can use the dominated convergence
theorem to obtain that \begin{eqnarray*}
\int_{Q\cap\left\{ x:-M\le f(x)\le M\right\} }\left|f-c\right|d\lambda & = & \lim_{k\to\infty}\int_{Q\cap\left\{ x:-M\le f(x)\le M\right\} }\left|f_{k}-c_{k}\right|d\lambda\\
 & = & \liminf_{k\to\infty}\int_{Q\cap\left\{ x:-M\le f(x)\le M\right\} }\left|f_{k}-c_{k}\right|d\lambda\\
 & \le & \liminf_{k\to\infty}\int_{Q}\left|f_{k}-c_{k}\right|d\lambda\,.\end{eqnarray*}
In view of (\ref{eq:newmib}) (cf.~Lemma \ref{lem:medianoptimal}),
this last expression is dominated by

\noindent \[
\lambda(Q)\liminf_{k\to\infty}\left\Vert f_{k}\right\Vert _{BMO}\,.\]
So, by applying the monotone convergence theorem, we deduce that \[
\int_{Q}\left|f-c\right|d\lambda=\lim_{M\to\infty}\int_{Q\cap\left\{ x:-M\le f(x)\le M\right\} }\left|f-c\right|d\lambda\le\lambda(Q)\liminf_{k\to\infty}\left\Vert f_{k}\right\Vert _{BMO}\]
which shows that $f\in BMO$ and establishes (\ref{eq:mtj}), so completing
the proof. $\qed$

\bigskip{}

\begin{rem}
To obtain the results mentioned in the preamble to the preceding proposition
we of course simply choose $\varphi_{k}$ defined by $\varphi_{k}(t)=k\arctan\left(\frac{t}{k}\right)$
or \[
\varphi_{k}(t)=\left\{ \begin{array}{ccc}
-k & , & t<-k\\
f(x) & , & -k\le t\le k\\
k & , & t>k\end{array}\right.\]
 respectively. 
\end{rem}

\subsection{\label{sub:convexboundary}The boundary of a convex set in $\mathbb{R}^{d}$}

It is surely very very well known that the boundary $\partial K$
of a convex $K$ set in $\mathbb{R}^{d}$ satisfies $\lambda\left(\partial K\right)=0$.
But let us give an explicit proof of this in this preliminary {}``lecture
notes'' version of the paper. By translating $K$ if necessary, we
can assume that the origin $\vec{0}=\left(0,0,...,0\right)$ of $\mathbb{R}^{d}$
is in $K$. We will also temporarily assume that $K$ is bounded.
If $\lambda(K)=0$ there is nothing to prove. This means that we can
assume that $K$ is not contained in any $d-1$ dimensional subspace
of $\mathbb{R}^{d}$. Consequently $K$ contains the convex hull of
$\vec{0}$ and $d$ more points whose {}``span'' is $\mathbb{R}^{d}$.
It follows that $K$ must contain a cube $Q$. Since we are permitted
to translate $K$ yet again if necessary, we may now assume that the
centre of $Q$ is at $\vec{0}$. For each $r\in(0,1)$ and each $x\in K$
the set $rx+(1-r)Q$ is contained in $K$. This means that the convex
set $rK$ is contained in the interior $K^{\circ}$ of $K$. Thus
$\lambda\left(K^{\circ}\right)\ge r^{d}\lambda(K)$ for every $r\in(0,1)$.
Consequently $\lambda(K^{\circ})=\lambda(K)$ which of course implies
that $\lambda\left(\partial K\right)=0$.

Finally, if $K$ is unbounded, we can set $K_{n}=\left\{ x\in K:\left\Vert x\right\Vert \le n\right\} $
and use the fact that $\partial K\subset\bigcup_{n\in\mathbb{N}}\partial K_{n}$.

\subsection{\label{sub:affinecommute}The functional $\mathbf{J}(f,E,s)$ {}``commutes''
with affine transformations.}

Here is the straightforward calculation which proves (\ref{eq:preaffine})
and therefore also (\ref{eq:affinestuff}) for the function $g$ defined
by $g(x)=f(rx+x_{0})$ where $r$ is a non zero real number and $x_{0}$
is a constant point in $\mathbb{R}^{d}$ and $E$ is an admissible
set in the domain of definition of $g$. \begin{eqnarray*}
\lambda\left(\left\{ x\in E:\left|g(x)-c\right|>\alpha\right\} \right) & = & \lambda\left(\left\{ x\in E:\left|f(rx+x_{0})-c\right|>\alpha\right\} \right)\\
 & = & \lambda\left(\left\{ rx\in rE:\left|f(rx+x_{0})-c\right|>\alpha\right\} \right)\\
 & = & \lambda\left(\frac{1}{r}\left\{ y\in rE:\left|f(y+x_{0})-c\right|>\alpha\right\} \right)\\
 & = & r^{-d}\lambda\left(\left\{ y\in rE:\left|f(y+x_{0})-c\right|>\alpha\right\} \right)\\
 & = & r^{-d}\lambda\left(\left\{ y+x_{0}\in rE+x_{0}:\left|f(y+x_{0})-c\right|>\alpha\right\} \right)\\
 & = & r^{-d}\lambda\left(\left\{ x\in rE+x_{0}:\left|f(x)-c\right|>\alpha\right\} -x_{0}\right)\\
 & = & r^{-d}\lambda\left(\left\{ x\in rE+x_{0}:\left|f(x)-c\right|>\alpha\right\} \right)\,.\end{eqnarray*}

This calculation, together with the fact that $\lambda(rE+x_{0})=r^{d}\lambda(E)$,
shows that the condition \[
\lambda\left(\left\{ x\in E:\left|g(x)-c\right|>\alpha\right\} \right)<s\lambda(E)\]
is equivalent to \[
\lambda\left(\left\{ x\in rE+x_{0}:\left|f(x)-c\right|>\alpha\right\} \right)<s\lambda(rE+x_{0})\]
which is exactly what we need to show that $\mathbf{J}(g,E,s)=\mathbf{J}(f,rE+x_{0},s)$
for each admissible set $E$.

\subsection{\label{sub:lrconj}Left continuity and right discontinuity of the
function $s\mapsto\mathbf{J}(f,Q,s)$.}
\begin{lem}
Suppose that $Q$ and $f$ are as in part (ii) of Proposition \ref{pro:prop}.
Then the function $s\mapsto\mathbf{J}(f,Q,s)$ is non increasing and
left continuous. 
\end{lem}
\noindent \textit{Proof.} Since we know that $s\mapsto\mathbf{J}(f,Q,s)$
is non increasing, in order to show that this function is also left
continuous, we only need to prove that \[
\lim_{n\to\infty}\mathbf{J}\left(f,Q,s-1/n\right)\le\mathbf{J}(f,Q,s)\]
for each fixed $s\in(0,1)$. As in the proof of Proposition \ref{pro:prop},
let $g:Q\to\mathbb{R}$ denote the function which is the restriction
of $f$ to $Q$. Fix an arbitrary positive number $\varepsilon$.
In view of (\ref{eq:ffnn}) there exists $u\in(0,s\lambda(Q))$ such
that \[
g^{*}(u)-g^{*}(u+(1-s)\lambda(Q))\le\mathbf{J}(f,Q,s)+\varepsilon/2\,.\]
Of course we have $u<(s-1/n)\lambda(Q)$ for all sufficiently large
$n$. Furthermore, because of the right continuity of $g^{*}$, the
inequality \[
g^{*}\left(u\right)-g^{*}(u+(1-s+1/n)\lambda(Q))\le\mathbf{J}(f,Q,s)+\varepsilon\]
also holds for all sufficiently large $n$. These last two inequalities
imply, again in view of (\ref{eq:ffnn}), that $\mathbf{J}(f,Q,s-1/n)\le\mathbf{J}(f,Q,s)+\varepsilon$
for all sufficiently large $n$, and this suffices to complete the
proof. $\qed$
\begin{example}
\label{rem:rexample}Suppose that $d=1$ and $Q=(0,1)$ and $a\in(0,1)$.
Let $N$ be the smallest integer for which $Na\ge1$ and let $f$
be the restriction to $(0,1)$ of the function \[
N\chi_{(0,a)}+(N-1)\chi_{[a,2a)}+(N-2)\chi_{[2a,3a)}...+\chi_{[(N-1)a,Na)}\,.\]
Then $f$ is non increasing and right continuous and $f^{*}=g^{*}=f$
on $(0,1)$. Consequently, by (\ref{eq:ffnn}), $\mathbf{J}(f,(0,1),s)\ge1/2$
for $s\in(0,1-a]$ and $\mathbf{J}(f,(0,1),s)=0$ for $s\in(1-a,1)$.
This shows that, in general, $s\mapsto\mathbf{J}(f,Q,s)$ is not right
continuous.
\end{example}

\subsection{\label{sub:ProveSimpleProps}The proof of Lemma \ref{lem:SimplePropsCubes}.}

For part (i), let $y$ be an arbitrary point of $Q(x_{*},r_{*})$.
For each $n\in\mathbb{N}$ we set $y_{n}=x_{n}+\frac{r_{n}}{r_{*}}(y-x_{*})$\,.
Since $\left\Vert y-x_{*}\right\Vert _{\ell_{d}^{\infty}}\le r_{*}$
it follows that $y_{n}\in Q(x_{n},r_{n})$. Therefore $y_{n}\in Q$.
Since $Q$ is closed it follows that $y=\lim_{n\to\infty}y_{n}$ is
also in $Q$ and this shows that $Q(x_{*},r_{*})\subset Q$.

For part (ii) we observe that \begin{equation}
\begin{array}{cl}
 & \left|\lambda(A\cap Q(x_{n},r_{n}))-\lambda(A\cap Q(x_{*},r_{*}))\right|\le\left|\int_{Q(x_{n},r_{n})}\chi_{A}d\lambda-\int_{Q(x_{\ast},r_{\ast})}\chi_{A}d\lambda\right|_{\phantom{Q}}\\
= & \left|\int_{Q(x_{n},r_{n})\setminus Q(x_{\ast},r_{\ast})}\chi_{A}d\lambda-\int_{Q(x_{\ast},r_{\ast})\setminus Q(x_{n},r_{n})}\chi_{A}d\lambda\right|_{\phantom{Q}}^{\phantom{T}}\\
\le & \lambda\left(Q(x_{n},r_{n})\setminus Q(x_{\ast},r_{\ast})\right)+\lambda\left(Q(x_{\ast},r_{\ast})\setminus Q(x_{n},r_{n})\right)^{\phantom{\begin{array}{cc}
\end{array}}}\,.\end{array}\label{eq:zerpluf}\end{equation}
Let $R_{n}=\left\Vert x_{n}-x_{*}\right\Vert _{\ell_{d}^{\infty}}+\max\left\{ r_{*},r_{n}\right\} $.
It is clear (from (\ref{eq:defqxr})) that the cube $Q(x_{*},R_{n})$
contains both of the cubes $Q(x_{*},r_{*})$ and $Q(x_{n},r_{n})$.
Therefore the expression in the last line of (\ref{eq:zerpluf}) is
dominated by \begin{eqnarray*}
 &  & \lambda\left(Q(x_{\ast},R_{n})\setminus Q(x_{\ast},r_{\ast})\right)+\lambda\left(Q(x_{\ast},R_{n})\setminus Q(x_{n},r_{n})\right)\\
 & = & R_{n}^{d}-r_{*}^{d}+R_{n}^{d}-r_{n}^{d}\,.\end{eqnarray*}
This last expression tends to $0$ as $n$ tends to $\infty$. This
establishes (\ref{eq:ttz}) and so completes the proof of the lemma.
$\qed$

\end{document}